%% file: manuscript_clean.tex
\documentclass[preprint]{elsarticle}
\usepackage[utf8]{inputenc}

\usepackage{todonotes}
\usepackage{caption}

\usepackage{hyperref}
\usepackage{amsfonts}
\usepackage{geometry}
\usepackage[ruled]{algorithm2e}
\usepackage{mathtools}

\usepackage{amsmath}

\usepackage{amssymb}   
\usepackage{amsthm}

\usepackage{pict2e}
\usepackage{keyval}
\usepackage{calc}
\usepackage{fp}
\usepackage{diagbox}
\usepackage{booktabs}
\usepackage{tabularx,colortbl, makecell, caption, multirow}
\usepackage{boldline}

\def\Tb{\mathcal{T}_B}
\def\Tbf{\mathcal{T}_{B}'}
\def\Tbl{\mathcal{T}_{B}''}
\def\Tl{\mathcal{T}_L}

\def\H{\mathcal{H}}
\def\S{\mathcal{S}}

\def\I{\mathcal{I}}
\newcommand{\R}{\mathbb{R}}

\def\sgn{\mbox{\rm sgn}}

\usepackage{rotating}
\usepackage{graphicx}
\usepackage{bm}
\usepackage{xcolor}
\usepackage{multirow}
\usepackage{nicematrix}

\usepackage{graphicx}
\graphicspath{{images/}}
\usepackage{float}
\usepackage{longtable}
\usepackage{adjustbox}

\newcommand{\ds}{\displaystyle}

\usepackage{colortbl} 

\graphicspath{{./figures/}}

\usepackage[titletoc,toc,title]{appendix}
\usepackage{comment}
\usepackage{graphics}
\usepackage{graphicx}
\usepackage{adjustbox}
\usepackage{tikz}
\usepackage{caption}
\usepackage{psfrag}
\usepackage{soul}
\usepackage{mwe}
\usepackage{pgfplots}
\pgfplotsset{compat=1.17} 

\graphicspath{{./figures/}}

\usepackage{comment}
\usepackage{graphics}
\usepackage{graphicx}
\usepackage{adjustbox}
\usepackage{tikz}
\usepgfplotslibrary{dateplot}
\usepgfplotslibrary{groupplots}
\usepackage{caption}
\usepackage{subcaption}
\usepackage{psfrag}
\usepackage{soul}
\usepackage{mwe}
\usepgfplotslibrary{fillbetween}

\usetikzlibrary{er,positioning, calc, shapes, arrows, chains, fit}

\definecolor{mygreen}{RGB}{80,176,50}

\usepackage{xcolor}

\makeatletter

\newcommand*{\@rowstyle}{}

\newcommand*{\rowstyle}[1]{
  \gdef\@rowstyle{#1}%
  \@rowstyle\ignorespaces%
}

\newcolumntype{=}{
  >{\gdef\@rowstyle{}}%
}

\newcolumntype{+}{
  >{\@rowstyle}%
}

\makeatother

\definecolor{red}{rgb}{0.0, 0.0, 0.0}
\definecolor{dred}{rgb}{0.55, 0.0, 0.0}
\definecolor{dblue}{rgb}{0.0, 0.2, 0.6}
\definecolor{dgreen}{RGB}{80,176,50}

\usepackage{colortbl}

\definecolor{mygreen}{RGB}{80,176,50}
\definecolor{myred}{RGB}{222,36,16}
\definecolor{myind}{RGB}{93,174,255}
\definecolor{mypink}{RGB}{219,129,255}
\definecolor{azzdark}{RGB}{42,72,180}
\definecolor{azzlight}{RGB}{53,146,203}
\definecolor{verdone}{RGB}{0,128,0}
\definecolor{indaco}{RGB}{68,137,206}
\definecolor{celestino}{RGB}{222,247,254}
\definecolor{verdeacceso}{RGB}{71,185,109}

\definecolor{pred}{RGB}{182,7,40}
\definecolor{pblue}{RGB}{58,76,189}

\newcommand{\red}[1]{\textcolor{red}{#1}}
\newcommand{\dblue}[1]{\textcolor{dblue}{#1}}

\newcommand{\mygreen}[1]{\textcolor{mygreen}{#1}}

\usetikzlibrary{er,positioning, calc, shapes, arrows, chains, fit, decorations.pathreplacing}

\definecolor{orange}{rgb}{0.8, 0.47, 0.13}
\definecolor{green}{rgb}{0.0, 0.62, 0.42}

\newcommand{\fnote}{\textcolor{black}}
\newcommand{\lnote}{\textcolor{black}}

\newcommand{\fnotee}{\textcolor{black}}
\newcommand{\lnotee}{\textcolor{black}}
\newcommand{\mnotee}{\textcolor{black}}

\usepackage{fancyhdr}

\begin{document}

\begin{frontmatter}

\title{
Margin Optimal Classification Trees}

\author[1]{Federico D'Onofrio}
\author[2]{Giorgio Grani}
\author[1]{Marta Monaci}
\author[1]{Laura Palagi}

\address[1]{Department of Computer, Control and Management Engineering Antonio Ruberti (DIAG), \\Sapienza University of Rome, Rome, Italy}

\address[2]{Department of Statistical Sciences, Sapienza University of Rome, Rome, Italy}

\begin{abstract}

In recent years\lnote{,} there has been growing attention to interpretable machine learning models which can give explanatory insights on  
their behaviour.
Thanks to their interpretability,
decision trees have been intensively studied for classification tasks and\mnotee{,} due to the remarkable advances in mixed integer programming (MIP), various approaches have been proposed to formulate the problem of training an Optimal Classification Tree (OCT) as a MIP model. 
We present a novel mixed integer quadratic formulation for the OCT problem, which exploits the generalization capabilities of Support Vector Machines 
for binary classification. 
\lnote{Our model, denoted as 
 Margin Optimal Classification Tree (MARGOT), encompasses} maximum margin multivariate hyperplanes nested in a binary tree structure. To enhance the interpretability of our approach, we analyse two alternative versions of MARGOT, which include feature selection constraints inducing \fnotee{sparsity of the hyperplanes' coefficients}. 
 First, MARGOT has been tested on non-linearly separable synthetic datasets in a 2-dimensional feature space to provide a graphical representation of the maximum margin approach. Finally, the proposed models have been tested on benchmark datasets from the UCI repository.
 The MARGOT formulation turns out to be easier to solve than other OCT approaches, and the generated tree better generalizes on new observations.
 \lnote{The two interpretable versions effectively select the most relevant features, maintaining good prediction quality.} 
 \end{abstract}


\end{frontmatter}

\parindent 0cm


\pagestyle{fancy}
\lhead{}
\renewcommand{\headrulewidth}{0pt}
\rhead{\textcolor[gray]{0.3}{Published at: \href{https://doi.org/10.1016/j.cor.2023.106441}{https://doi.org/10.1016/j.cor.2023.106441}}}
\thispagestyle{fancy}  

\textcopyright{ 2023. Licensed under the Creative Commons  \href{https://creativecommons.org/licenses/by-nc-nd/4.0/}{CC-BY-NC-ND 4.0}.} 

\section{Introduction}

\subsection{Related Work}
In recent years\lnote{,} there has been growing interest in interpretable Machine Learning (ML) models \cite{rudin2022interpretable}. Decision trees are among the most popular Supervised ML tools used for classification tasks. They are famous for being easy to manage, \mnotee{having low} computational requirements, and the final model is easily understandable from a human perspective as opposed to other ML methods that are seen as black boxes.
Given a set of points and class labels, a classification tree method builds up a binary tree structure of a maximum predefined depth. Trees are composed of branch and leaf nodes, and the branch nodes apply a sequence of dichotomic rules, called \textit{splitting rules}, to partition the training samples into disjoint subsets. {Splitting rules route samples to the left or right child node, and they are usually defined by hyperplanes. In a univariate tree, these hyperplanes are orthogonal\lnote{,} involving one single feature, while in a multivariate tree, they can be oblique\lnote{,} involving more than one feature. A value for the predicted class label is assigned to each leaf node according to some simple rule, for instance, the most common label rule.}
The key advantage of tree methods lies in their interpretability. The process behind a decision tree is transparent, and the sequential tree structure mimics \mnotee{the} human decision-making process. These properties can be crucial factors in many applications, ranging from business and criminal justice to healthcare and bioinformatics. \lnote{Indeed, in these domains, it is of great interest to use explainable approaches to help humans understand the model's decisions and identify a subset of the most prominent features that influence the classification outcome.}
To this aim, it is preferable to build and manage shallow trees with small depth; indeed, if allowed to grow large, decision trees lose their interpretability aspect.

It is well-known that constructing a binary decision tree in an optimal way is an \emph{NP}-complete problem \cite{hyafil1976constructing}. 
For this reason, traditional approaches for finding decision trees rely on heuristics. In general, they are based on a top-down greedy strategy for growing the tree by generating splits at each node, and, once the tree is built, a bottom-up pruning procedure is applied to handle the complexity of the tree, i.e. the number of splits. 
Breiman et al.~\cite{Breiman1984CART} proposed a heuristic algorithm known as CART (Classification and Regression Trees) for learning univariate decision trees.
Starting from the root node, each hyperplane split is generated by minimizing a local impurity function, e.g. the Gini impurity for classification tasks.

Other univariate approaches employing different impurity functions were later proposed by Quinlan \cite{quinlan1986induction, quinlan1993c45programs} in his ID3 and C4.5 algorithms. 
\lnote{These heuristic procedures produce a solution in fast computational time but may also generate tree models with poor generalization performances. }
In order to overcome these drawbacks,
tree ensemble methods, such as Random Forests \cite{Breiman2001RandomForests}, TreeBoost \cite{FriedmanBoostedTrees} and XGBoost \cite{XGBoostChen}, have been proposed. These approaches combine decision trees using some kind of randomness; however\mnotee{,} using multiple trees leads to a lack of interpretability of the final model.

Another way to improve prediction quality is to use multivariate decision trees\mnotee{,} which employ oblique hyperplane splits. These methods involve more features per split, thus producing smaller trees but at the expense of the computational cost. Several approaches for inducing multivariate trees have been proposed (see \cite{MurthyObliqueDT, BrodleyMultivariateDT, Vercellis2003Multivariate, Wickramarachchi2016HHCART}). For instance, OC1 \cite{MurthyObliqueDT} is a greedy algorithm that searches for the best hyperplane at each node by applying a randomized perturbation strategy\lnote{. In contrast}, \cite{Vercellis2003Multivariate} presented a heuristic procedure that, at each step, solves a variant of the Support Vector Machine problem, where the empirical error is discretized by counting the number of misclassified samples.

Recently, several papers have been devoted to global exact optimization approaches to find an Optimal Classification Tree (OCT) using mathematical programming tools and, in particular\mnotee{,} Mixed Integer Programming (MIP) models (see the recent surveys \cite{GAMBELLA2021807, Morales2021MathematicalOptimizationOCRT} and references therein). Indeed, the significant improvement in the last thirty years of both algorithms for integer optimization and computer hardware has led to an incredible increase in the computational power of mixed integer solvers, as shown in \cite{Bixby2012MIP}. Thus, MIP approaches became viable in the definition of ML methods, being \cite{Bertsimas2017OptimalClassification} the seminal paper inaugurating a new era in \lnote{using} mixed integer based optimization to learn OCTs.
Such approaches find the decision tree in its entirety through the resolution of a single optimization model, defining each branching rule with full knowledge of all the remaining ones.
In \cite{Bertsimas2017OptimalClassification}, Bertsimas and Dunn proposed two Mixed Integer Linear Programming (MILP) models to build optimal trees with a given maximum depth based on univariate and multivariate splits.
Along these lines, Günlük et al.~\cite{Gunluk2021OptimalDT} proposed a MIP formulation for binary classification tasks by exploiting the structure of categorical features and \mnotee{modelling} combinatorial decisions. Further, in order to circumvent the problem of the curse of dimensionality related to the MIP approaches, Verwer and Zhang \cite{Verwer2019Learning} presented BinOCT, a binary linear programming model,
where the size is independent of the training set dimension.
In \cite{Aghaei2021StrongOCT}, Aghaei et al.~proposed a flow-based MILP model for binary features with a stronger linear relaxation and, by exploiting the decomposable and combinatorial structure of the model, derived a Benders' decomposition method to deal with larger instances.
\red{Boutilier et al. \cite{Boutilier2022ShatterinInequalities} presented a new formulation for learning multivariate optimal trees. Moreover, they introduced a new class of valid inequalities and leveraged them within a Benders-like decomposition to improve the optimization process.}
In addition to expressing the combinatorial nature of the decisions involved in the process, the mixed integer framework is suitable to handle global objectives and constraints to embed desirable properties such as fairness, sparsity, cost-sensitivity, robustness, as it has been addressed in \cite{Morales2021MathematicalOptimizationOCRT}, \cite{Verwer2017LearningDT}, \cite{Aghaei2019LearningOCfair}, \cite{Aghaei2021StrongOCT} and \cite{BlancoRobust}.

Alongside integer optimization, continuous optimization paradigms have also been investigated \lnote{in the context of optimal trees.} In \cite{Morales2021OptimalRandomized}, Blanquero et al.~proposed a nonlinear 
programming model for learning an optimal "randomized" classification tree with oblique splits, where at each node\lnote{,} a random decision is made according to a soft rule, induced by a continuous cumulative density function. Later, in \cite{Morales2020Sparsity}, they addressed global and local sparsity in the randomized optimal tree model \red{(S-ORCT)} by means of regularization terms based on polyhedral norms \red{($\ell_1$-norm and $\ell_{\infty}$-norm)}. In their randomized framework, a sample is not assigned to a class in a deterministic way but only with a given probability.
\red{Following this research line, Amaldi et al. \cite{Amaldi2023OnMultivariateRandomized} investigated additional versions of the S-ORCT model based on concave approximations of the $\ell_0$-norm and proposed a general proximal point decomposition scheme to tackle larger datasets}. 

Following a different viewpoint, approaches using a Support Vector Machine (SVM) (see \cite{VapnikSVM1995}, \cite{WangSVM2005}, \cite{Piccialli2018NonlinearOA}) for each split in the tree have been investigated. First, Bennett et al.~\cite{Bennett1998ASupport} 
provided a primal continuous formulation with a non-convex objective function and a dual convex quadratic model to train optimal trees where each decision rule is based on a modified SVM problem. 
\mnotee{Their model can involve kernel functions to construct nonlinear splitting rules}.
\mnotee{The resulting problems are computationally hard to solve, and a tabu search algorithm is used to approximately find solutions.}
Recently, margin-based splits of the SVM type have been proposed by Blanco et al. in \cite{BlancoRobust}. The authors introduced a Mixed Integer \mnotee{Nonlinear} formulation for the OCT problem to solve binary classification tasks.
The aim \mnotee{was} to build a robust tree classifier, where during the training phase\lnote{,} some of the labels of the dataset are allowed to be changed \mnotee{in order} to detect the label noise. Observations are relabelled based on misclassification errors, as   
described in \cite{BlancoRelabelingSVM}.
The method aims to seek a trade-off between four objectives, the first being the maximization of the minimum margin among all the margins of the hyperplane splits in the tree. 
\mnotee{In addition, it minimizes the misclassification cost at the branch nodes, the number of relabelled observations and the complexity of the tree.}
The model is formulated as a Mixed Integer Second Order Cone Optimization problem.
In \cite{blanco2023multiclass}, an extension of the model in \cite{BlancoRobust} for handling multiclass instances has been proposed. 

\subsection{Our contribution}

Our approach falls \mnotee{within} the basic framework of using Support Vector Machines to define optimal classification trees with multivariate splits for binary classification tasks. In particular, our model employs maximum margin hyperplanes obtained using a linear soft SVM paradigm in a nested binary tree structure. The maximum depth of the tree is fixed, as usual in OCT approaches. 

The main contribution of this paper is a novel Mixed Integer Quadratic Programming (MIQP) formulation, denoted as Margin Optimal Classification Tree (MARGOT), for learning classification trees. Our formulation differs from others in the literature as we exploit the statistical learning properties of the $\ell_2$-regularized soft SVM formulation. In particular, the SVM quadratic convex function is retained\lnote{,} and the collective measure of performance of the OCT is obtained as the sum of the objective functions of each SVM-based problem over all the branch nodes of the tree. Differently, in \cite{BlancoRobust}, only the minimum among all the hyperplane margins is maximized.
Further, exploiting both the SVMs properties and the binary classification setting \mnotee{allows} us to drastically reduce the overall number of binary variables needed in our MIQP model. 
 Specifically, we need to introduce as many binary variables as half the number of leaf nodes, which is much less than other OCT MIP approaches. 
We show, both on synthetic datasets in a 2-dimensional feature space and on datasets
selected from the UCI Machine Learning repository \cite{UCI2019}, that MARGOT formulation requires less computational effort than other state-of-the-art MIP models for OCT, and it can be solved to certified optimality on nearly all the considered problems with a limited computational time using off-the-shelf solvers.
As a consequence of the maximum margin approach, our model produces OCTs with a higher out-of-sample accuracy.

As a second contribution, we aim to reduce the number of features used in each split to enhance the interpretability of the model itself.
Indeed, for tabular data, sparsity is a core component of interpretability 
\cite{rudin2022interpretable}, and having fewer features selected at each branching node allows the end user to identify the key factors affecting the classification of the samples.
  
 Actually, due to the intrinsic statistical properties of SVMs, such a model tends to use a large number of features to define each split of the classification tree.  For this reason,  we propose two embedded models that simultaneously train the OCT and perform feature selection.
Embedded models for feature selection in SVMs have been studied in several papers (see e.g. \cite{JIMENEZCORDERO202124,Carrizosa2011,MaldonadoSVMFS2014,LabbeMixedInteger,LEE20221055}). To control the sparsity of the oblique splits, we draw inspiration from the model proposed in \cite{LabbeMixedInteger}\mnotee{,} and we use additional binary variables and a budget on the number of features used. 
We consider two different \mnotee{modellings} of the budget on the number of features: hard constraints and soft penalization, respectively implemented in HFS-MARGOT and SFS-MARGOT. Numerical \mnotee{results} on the UCI datasets are reported, showing that the hard version seems to be easier to solve to certified optimality in a reasonable CPU time, resulting in a more sparse solution too. For all the formulations, we propose a simple greedy heuristic to obtain a first incumbent which exploits the SVM-based tree structure. \medskip

The rest of the paper is organized as follows. In Section \ref{sec: preliminaries}, a brief introduction about Multivariate Optimal Classification Trees, proposed in \cite{Bertsimas2017OptimalClassification}, and Support Vector Machines is provided.
In Section \ref{sec: margot model}, we introduce our approach and its formulation, denoted as MARGOT. In Section \ref{sec: margot with fs}, we present the two interpretable versions of the model with hard and soft feature selection techniques to address the sparsity of the \mnotee{hyperplanes' weights}. In Section \ref{sec: heuristic_ws}, we provide a heuristic to generate starting feasible solutions for the \mnotee{analysed} MIP problems. Then, in Section \ref{sec: results}, we first evaluate MARGOT on 2-dimensional synthetic datasets, and we report a graphical representation of the generated trees. Finally, computational experiments on benchmark datasets from the UCI repository are presented for all the proposed models.




\section{Preliminaries}\label{sec: preliminaries}

\subsection{Multivariate Optimal Classification Trees}
\label{sec:prelim_OCT}

In this section, we introduce in more detail multivariate optimal classification trees.
Given a dataset $\{(x^i,y^i) \in \mathbb{R}^n\times\{1,\dots,K\},\ i\in\I\}$, and a maximum depth $D$, an optimal classification tree is made up by $2^{(D+1)}-1$ nodes, divided in:
\begin{itemize}
    \item \textit{Branch nodes},
    $\Tb = \{0,\dots, 2^D-2\}$: a branch node applies a splitting rule on the feature space defined by a separating hyperplane $\H_t(x):= \{x: h_t(x)=0\}$, where $h_t(x) = w^T_t x + b_t$ is the hyperplane function and $w_t \in  \mathbb{R}^n$ and $ b_t\in\mathbb{R}$. If $h_t(x^i) \geq 0 $, sample $i$ will follow the right branch of node $t$, otherwise it will follow the left one;

    \item \textit{Leaf nodes}, 
    $\Tl = \{2^D - 1  , \dots , 2^{D+1} - 2\}$: leaf nodes act as collectors\mnotee{,} and the samples which end up in the same leaf are classified with the same class label.
    
    
\end{itemize}

The training phase aims at building a classification tree by finding coefficients $w_t$ and the intercept $b_t$ for each $t\in\Tb$ and by assigning class labels to the leaf nodes. According to the hierarchical tree structure, the feature space will be partitioned into disjoint regions, \mnotee{each} corresponding to a leaf node of the tree. 
The obtained tree is then used to classify out-of-sample data: every new sample will follow a unique path within the tree based on the splitting rules, ending up in a leaf node that will predict its class label. 

In \cite{Bertsimas2017OptimalClassification}, Bertsimas and Dunn proposed a MILP optimization model for training multivariate OCTs, denoted as OCT-H, where, in the objective function, the misclassification error together with the number of features used at each split is minimized. In the optimization model, each sample is forced to end up in a single leaf, a class label for each leaf node is chosen according to the most common label rule and the classification error is computed according to the assignment of each sample to a leaf. 
\textit{Routing constraints} enforce each sample to follow a unique path, while other constraints control the complexity of the tree by imposing pruning conditions and a minimum number of points accepted by each \mnotee{non-empty} leaf.  

\subsection{A brief overview on Support Vector Machines}

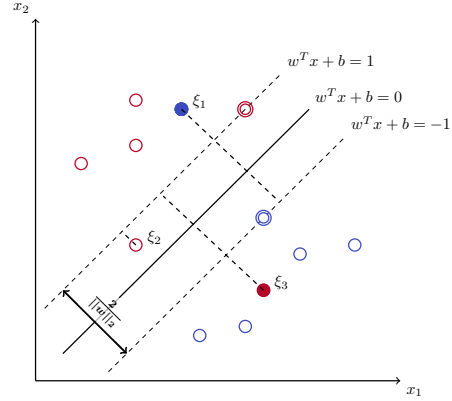
\begin{figure}[ht]
\centering
\scalebox{.6}{\input{SVMplot}}
\caption{Example showing how errors are upper-bounded in the SVM approach by variables $\xi_i$: values  $0<\xi_i\le 1$ correspond to samples that lie within the margin but are correctly classified;  \fnotee{values} $\xi_i>1$ \mnotee{correspond to} samples wrongly classified.
}
\label{fig:SVM}
\end{figure}

Given a binary classification instance $\{(x^i,y^i) \in \mathbb{R}^n\times\{-1,1\},\ i\in\I\}$, the linear \mnotee{soft-margin} SVM classification problem defines a linear classifier as the function 
$f:\R^n\to \{-1,1\}$, 
$$f(x) = \sgn({w^*}^T x + b^*),$$
using the structural risk minimization
principle (\cite{VapnikSVM1995, vapnik1999nature}). In the traditional $\ell_2$-regularized $\ell_1$-loss linear SVM, coefficients $(w^*,b^*)\in\R^n\times\R$ identify a separating hyperplane \lnote{that} maximizes its margin, i.e.
the distance \mnotee{between two parallel
hyperplanes, each supporting samples belonging to one of the two classes.}
The tuple is found by solving the following convex quadratic problem:


\begin{align}
    \text{(SVM}) \quad
    \min \limits_{{w,b,\xi}} \quad & \ds\frac 1 2  \|w\|^2_2 +C \sum_{ i\in \I} \xi_{i} \label{cons:SVM_obj}\\
    \text{s.t.} \quad 
    & y^i(w^Tx^i+b)\ge  1 - \xi_{i} && \forall i \in \I \label{cons: SVM margin cons}\\
    & \xi_i \geq 0 &&\forall i \in \I, \nonumber
\end{align}

where $C$ is \mnotee{a} hyperparameter that balances the two objectives:
 the maximization of the margin $2\|w\|_2^{-1}$  
 and the minimization of the misclassification cost. 
Variables $\xi_i$ allow for violation of the \textit{margin constraints} \eqref{cons: SVM margin cons} and a sample $i$ is misclassified when $\xi_i> 1$, while values $0 < \xi_i \leq 1$ correspond to correctly classified samples inside the margin. 
The further a misclassified data point $x^i$ \mnotee{is} from a feasible hyperplane, the greater the value of variable $\xi_i$ \mnotee{will be} (see Fig. \ref{fig:SVM}).
Thus, $\sum_{i\in\I}\xi_i$ is an upper bound on the number of samples misclassified by the hyperplane. 
 

Although the objective function of the SVM problem is loosely convex, in \cite{burges1999uniqueness}, necessary and sufficient conditions are given for the support vector solution to be unique. 
In particular, with reference to \eqref{cons:SVM_obj} where $C$ is the same for all $i$, a necessary condition for the solution
to be non-unique is that the negative and positive polarity support vectors are equal in number.
Further, it \fnotee{has been proven} that\mnotee{,} even when the solution is not unique,  all solutions share the same $w$. Thus, among the infinite separating hyperplanes, SVM selects the unique one that 
maximizes the margin.

Minimizing the $\ell_2$-norm of the vector $w$ 
has little effect on \fnotee{its sparsity}, namely in reducing the number of components different from zero. 
In this regard, in the literature, many papers adopted the SVM version where the $\ell_2$-regularization term is replaced by the $\ell_1$ one (see e.g. \cite{bradley2000massive,mangasarian2006exact,fan2008liblinear}), because $\ell_1$-norm acts on the sparsity of the vector. Some references \mnotee{also suggest} the combined use of the two terms \cite{10.2307/24307560,hajewski2018smoothed}.
In the \fnotee{models proposed in this paper,} we consider the $\ell_2$-regularized $\ell_1$-loss linear SVM \mnotee{and, when addressed,} the sparsity of the oblique splits is \mnotee{modelled} by constraints involving additional binary variables, following the idea of \cite{LabbeMixedInteger}.

\section{The Margin Optimal Classification Tree}
\label{sec: margot model}

In this section, we propose a novel MIQP model for constructing optimal classification trees which encompass\fnotee{es} multivariate hyperplanes. The aim is to exploit the generalization capabilities of the soft SVM approach using maximum margin hyperplanes,
thus the name 
Margin Optimal Classification Tree (MARGOT). For the sake of interpretability, we also analyse two alternative versions of MARGOT which \mnotee{reduce} 
the number of features used at each split. These additional models are addressed in Section \ref{sec: margot with fs}.

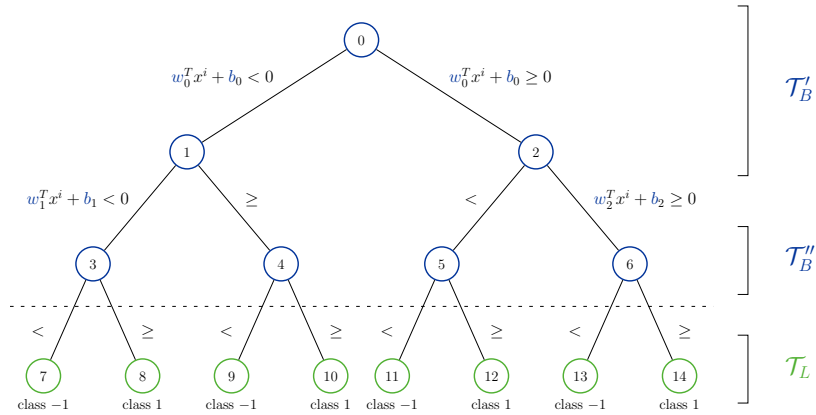
\begin{figure}[ht]
\centering
\hspace{1.2cm}
\scalebox{.45}{{\input{margottona_simple_def}}}
\caption{Representation of a tree model with depth $D=3$ and of the sets $\Tbf,\Tbl, \Tl$.}
\label{fig: margottona}
\end{figure}

In order to formally present the MARGOT formulation, besides  
 the sets of branch and leaf nodes,  we use the following additional notation (see Figure \ref{fig: margottona}). The set of branch nodes $\Tb$ is partitioned \mnotee{into}:
    \begin{itemize}
        \item $\Tbl$, the set of nodes in the last branching level; 
         \item $\Tbf$, the set $\Tb \setminus{\Tbl}$.
    \end{itemize}
    
We also define the following sets:
\begin{itemize}
    \item $\S(t)$\fnotee{,}  the set of nodes of the subtree rooted at node $t\in \Tb$;
    \item $\S^{\prime\prime}(t)\coloneqq\S(t) \cap \Tbl$\fnotee{,}  the subset of nodes of $\S(t)$ belonging to the last branching level $\Tbl$;
    \item $\S^{\prime\prime}_L(t)$ and $\S^{\prime\prime}_R(t)$\fnotee{,} the set of nodes in $\Tbl$ under the left and right branch of node $t\in \Tbf$  such that:
    \begin{itemize}
        \item[-] $\S^{\prime\prime}(t) = \S^{\prime\prime}_L(t) \cup \S^{\prime\prime}_R(t)$
        \item[-] $\S^{\prime\prime}_L(t) \cap \S^{\prime\prime}_R(t) = \emptyset$.
    \end{itemize}
\end{itemize}


\fnotee{The formulation needs to model the fact that} the hyperplane at each node $t$ needs to be trained on just a subset of samples. 
To this aim, let $\I_t\subseteq \I$ be the index set of samples routed to $t\in\Tb$. The definition of the hyperplane at each node $t$ is obtained by means of 
an SVM-type problem. This means that, for each node $t \in \Tb$, we will have standard variables $(w_t,b_t, \xi_t)\in \R^{n}\times \R \times \R^{|\I_t|}$
which must satisfy the soft SVM margin constraints:
\begin{align}
& y^i ( \ds w_t^T x^i + b_t) \ge 1 - \xi_{t,i} &&\qquad \forall i\in\I_t \label{cons: margin cons I_t} \\
& \xi_{t,i} \ge 0 &&\qquad \forall i\in\I_t  \nonumber.
\end{align}

Samples $(x^i, y^i)$, $i\in \I_t$, can be split among the right or left child node of $t$ depending on the rules:
 
\begin{equation}\label{eq:routingrule}
   w_t^T x^i + b_t \geq 0 \;  \text{ if }\; i \in \I_{R(t)}  \; \text{ or } \; w_t^T x^i + b_t < 0  \; \text{ if } \;i \in \I_{L(t)},
\end{equation}

where sets $\I_{R(t)}$ and  $\I_{L(t)}$ are the index sets of samples assigned to the right and left child
nodes of $t$, respectively, thus $\I_{R(t)}\cup \I_{L(t)}=\I_t$ and $ \I_{R(t)} \cap \I_{L(t)}=\emptyset$. A set of routing constraints is therefore needed, 
for each sample $i\in \I_t$, in order to impose the correct sign of the hyperplane function in $x^i$, $h_t(x^i) = w_t^T x^i + b_t$.

The objective function of the single SVM-type problem at node $t$ optimizes the \mnotee{trade-off} between the maximization of the hyperplane margin and the minimization of the upper bound on the misclassification cost given by the sum of the slack variables $\xi_{t,i}$ for all $i\in\I_t$, weighted by a positive coefficient $C_t$: 
$$   \ds\frac 1 2  \|w_t\|^2_2 +  C_t \sum_{i\in\I_t} \xi_{t,i}.  
$$

The aim is to train all the hyperplanes with a single optimization model. Thus, in the objective function, we sum the previous terms over all branch nodes $t \in \Tb$ and all samples $i\in\I_t$:
\begin{equation}
\min \sum_{t \in \Tb } \Big ( \ds\frac 1 2 \|w_t\|^2_2 + C_t \sum_{ i\in\I_t} \xi_{t,i} \Big ).  \nonumber
\end{equation}

However, the route of the samples in the tree, and consequently the definition of subsets $\I_t$ with $t\in\Tb$, is not preassigned, but it is defined by the optimization procedure. 
Hence, we need to define variables $\xi_{t,i}$ and all margin and routing constraints for every sample in $\I$, \lnotee{considering} that constraints at node $t$ must activate only on the subset of samples $\I_t$. 
In order to model the activation/deactivation of these constraints, we need to introduce binary variables which determine the unique path of each sample in the tree. In \mnotee{state-of-the-art} OCT models, such variables model the assignment of each sample in $\I$ to either leaf nodes (resulting in $|\Tl|\cdot|\I|=2^D |\I|$ variables), as in \cite{Bertsimas2017OptimalClassification},  
 or to branch nodes (resulting in  $\ds|\Tb|\cdot|\I|= (2^D-1) |\I|$ \mnotee{variables)}, as in \cite{BlancoRobust}. Routing constraints are defined using these variables, often leading to large MIP models. We aim to reduce as much as possible both the number of binary variables and the constraints used in the model to obtain a more tractable problem. 
For each sample in $\I$ we can introduce such binary variables only for the branch nodes in $\Tbl$, resulting in $|\Tbl|\cdot|\I| = 2^{D-1}|\I|$ variables 
which are half the value $|\Tl|\cdot|\I|$ and less than  $|\Tb|\cdot|\I|$. \mnotee{Indeed}, following the SVM approach, we do not need to model the assignment of labels to the leaf nodes.  \mnotee{This is because,} once hyperplanes $\H_t$ for $t\in\Tbl$ are defined, labels are then implicitly assigned to the leaves, with positive labels always assigned to right leaf nodes and negative labels to the left ones, as shown in Fig. \ref{fig: margottona}. Moreover, the modelling of the leaf level is usually needed to evaluate the misclassification error, which is usually computed "inside" the leaves, counting, with appropriate binary variables, the number of misclassified samples assigned to each leaf.  Nonetheless, in our case\mnotee{,} the misclassification cost is controlled by its upper bound defined by the sum of slack variables 
which do not depend on the leaf nodes. Thus, we will model the assignment of a sample $i$ only to a node in $\Tbl$, and this will be sufficient to reconstruct the unique path of the sample within the tree.


We can now define binary variables $ z_{i,t} $ for all $i\in \I$ and $t\in \Tbl$, as follows:

$$
\begin{array}{llll}
    z_{i,t} =
    \begin{cases} 
      1 & \text{if sample $i$ is \textit{assigned} to node $t\in \Tbl$} \\
      0 & \text{otherwise} 
   \end{cases}.
\end{array}
$$
Each sample has to be assigned to exactly one node $t\in\Tbl$, so we must impose that

\begin{align}
  &  \sum_{t\in\Tbl} z_{i,t} = 1  && \qquad \forall i\in\I \label{cons: MARGOT sum z} \\
  & z_{i,t} \in \{0,1\} && \qquad \forall i\in\I, \quad t\in\Tbl \label{cons: MARGOT z},
\end{align}

where constraints \eqref{cons: MARGOT sum z} will be also referred to as \textit{assignment constraints}. To model the SVM margin constraints \eqref{cons: margin cons I_t}, we observe that
whenever a sample $i$ belongs to $\I_t$, it must be assigned to a node in $\S^{\prime\prime}(t)$, hence we must have

$$\sum_{\ell\in\S^{\prime\prime}(t)} z_{i,\ell} =1 \qquad \qquad  \forall i\in\I_t.$$ 

We use this condition to activate or deactivate the SVM margin constraints when $i\in \I_t$ or $i\not\in \I_t$ respectively by means of a Big-M term.
Hence, we can  express the SVM constraints 
as

\begin{align}
  & y^i ( \ds w_t^T x^i + b_t) \ge 1 - \xi_{t,i} - M_{\xi} \Big (1-\ds \sum_{\ell\in\S^{\prime\prime}(t)} z_{i,\ell}  \Big ) && \forall i\in\I, \quad \forall t\in\Tb \label{cons: MARGOT class} \\
    & \xi_{t,i} \ge 0 && \forall i\in\I, \quad \forall t\in\Tb, \label{cons: xi non neg}
\end{align}

where $M_{\xi}>0$ is a sufficiently large value such that $  M_{\xi} \ge 1 - y^i ( \ds w_t^T x^i + b_t)$ is  satisfied for all $i\in\I$. When a sample $i\not\in \I_t$, margin constraints in \eqref{cons: MARGOT class} will always be satisfied, and variables $\xi_{t,i}$ at the optimum will be set to 0 because their sum is minimized in the objective function.

It remains to force each sample $i\in \I$ to follow a unique path from the root node to the node in $\Tbl$. As we have already commented, we must impose routing constraints only for the branch nodes in $\Tbf$. Indeed, the hyperplane at each node $t \in \Tbl$ is defined according to the soft SVM-type model using  \fnotee{$\xi$ variables} to measure the misclassification cost, and it does not depend on how the samples are finally routed in the leaves (namely on the predicted label). Thus, for each $t\in\Tbf$, we introduce the routing constraints modelling rules in \eqref{eq:routingrule} observing that 
 a sample $i\in\I_t$  following either the left or right branch from $t$, must satisfy 
$$\text{either}\quad\sum_{\ell\in\S^{\prime\prime}_L(t)} z_{i,\ell}=1 \quad\text{or}\quad \sum_{\ell\in\S^{\prime\prime}_R(t)} z_{i,\ell}=1 .$$ 

We can model the routing conditions
for each $ t \in \Tbf$ and for each $i\in\I$
using big-M constraints  as follows:
\begin{align}
 & \displaystyle w_t^T x^i + b_t \geq -M_{\H} \Big (1 - \sum_{\ell\in\S^{\prime\prime}_R(t)} z_{i,\ell}\Big) && \forall i\in\I,  \ \forall \  t \in \Tbf
 \label{cons: MARGOT right}
 \\
 & \displaystyle w_t^T x^i + b_t + \varepsilon \leq M_{\H} \Big (1 - \sum_{\ell\in\S^{\prime\prime}_L(t)} z_{i,\ell} \Big) && \forall i\in\I ,\   \forall \ t \in \Tbf \label{cons: MARGOT left},
\end{align}

where $ \varepsilon >0 $ is a sufficiently small positive value to model closed inequalities. We observe that when a sample $i\not\in \I_t$, we have that $\sum_{\ell \in \mathcal{S}^{\prime\prime}_L(t)} z_{i,\ell} = 
\sum_{\ell \in \mathcal{S}^{\prime\prime}_R(t)} z_{i,\ell}=0$ and both the constraints do not force any restriction on \mnotee{the sample}. Thus, each separating hyperplane $(w_t,b_t)$ \fnotee{is constructed using} a subset of samples.

The final MARGOT formulation is
the following:

\begin{align}
    (\text{MARGOT}) \quad
    \min \limits_{{w,b,\xi ,z} } \quad & \sum_{t \in \Tb } \Big ( \ds\frac 1 2 \|w_t\|^2_2 + C_t \sum_{ i\in\I} \xi_{t,i} \Big ) \nonumber\\ 
       \text{s.t.} \quad
    & y^i ( \ds w_t^T x^i + b_t) \ge 1 - \xi_{t,i} - M_{\xi} \Big (1-\ds \sum_{\ell\in\S^{\prime\prime}(t)}z_{i,\ell} \Big )  && \forall i\in\I, \quad \forall t \in \Tb \nonumber \\\  
     & w_t^T x^i + b_t \geq -M_{\H} \Big (1 - \sum_{\ell\in\S^{\prime\prime}_R(t)} z_{i,\ell} \Big) && \forall i\in\I, \quad \forall t \in \Tbf\nonumber  \\ 
     & w_t^T x^i + b_t + \varepsilon \leq M_{\H}\Big (1 - \sum_{\ell\in\S^{\prime\prime}_L(t)} z_{i,\ell} \Big) && \forall i\in\I, \quad \forall t \in \Tbf\nonumber  \\  
     & \ds\sum_{t\in\Tbl} z_{i,t} = 1 && \forall i\in\I \nonumber\\
    & \xi_{t,i} \geq 0 && \forall i\in\I, \quad \forall t \in \Tb \nonumber\\
     & z_{i,t} \in\{0,1\} && \forall i\in\I, \quad \forall t \in \Tbl. \nonumber
\end{align} 

\red{It is important to observe that the complexity of the tree, namely the number of effective splits, is implicitly controlled in MARGOT. 
Indeed, let us assume that the node $t$ does not split, namely that $w_t = 0$. 
The value of $b_t$ and $ \xi_{t,i}$ are set to appropriate values according to the SVM constraints \eqref{cons: margin cons I_t} that read as
$$y^i  b_t \ge 1 - \xi_{t,i}=\begin{cases}
    b_t \ge 1 - \xi_{t,i} &\text{if }y^i=1\cr 
    b_t \le  \xi_{t,i} -1 &\text{if }y^i=- 1 
\end{cases}.$$
Hence, the minimization of the misclassification cost will lead to 
$\xi_{t,i}=0$ for samples $i\in\I_t$ labelled with the most common label $\hat y$ in the node $t$, and, accordingly, $\hat y b_t\ge 1$.
For the samples $i\in\I_t$ with the minority label, the minimization of $ \xi_{t,i}$,  and the constraints $  \xi_{t,i} \ge 1+|b_t| $ and $ |b_t|\ge 1$
will lead to
 $b_t=\hat y$ 
 and $\xi_{t,i} = 2$.
Thus, if a sample $i$ is correctly classified, the misclassification cost related to node $t$ is $C_t\xi_{t,i} = 0$, otherwise $C_t\xi_{t,i} = 2C_t$.
In the special case when samples $i\in\I_t$
belong to the same class $\hat y$, then we get $w_t=0$, $\xi_{t,i}=0$ and we do not incur any cost in the objective function for that node. 
In this case, 
we 
have multiple solutions for $b_t$ that must satisfy $\hat y b_t\ge 1 $ for $i\in\I_t$.}


In Table \ref{tab: dimension}, we report a summary of the number of variables and constraints as a function of the depth $D$ of the tree. For the sake of simplicity, all the notation used in the definition of the model is reported in Table \ref{tab: notation} in the Appendix.


\begin{table}[ht]
\renewcommand\arraystretch{1.2}
\centering
\begin{tabular}{l|lc}

                                                       & Class & {Cardinality}                      \\ \hline \hline
\parbox[t]{2mm}{\multirow{3}{*}{\rotatebox[origin=c]{90}{\footnotesize Variables}}}&
\multicolumn{1}{l}{Continuous variables $(w,b,\xi)$} & \multicolumn{1}{c}{$(n+1+|I|)(2^{D}-1)$} \\[.8em]
&\multicolumn{1}{l}{Integer variables $z$}            & \multicolumn{1}{c}{$|I|\cdot 2^{D-1}$}   \\[.8em]
\hline
\parbox[t]{2mm}{\multirow{4}{*}{\rotatebox[origin=c]{90}{\footnotesize Constraints}}}&
\multicolumn{1}{l}{Routing constraints}               & \multicolumn{1}{c}{$2(2^{D-1}-1)|I|$}    \\[.8em]
&\multicolumn{1}{l}{Margin constraints}                & \multicolumn{1}{c}{$(2^{D}-1)|I|$}       \\[.8em]
&\multicolumn{1}{l}{Assignment constraints}            & \multicolumn{1}{c}{$|I|$}                \\\hline
\end{tabular}\caption{Summary of the dimensions of MARGOT.}
    \label{tab: dimension}
\end{table}
It is worth noticing that\mnotee{,} even in the case in which the binary variables $z_{i,t}$ are fixed to values $\hat z_{i,t}$ (e.g. by setting the values returned by another classification tree method such as CART), the subproblems solved at each $t\in\Tb $ are not pure SVM problems unless $t\in\Tbl$. Let us first note that sets $\I_t$, \red{for all $t\in\Tb$}, can be equivalently redefined as:

\begin{equation} 
     \I_t \coloneqq \Big\{ i\in\I : \sum_{\ell \in  \S^{\prime\prime}(t)}\hat{z}_{i,\ell} = 1 \Big\}.\nonumber
\end{equation} 
Similarly, sets $\I_{R(t)}$ and $\I_{L(t)}$, for all $t\in\Tbf$,
can be redefined as:
\begin{equation} \nonumber
    \I_{R(t)} \coloneqq \Big\{ i\in\I : \sum_{\ell \in  \S^{\prime\prime}_R(t)}\hat{z}_{i,\ell} = 1 \Big\}
    \quad \text{and} \quad
    \I_{L(t)} \coloneqq \Big\{ i\in\I : \sum_{\ell \in  \S^{\prime\prime}_L(t)}\hat{z}_{i,\ell} = 1 \Big\} .
\end{equation}
Thus, the MARGOT optimization problem  
decomposes into the resolution of $|\Tb|$ problems, where the first $|\Tbf|$ problems, one for each $t \in \Tbf$, are of the type:

\begin{align*}
    \min \limits_{{w_t,b_t,\xi_t} } \quad & \ds\frac 1 2 \|w_t\|^2_2 + C_t \sum_{ i\in\I_t} \xi_{t,i}  \nonumber\\ 
    \text{s.t.} \quad & y^i ( \ds w_t^T x^i + b_t) \ge 1 - \xi_{t,i}  && \forall i \in \I_t \nonumber \\\  
    & w_t^T x^i + b_t \geq 0 && \forall i\in\I_{R(t)} \nonumber\\ 
    & w_t^T x^i + b_t + \varepsilon \le 0   && \forall i\in \I_{L(t)} \nonumber  \\  
    & \xi_{t,i} \geq 0 && \forall i \in  \I_t \nonumber.
\end{align*}


Only the remaining $|\Tbl|$ problems are pure SVM problems in that routing constraints 
are not defined for the nodes of the last branching level.

\section{MARGOT with feature selection}\label{sec: margot with fs}

In Machine Learning, Feature Selection \mnotee{(FS)} is the process of selecting the most relevant features of a dataset. 
\mnotee{Among FS approaches, embedded methods} integrate feature selection in the training process.
In optimization literature, a solution is defined as sparse when the cardinality of the variables not equal to 0 is "low". The sparsity of an optimal solution is a requirement that is highly desirable in many application contexts.
As \mnotee{a} matter of fact, the concept of embedded feature selection translates \mnotee{into} the sparsity requirement for the solution of the optimization model used for the training process.

The MARGOT formulation does not
take into account the sparsity of the hyperplane coefficients variables $w_{t}$ for each node $t\in \Tb$. Thus, in order to improve the interpretability of our method, we propose two alternative versions of the MARGOT model where the number of features used at each branch node of the tree is either limited ("hard" approach) or penalized ("soft" approach). 
\mnotee{This way, the hyperplane at each branch node is induced to use only a subset of features.} 
This, together with the tree structure of the model, yields a hierarchy scheme on the subset of features which mostly affect the classification. 
In more detail, we introduce, for each node $t\in\Tb$ and for each feature \mnotee{$j=1,\dots,n$}, a new binary variable $s_{t, j} \in \{0,1\}$ such that:
$$
\begin{array}{llll}
    s_{t,j} =
    \begin{cases} 
      1 & \text{if feature $j$ is \emph{selected} at node $t$ ($w_{t,j}\neq0$)} \\
      0 & \text{otherwise}
   \end{cases}.
\end{array}
$$

Classical Big-M constraints 
on the variables $w_{t,j}$ must be added to model the above implication: 
\begin{align}
& -M_w s_{t,j} \leq w_{t,j} \leq M_w s_{t,j} && \forall t\in\Tb, \quad \forall j =1,\dots,n \label{cons: MARGOT HFS fs cons on w}  \\
& s_{t,j} \in\{0,1\} && \forall t\in\Tb, \quad \forall j =1,\dots,n,  \label{cons: MARGOT HFS s}
\end{align}
where $M_w$  is set to a sufficiently large value. 

Similarly to \cite{LabbeMixedInteger}, where a MILP feature selection version of the $\ell_1$-regularized SVM primal problem is proposed, in the hard features selection approach, we restrict the number of features used at each node to be not greater than a budget value.
We do that by introducing a hyperparameter $B_t$ and a budget constraint for each branch node $t\in\Tb$:

\begin{equation}
\ds\sum_{j=1}^n s_{t,j} \leq B_t.   \nonumber
\end{equation}

The resulting formulation for the hard version\mnotee{,} denoted as HFS-MARGOT\mnotee{,} is the following:

\begin{align}
    (\text{HFS-MARGOT}) \quad
    \ds \min \limits_{{w,b,\xi ,z, s} } \quad  & \sum_{t \in \Tb } \Big ( \ds\frac 1 2 \|w_t\|^2_2 + C_t \sum_{ i\in\I} \xi_{t,i} \Big ) \nonumber\\
    \text{s.t.} \quad 
    & \eqref{cons: MARGOT sum z} - \eqref{cons: MARGOT left} \nonumber \\
    & -M_w s_{t,j} \leq w_{t,j} \leq M_w s_{t,j} && \forall t \in \Tb, \quad\forall j=1,\dots,n 
    \nonumber\\
        & \ds\sum_{j=1}^n s_{t,j} \leq B_t && \forall t \in \Tb  \nonumber\\
& s_{t,j} \in\{0,1\} && \forall t\in\Tb, \quad \forall j =1,\dots,n \nonumber.  
\end{align}
 

In the soft approach, we remove the budget constraints and control their violations by adding a penalization term in the objective function weighted by an appropriate hyperparameter $\alpha$:

$$ \ds \sum_{t \in \Tb} \max\Big \{0, \sum_{j=1}^n s_{t,j} - B_t \Big \}.$$

The resulting version allows for more than $B_t$ features to be selected at each splitting node $t$ by penalizing the number of the features that exceed the budget. The $\max$ functions can be linearized with the introduction of a new continuous variable $u_t$, for all $t\in\Tb$, thus obtaining the following MIQP formulation denoted as SFS-MARGOT:

\begin{align}
    (\text{SFS-MARGOT}) \quad
    \min \limits_{{w,b,\xi ,z, s, u} } \quad & \sum_{t \in \Tb } \Big ( \ds\frac 1 2 \|w_t\|^2_2 + C_t \sum_{ i\in\I} \xi_{t,i}  + \alpha  u_{t} \Big )  \nonumber\\
    \text{s.t.} \quad & \eqref{cons: MARGOT sum z} - \eqref{cons: MARGOT left} \nonumber\\
    & -M_w s_{t,j} \leq w_{t,j} \leq M_w s_{t,j} && \forall t \in \Tb, \quad\forall j=1,\dots,n  \nonumber\\
    & u_t \geq \sum_{j=1}^n s_{t,j} - B_t  && \forall t \in \Tb \nonumber\\ 
    & u_t \geq 0 && \forall t\in\Tb  \nonumber \\
    & s_{t,j} \in\{0,1\} && \forall t\in\Tb, \quad \forall j =1,\dots,n \nonumber.  
\end{align}


Inducing \emph{local} sparsity on each vector $w_t$ may be preferable rather than addressing \emph{global} sparsity on the full vector $w$, as done in the OCT-H model proposed in \cite{Bertsimas2017OptimalClassification}. Indeed, global sparsity of the vector $w$ has little effect on the "spreadness" of the features among the splitting rules in the tree, often leading to trees with fewer and less interpretable splits, usually at the higher levels. 
\lnote{This way, a local approach can generate models that better exploit the tree's hierarchical structure.} A solution that is more "spread" in terms of features can thus result in a more interpretable machine learning model because it yields a hierarchy scheme among the features which mostly affect the classification. In this sense, HFS-MARGOT attains a locally sparse tree classifier, while SFS-MARGOT is more of a hybrid between HFS-MARGOT and OCT-H, and, depending on the choice of parameters $\alpha$ and $B_t$, $t \in \Tb$, it can be regarded as a more local or global approach. Of course, other constraints facing additional requirements on the selected features can be added to these formulations. Indeed, the sparsity of the $w_t$ variables may not be the only interesting property when addressing the interpretability of the decision.

Table \ref{tab: variables} provides an overview of all the variables used in the MARGOT formulations.

{\renewcommand\arraystretch{1.2} 
\begin{table}
    \centering
\begin{tabular}{lll}
Variable & Meaning &Model\\
\hline \hline
${w}\in \mathbb{R}^{|\Tb| \times n}$ & split coefficients &all\\

${b} \in \mathbb{R}^{|\Tb|}$ & split biases &all \\
${\xi} \in \mathbb{R}^{|\Tb| \times |\I|}$ & slack variables &all\\
${z} \in \{0,1\}^{|\I| \times |\Tbl|}$ & {samples assignment  to nodes in $\Tbl$} & all \\
${s} \in \{0,1\}^{|\Tb| \times n}$ & feature selection  &HFS/SFS-MARGOT\\
${u} \in \mathbb{R}^{|\Tb|}$ & soft FS penalization parameter
&SFS-MARGOT \\ \hline
\end{tabular}
\caption{Overview of all the variables used in the MARGOT formulations and their meaning.}
\label{tab: variables}
\end{table}}

\section{Heuristic for a starting feasible solution}\label{sec: heuristic_ws}

We  develop a simple greedy heuristic algorithm to find a feasible solution to be used as a good-quality warm start
for the optimization procedure.
As well known, the value of the warm start solution is an upper bound on the optimal one, and it can be used to prune nodes of the branch and bound tree of the MIP solver, eventually yielding shorter computational times.
Thus, developing a good warm start solution is usually addressed in MIP formulations and implemented in off-the-shelf \fnotee{solvers} at the root node of the branching tree. 
In \cite{Bertsimas2017OptimalClassification}, several warm start procedures are proposed, from the simplest one, which consists \mnotee{of} using the solution provided by CART, to more tailored ones. 

The general heuristic scheme, denoted as Local SVM Heuristic, exploits the special structure of the MIQP models addressed and can be applied to obtain feasible solutions for MARGOT, HFS-MARGOT, and SFS-MARGOT models. 
The Local SVM Heuristic is based on a greedy recursive top-down strategy. Starting from the root node, 
the maximum margin hyperplane is computed using an SVM model embedding, when needed, feature selection constraints and penalization. 
More in detail, for each node $t \in \Tb$, 
given a predefined index set $ \I_t \subseteq \I$, the heuristic solves the following problem:

\begin{align*}
    (\mnotee{\text{WS-SVM$_t$}}) \quad
    \min \limits_{{w_t,b_t,{\xi_t}},s_t,u_t} \quad & \ds\frac 1 2  \|w_t\|^2_2 +C_t \sum_{i\in \I_t} \xi_{t,i} + \alpha u_t \nonumber\\
    \text{s.t.} \quad 
    & y^i(w_t^Tx^i+b_t)\ge  1 - \xi_{t,i} && \forall i \in \I_t \nonumber\\
    & -M_w s_{t,j} \leq w_{t,j} \leq M_w s_{t,j} && \forall j=1,\dots,n \nonumber\\
    & u_t \geq \sum_{j=1}^n s_{t,j} - B_t\nonumber\\
    & u_t \geq 0  \nonumber\\
    & \xi_{t,i} \geq 0 &&\forall i \in \I_t \nonumber\\
    & s_{t,j} \in\{0,1\} && \forall j=1,\dots,n, \nonumber
\end{align*}

where \mnotee{hyperparameters} $ B_t, \alpha$ and variable \fnotee{$u_t$} may be fixed to specific values to get warm start solutions
for the three different models. 
In particular, when $B_t=n$ we do not impose restrictions on the number of features, and variable \fnotee{$u_t$} will be automatically set to $0$. When  $B_t<n$ and $u$ is set to zero\mnotee{,} we obtain an $\ell_2$-regularized SVM model with a hard constraint on the number of features, \mnotee{similar} to the approach in \cite{LabbeMixedInteger}. Finally, when $\alpha>0, B_t<n$ and \fnotee{variable $u_t$} 
\fnotee{is not fixed}, we impose a soft constraint on the number of features.
Given the optimal tuple $(\widehat w_t,\widehat b_t, \widehat \xi_t, \widehat s_t) \in \R^{n}\times \R \times \R^{|\I|} \times \{0,1\}^n  $ obtained at node $t$, the samples are partitioned to the left or right child node in the subsequent level of the tree according to the routing rules defined by the hyperplane $ \H_t = \{x \in \R^n : \widehat w_t^T x+\widehat b_t = 0\}$. Thus, each node $t$ works on a different subset of samples $\I_t\subseteq \I$, and $\I_{L(t)}$ \mnotee{and} $\I_{R(t)}$  are the index sets of samples assigned to the left and right child node of $t$\mnotee{, respectively}. At the end of the procedure, for each $t\in\Tb$, the solutions $(\widehat w_t,\widehat b_t,\widehat \xi_{t}) \in \R^{n}\times \R \times \R^{|\I|} $, together with solutions $\widehat s_t \in \{0,1\}^n$, when needed, constitute a feasible solution for the original problem. In the very last step, variables $z_{i,t}$ are set according to the definitions of sets $\I_t, \, t\in\Tbl$.
The general scheme encompassing the three different strategies is reported in Algorithm \ref{alg: local_svm}.

\begin{algorithm}[ht]\caption{Local SVM Heuristic}

	\DontPrintSemicolon
	\smallskip
	
	
	{\bf Data:}  $\{(x^i,y^i) \in \mathbb{R}^n \times \{-1,1\}, i \in \I\}$; \medskip
	
	\smallskip
	{\bf Parameters:}  $\{C_t>0, t\in\Tb\}$, $\widehat \alpha>0$,  $M_w>0$, $D$, $\varepsilon>0$, $\{\widehat B_t>0 , t\in\Tb\}$;

	\smallskip
	{\bf Input:} Model $\in \{\text{MARGOT, HFS-MARGOT, SFS-MARGOT}\}$;
	
	\medskip
	{\bf Initialize:} $\I_0 = \I$,  $\I_t  = \emptyset \; \forall t\in\Tb \setminus\{0\}$, $\widehat z_{i,t} = 0,$ $\forall t \in \Tbl, \forall i\in\I$, $\widehat \xi_{t,i} = 0,$ $\forall t \in \Tb, \forall i\in\I$;
	\BlankLine
	\For{level $k = 0,\dots,D-1$}{
	    \For {node $t =2^{k}-1,\dots, 2^{k+1}-2$}{\medskip
\If {model = \textup{MARGOT} }{
set $B_t=n$}

\If {model = \textup{HFS-MARGOT} }{
set $B_t=\widehat B_t$ and \fnotee{$u_t=0$}
}
\If {model = \textup{SFS-MARGOT}}{
set $B_t=\widehat B_t$ and $\alpha=\widehat \alpha$} 

\medskip

	       Find $(\widehat w_t, \widehat b_t,\widehat \xi_t, \widehat s_t)$ optimal solution of WS-SVM$_t$ 
	       \medskip
            
Set $\I_{L(t)} = \{i\in\I_t: \widehat{w}_t^Tx^i+\widehat{b}_t +\varepsilon \leq 0\}$ and $ \I_{R(t)} = \{i\in\I_t:  \widehat{w}_t^Tx^i+\widehat{b}_t \geq 0\}$
    }
	}
 \For {$t \in \Tbl$}{
\If {$i\in\I_t$}{
      Set $\widehat z_{i,t} = 1$ 
}
}\medskip

{\bf Output:} 
A feasible solution $(\widehat w, \widehat b,\widehat \xi, \widehat s, \widehat z)$ for all the input models.

\label{alg: local_svm}
\end{algorithm}

The heuristic procedure requires the solution of $2^D-1$ MIQP problems with a decreasing number of constraints (depending on the size of $\I_t$) that can be easily handled by off-the-shelf MIP solvers. We show in the computational experiments that the use of the \mnotee{proposed} heuristics improves the quality of the first incumbent solution with respect to the chosen optimization solver.

\section{Computational Results}\label{sec: results}

In this section\mnotee{,} we present different computational results where models MARGOT, HFS-MARGOT and SFS-MARGOT are compared to other {three} benchmark OCT models:
 \begin{itemize}
    \item  \red{OCT-1 and OCT-H, the \mnotee{traditional} univariate and multivariate optimal classification tree models proposed in \cite{Bertsimas2017OptimalClassification}};
     \item MM-SVM-OCT, as proposed in \cite{BlancoRobust},
     where no \mnotee{relabelling} is allowed.
\end{itemize}

\mnotee{Note that there are alternative methods for constructing optimal univariate trees based on dynamic programming algorithms (\cite{Aglin2020LearningOptimalDecision, Lin2020GeneralizedAndScalable}). However, to ensure a fair comparison, we select OCT-1 as a standard benchmark for univariate trees so that all the tested approaches rely on solving MIP formulations using the same optimization solver.}
\mnotee{Moreover, in the case of MM-SVM-OCT,} we did not allow relabelling, which is used to make the method robust against noisy data. Indeed, our aim is to evaluate the isolated effect of their way of constructing margin-based splits in the tree without considering the potential effects relabelling may introduce.


All mathematical programming models have been implemented on our own. They were coded in Python and solved using Gurobi 9.5.1 on a server Intel(R) Xeon(R) Gold 6252N CPU processor at 2.30 GHz and 96 GB of RAM. 
The source code and the data of the experiments are available at \url{https://github.com/m-monaci/MARGOT}, and additional implementation details are provided in the next sections.

\noindent We used two groups of datasets:
 

 \begin{itemize}
 \item  3 non-linearly separable synthetic datasets in a 2-dimensional feature space, in order to give a graphical representation of the maximum margin approach {(presented in section \ref{sec: results_synthetic})};
     \item  10 datasets from UCI Machine Learning Repository \cite{UCI2019}, to assess the effectiveness of the formulations as regards both the predictive and the optimization performances {(presented in section \ref{sec: results_UCI})}.
 \end{itemize}

{As \mnotee{regards} categorical data, we treated ordinal attributes as numerical ones, while we \mnotee{applied} the standard one-hot encoding for nominal features. We normalized the feature values of each dataset to the 0-1 interval}. For the results on the UCI datasets, we performed a $4$-fold \mnotee{cross-validation} \red{to select the best hyperparameters} 
which is detailed in the specific sections below. In Section \ref{sec: warm start results}, we eventually present a brief analysis on the Local SVM algorithm presented in Section \ref{sec: heuristic_ws}, in order to motivate the use of the warm start solution in input to the solver. \red{In all tables reporting predictive and optimization performances, the best result is highlighted in bold, while, when the time limit was reached, the time value is underlined.}
\red{The interested reader can also refer to the \ref{app: appendix1} for other insightful results that were omitted here to avoid excessive verbosity.}

\subsection{Results on 2-dimensional synthetic datasets}\label{sec: results_synthetic}

As regards the 2-dimensional datasets, 
we used two artificially constructed problems, \texttt{4-partitions} and \texttt{6-partitions}, and the more complex synthetic \texttt{fourclass} dataset 
 \cite{ho1996building} as reported in LIBSVM Library \cite{LIBSVM}.
Our aim is to offer a glimpse of the differences in the hyperplanes generated by the different \red{multivariate} optimal tree models, \red{thus OCT-1 was not compared here}. We also reported the solution returned by the Local SVM Heuristic (Algorithm \ref{alg: local_svm}) to show how far the greedy solution is from the optimal ones. For all three synthetic datasets, there \mnotee{exists} a set of hyperplanes that can reach perfect classification on the training data. In particular,  \texttt{4-partitions} and \texttt{6-partitions} were constructed by defining hyperplanes with margins and plotting 108 and 96 random points, respectively, in regions outside the margin (see Figures \ref{fig: ground truth 4 parts}, \ref{fig: ground truth 6 parts}).
Although reaching zero classification error on the training data is not desirable in ML models, in these cases, we want to highlight the power of using hyperplanes with margins to derive more robust classifiers. We did not account for the out-of-sample performance; therefore, the \mnotee{entire} datasets \mnotee{were} used to train the optimal tree.

Of course, because these datasets are in a 2-dimensional space, we do not consider HFS-MARGOT and SFS-MARGOT. The results are commented below, and a cumulative view is \mnotee{provided} in Table \ref{tab:2D}. \mnotee{For} all the experiments, the time limit of the solver has been set to 4 hours.

In Figures  \ref{fig: OCTs 4 parts}, \ref{fig: OCTs 6 parts}, \ref{fig: fourclass complete}, we report the hyperplanes generated by the Local SVM Heuristic, OCT-H, MM-SVM-OCT, and MARGOT. 
Different \mnotee{colours} correspond to different branch nodes of the tree, as reported in the legend.

{In the case of MM-SVM-OCT and MARGOT, at each last splitting node, we plotted the two supporting hyperplanes at a distance of $2||w||^{-1}$ to highlight the margins of the hyperplanes that define the predicted class of samples. Of course, such supporting hyperplanes can \mnotee{also be plotted} for the other splitting nodes. 
\fnotee{However, we omit them because the hyperplane margins at these nodes are very wide, and highlighting them may be confusing and not provide much insight for the reader.}

\begin{figure}[ht!]
\centering
\begin{subfigure}[c]{.264\textwidth}
  \centering
  \includegraphics[width=1\linewidth]{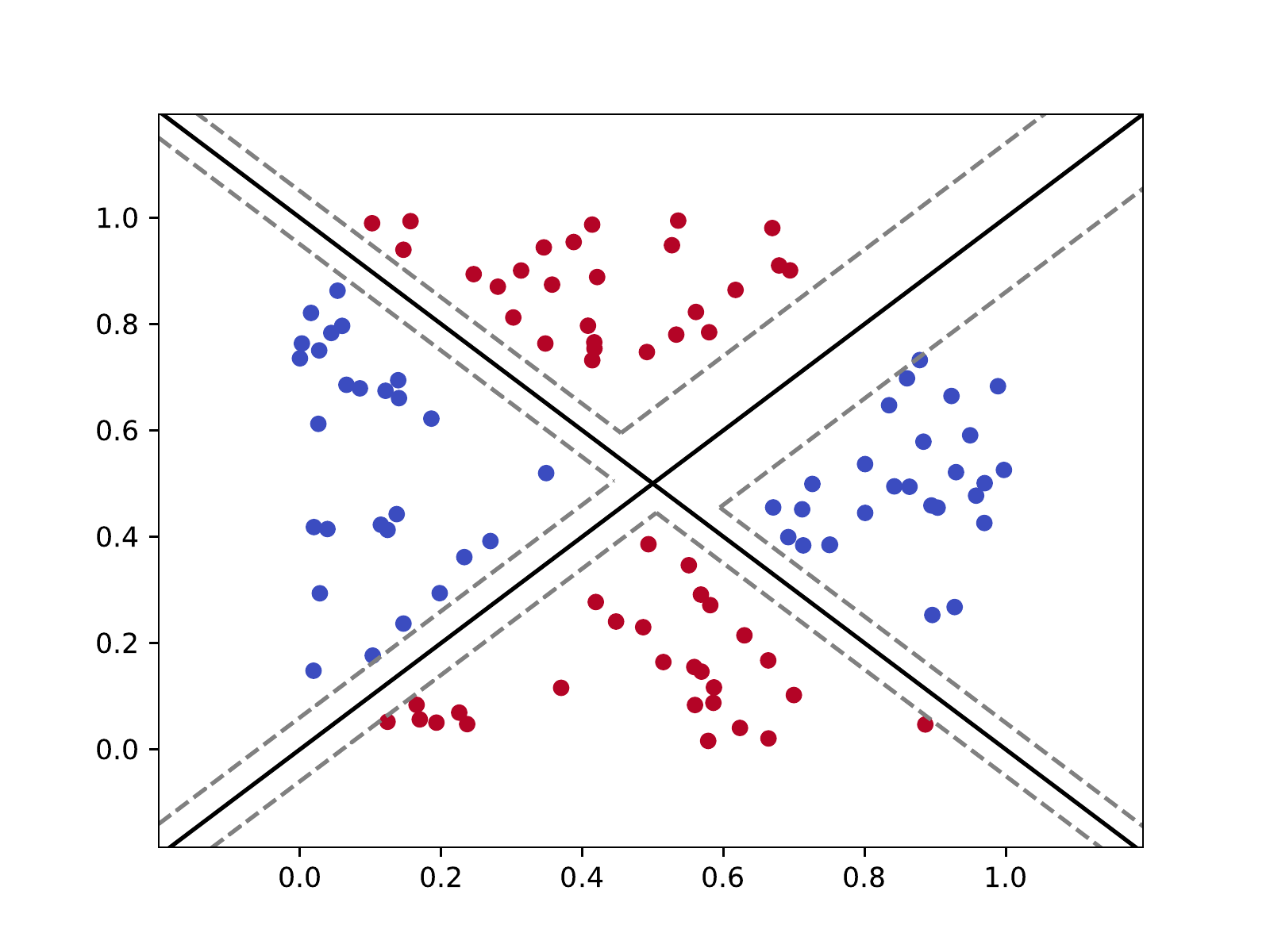}
  \caption{Ground truth}
  \label{fig: ground truth 4 parts}
\end{subfigure}%
\begin{subfigure}[c]{.528\textwidth}
  \centering 
  \begin{tabular}{cc}
    \includegraphics[width=.5\linewidth]{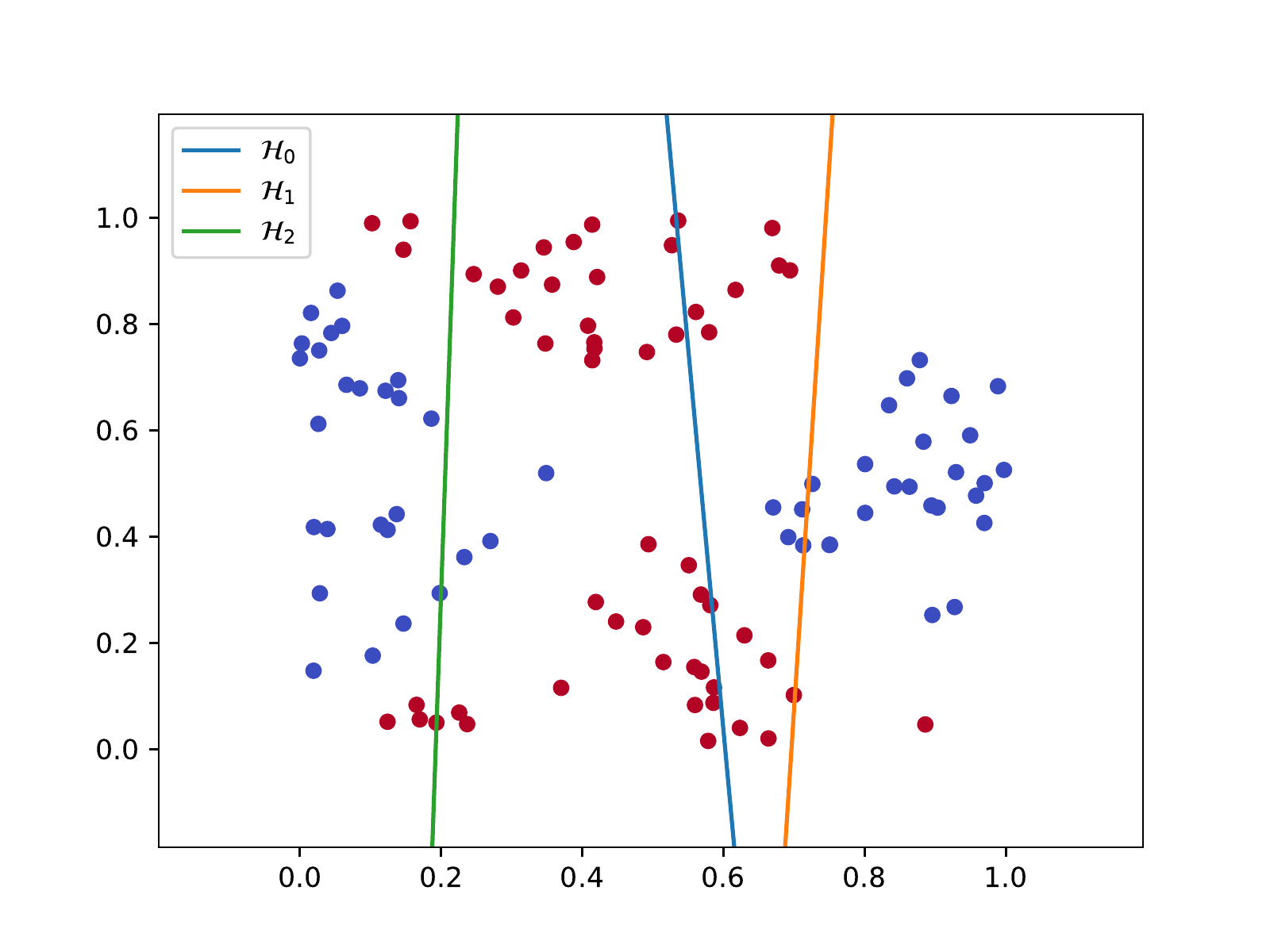} & 
    \includegraphics[width=.5\linewidth]{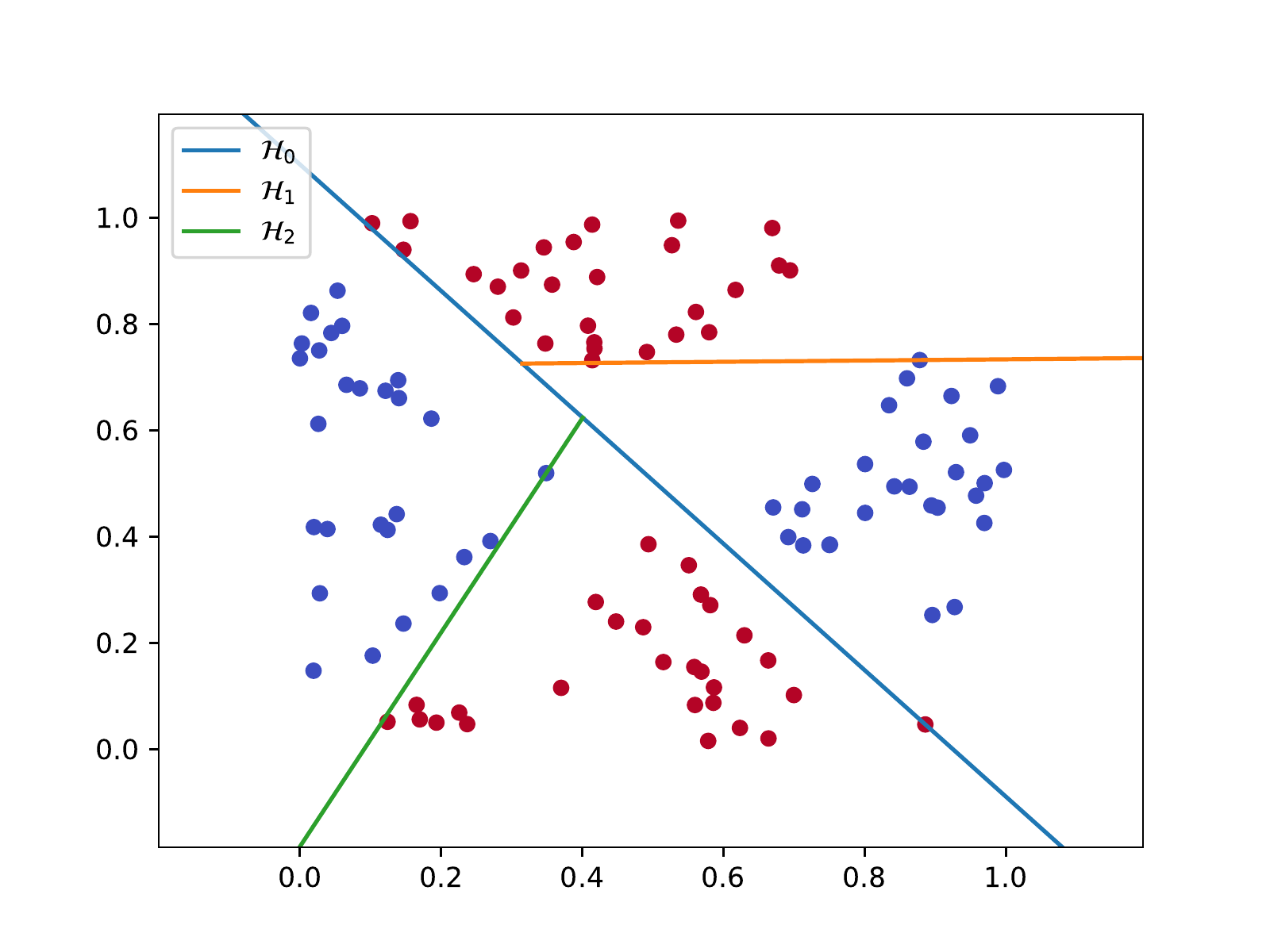}\\
    \footnotesize{(i) Local SVM Heuristic} & \footnotesize{ (ii) OCT-H} 
  \end{tabular}  
  \begin{tabular}{cc}
    \includegraphics[width=.5\linewidth]{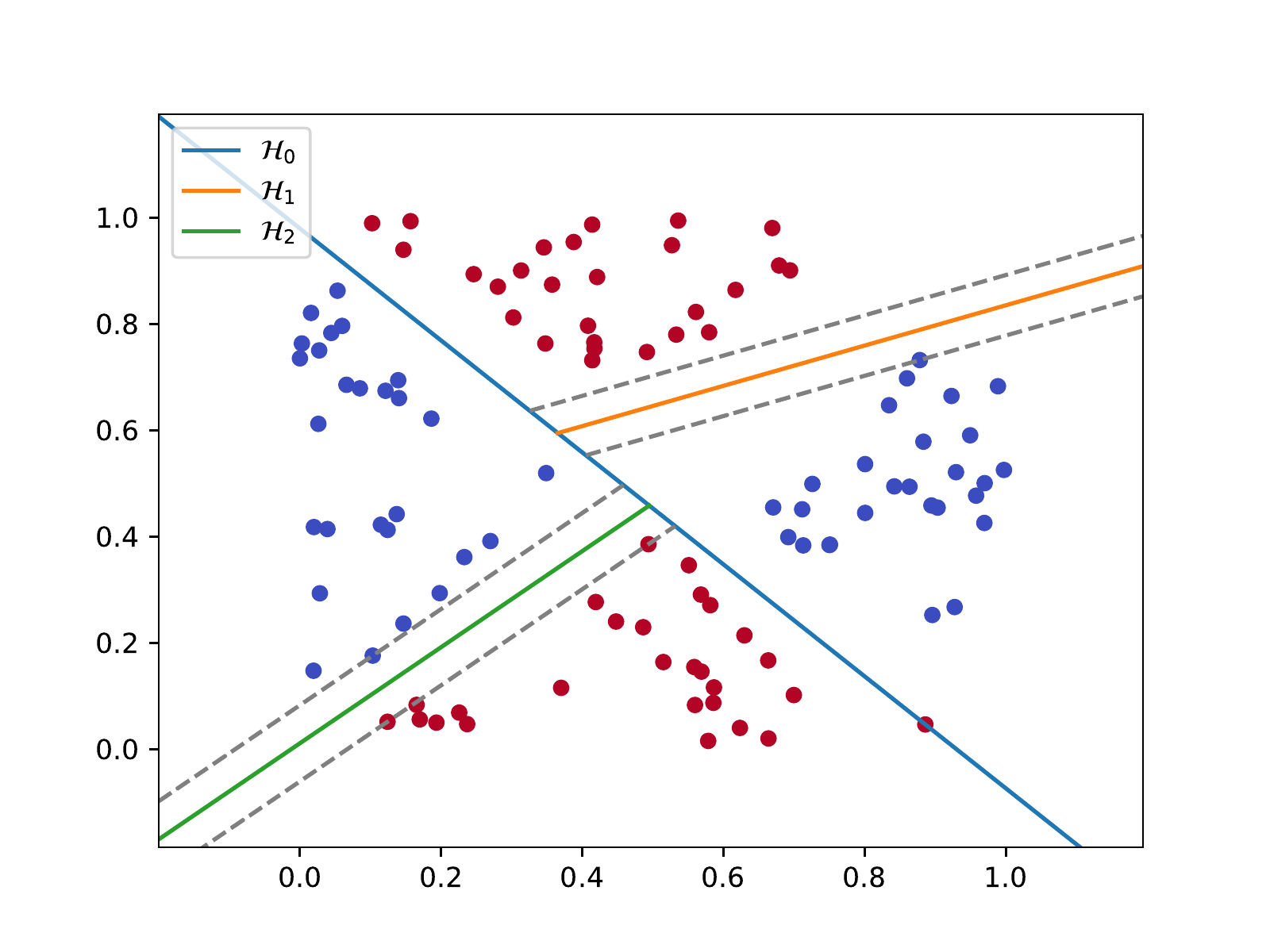}& 
    \includegraphics[width=.5\linewidth]{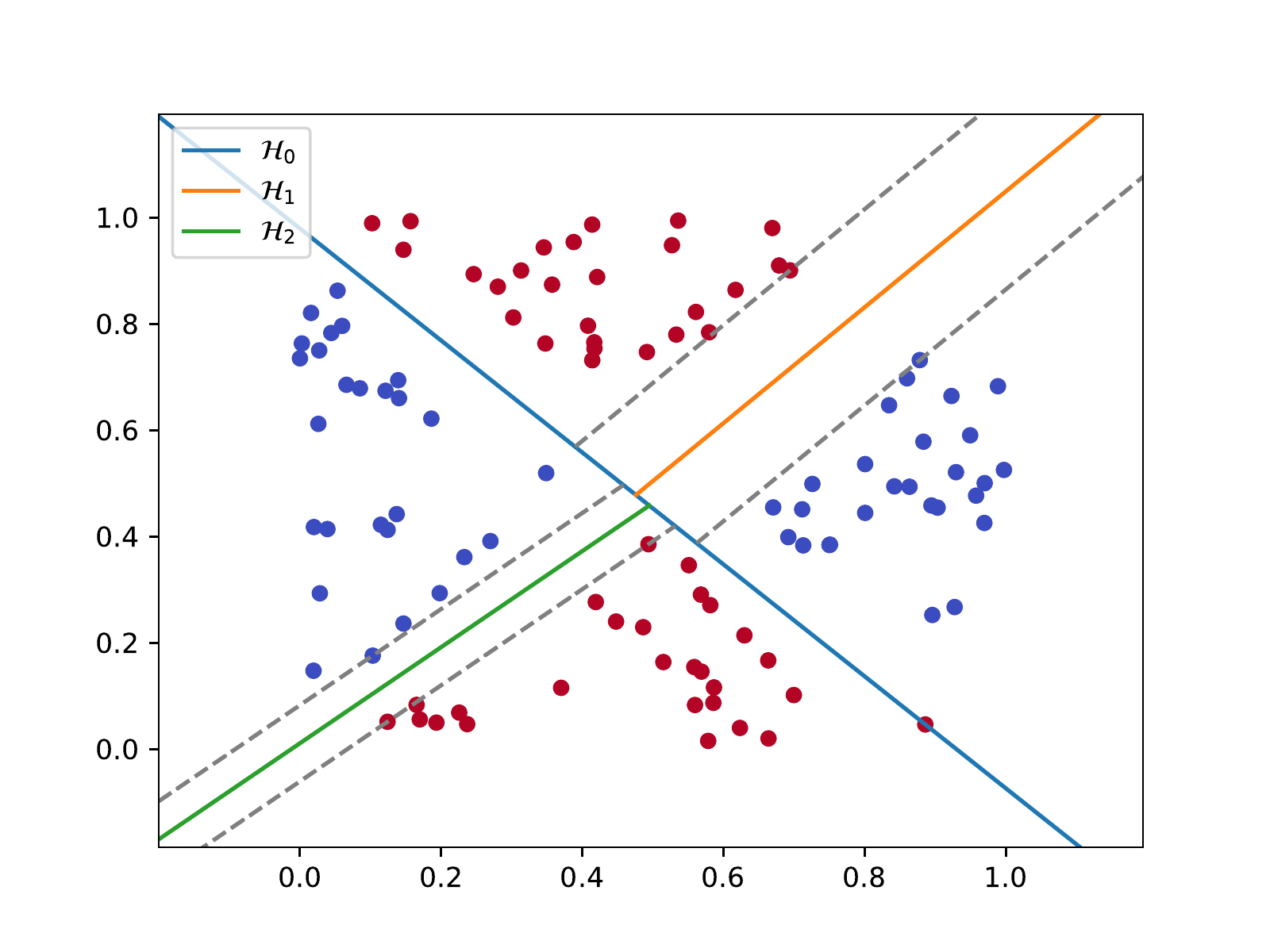}\\
    \footnotesize{(iii) MM-SVM-OCT} & \footnotesize{(iv) MARGOT}
  \end{tabular}
  \caption{}
  \label{fig: OCTs 4 parts}
\end{subfigure}
\caption{Results on the \texttt{4-partitions} synthetic dataset.}
\label{fig: 4 parts complete}
\end{figure}

\begin{figure}[ht!]
\centering
\begin{subfigure}[c]{.264\textwidth}
  \centering
  \includegraphics[width=1\linewidth]{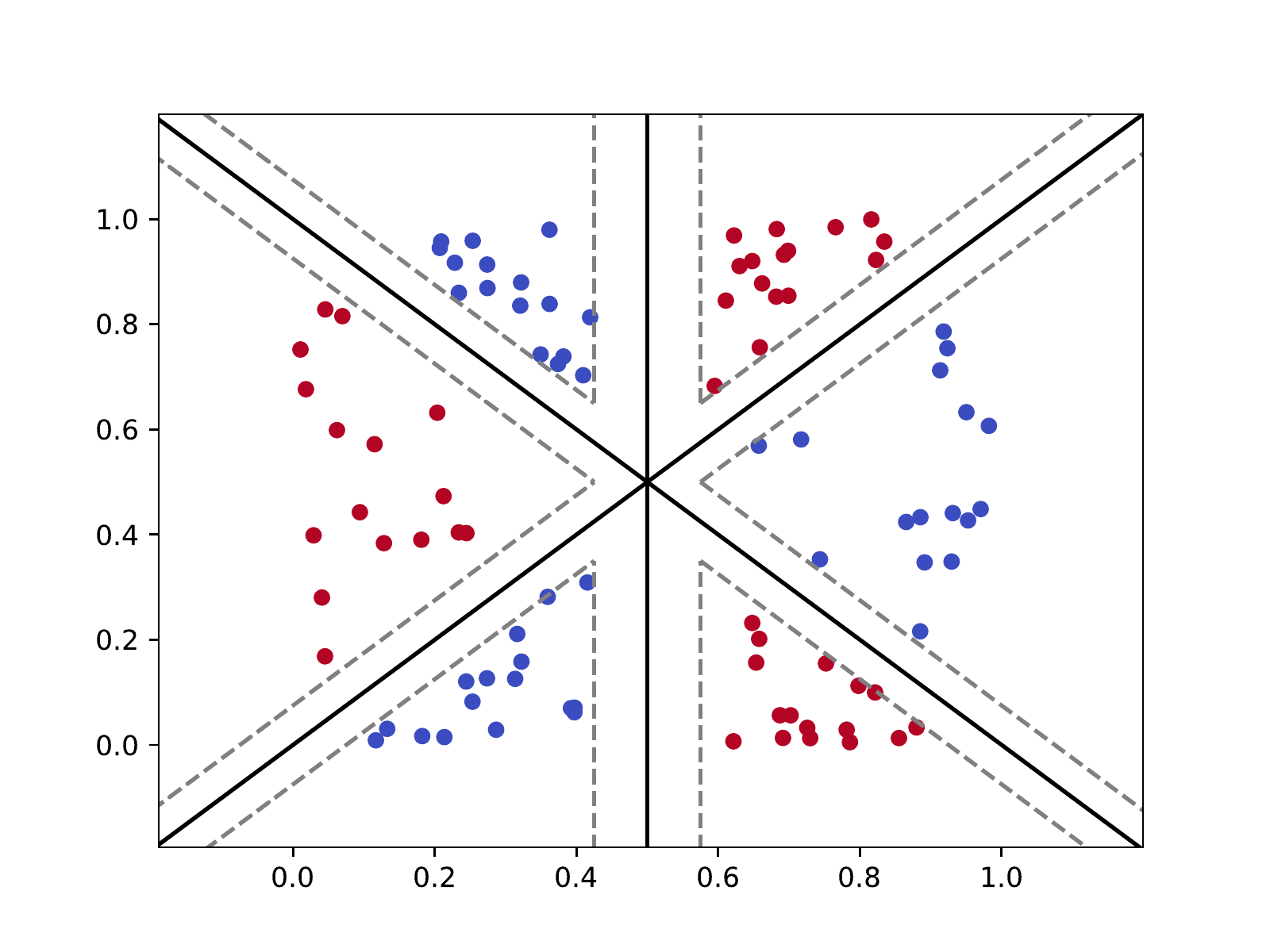}
  \caption{Ground truth
  }
  \label{fig: ground truth 6 parts}
\end{subfigure}%
\begin{subfigure}[c]{.528\textwidth}
\centering 
    \begin{tabular}{cc}
    \includegraphics[width=.5\linewidth]{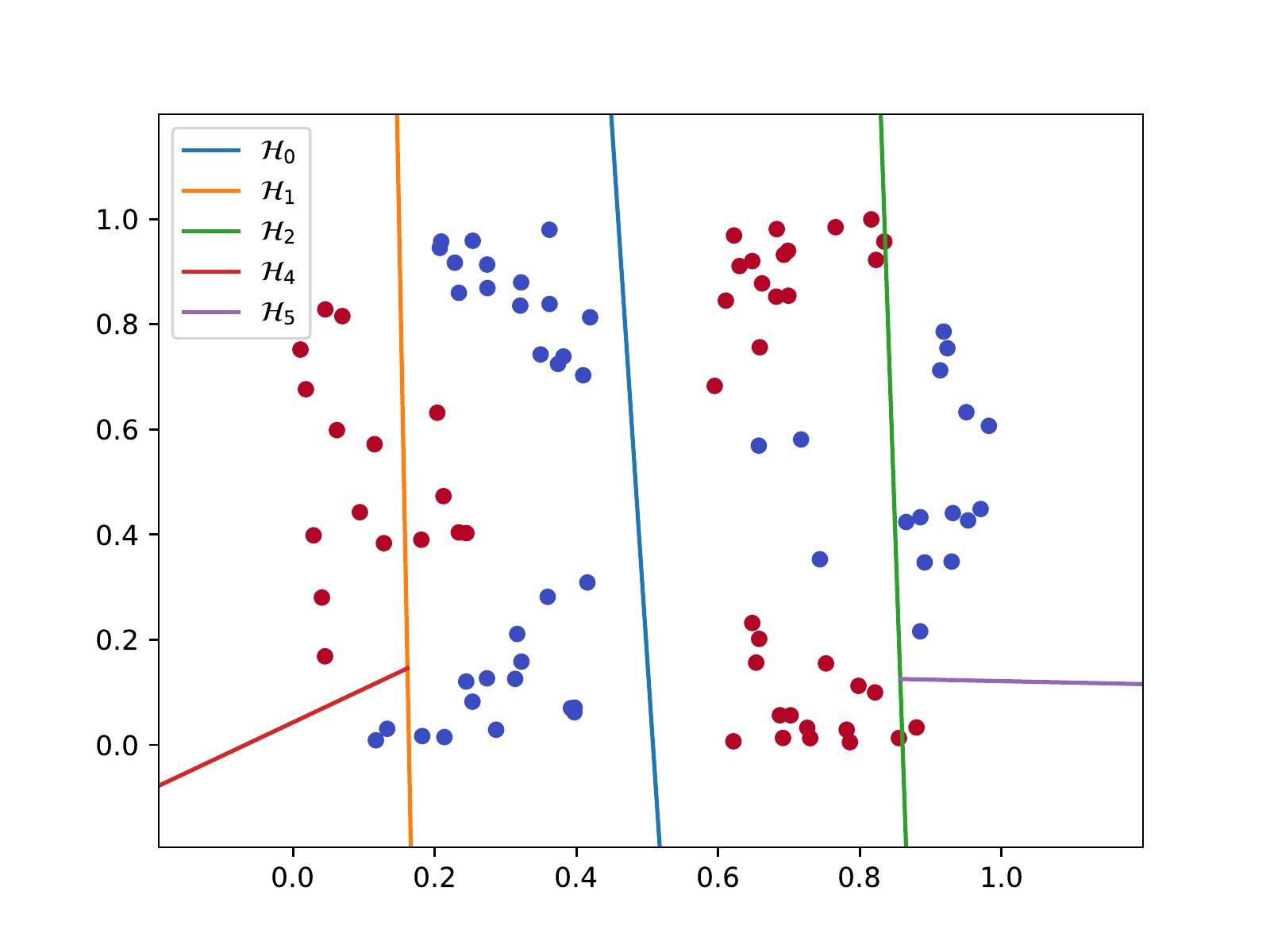}  & 
        \includegraphics[width=.5\linewidth]{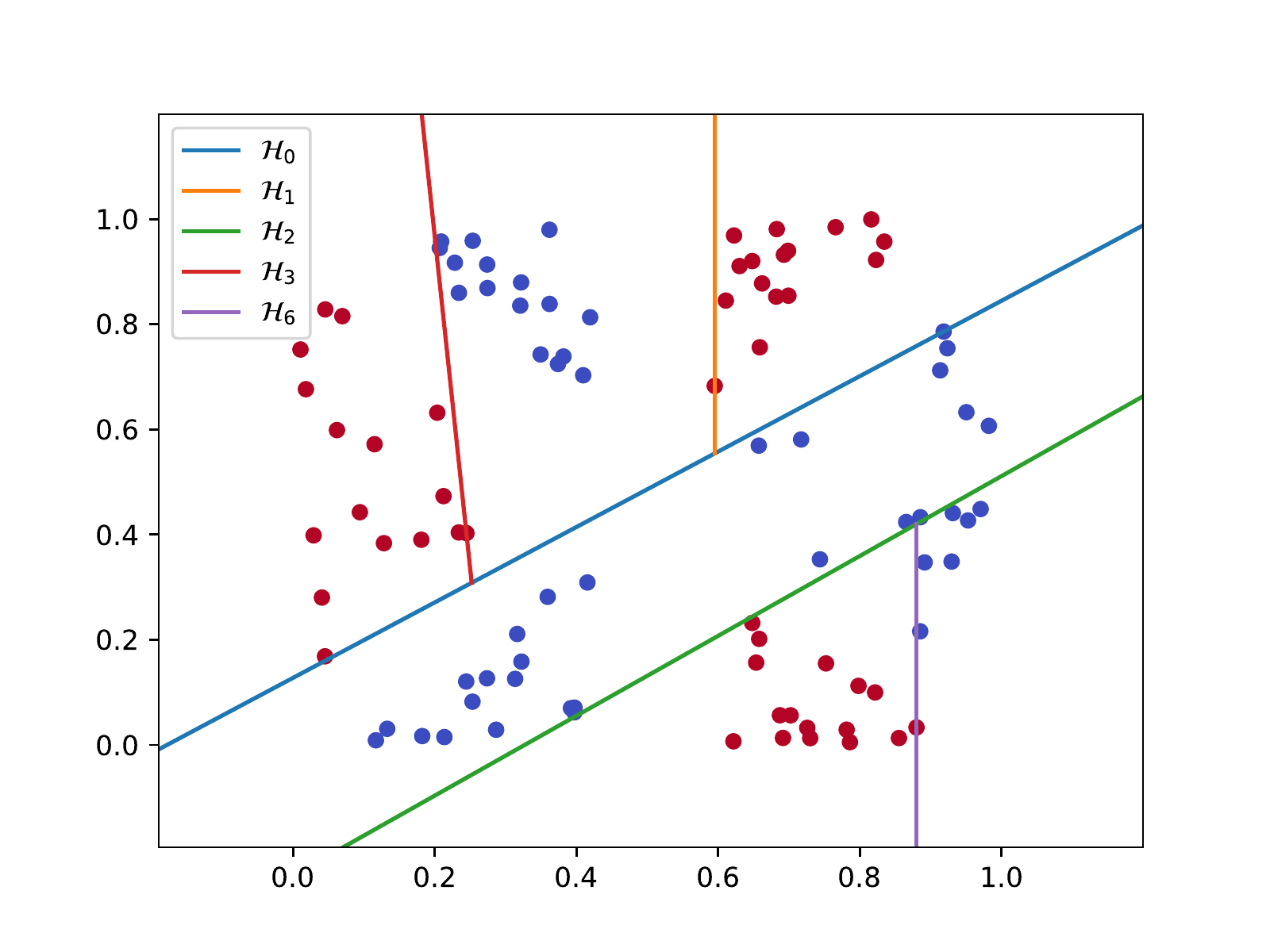}\\
            \footnotesize{(i) Local SVM Heuristic} & \footnotesize{ (ii) OCT-H} 
    \end{tabular}
    
    \begin{tabular}{cc}
    \includegraphics[width=.5\linewidth]{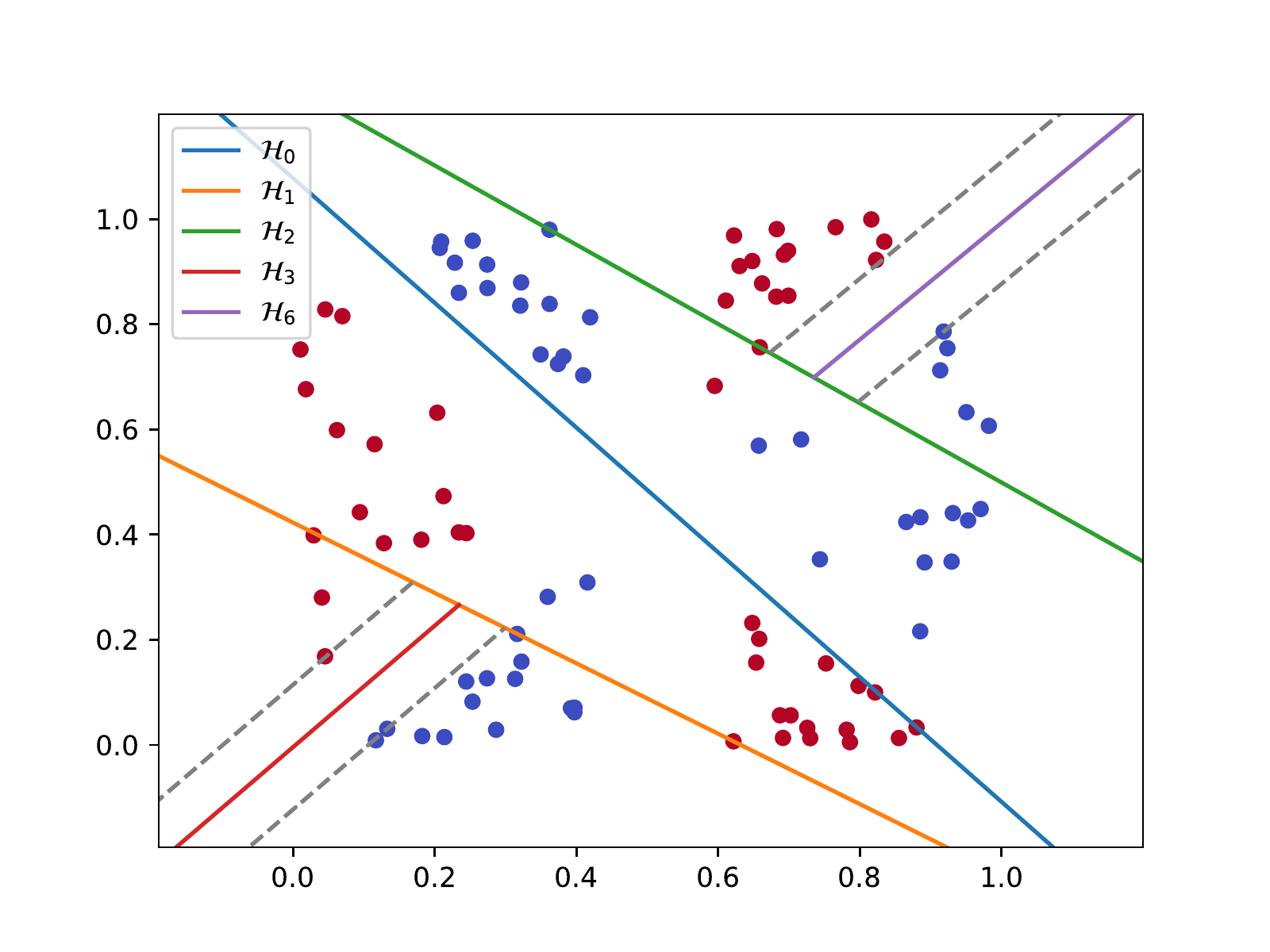}& 
    \includegraphics[width=.5\linewidth]{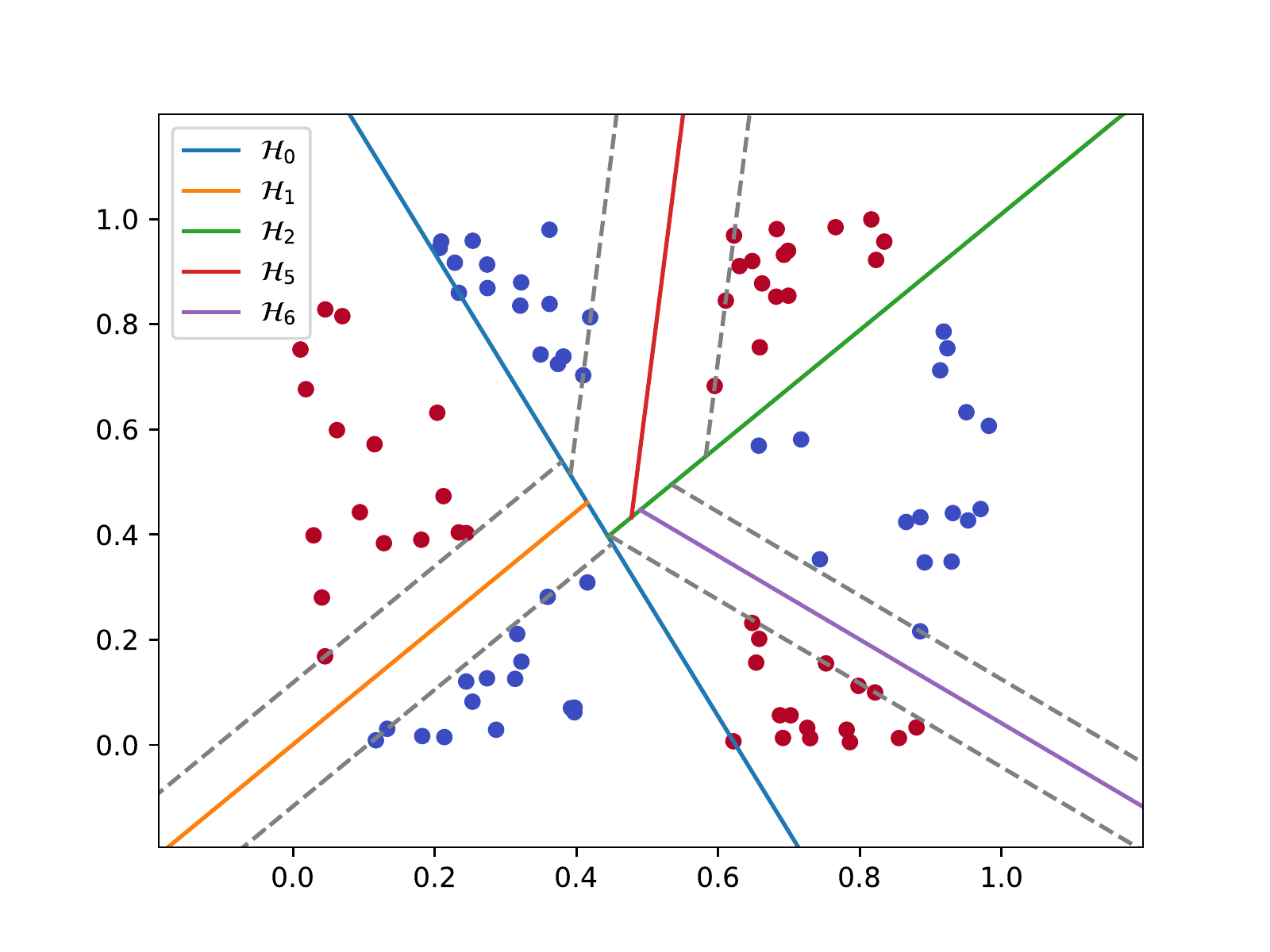}\\
    \footnotesize{(iii) MM-SVM-OCT} & \footnotesize{(iv) MARGOT}
    \end{tabular}
    \caption{}
    \label{fig: OCTs 6 parts}
 
\end{subfigure}
\caption{Results on the \texttt{6-partitions} synthetic dataset.}
\label{fig: 6 parts complete}
\end{figure}

For the \texttt{4-partitions} dataset\mnotee{,} we consider trees with depth $D = 2$. Fig. \ref{fig: OCTs 4 parts} (i) represents the tree obtained by the Local SVM Heuristic. Hyperplane $\H_0$ at the root node \mnotee{coloured} in blue is obtained on the whole dataset, while the hyperplanes at its child nodes, $\H_1$ in green and $\H_2$ in orange, are obtained considering the partition of the points given by $\H_0$ as the splitting rule on the whole dataset. 
The heuristic returns a classification tree that does not classify all data points correctly, obtaining an accuracy of $86.1\%$.
Concerning the OCT approaches in Fig. \ref{fig: OCTs 4 parts} (ii), (iii), and (iv), the solver 
certified \mnotee{the} optimal solution on all three models, thus obtaining 0\% MIP gap in different computational times.
All three OCT models reach an accuracy of $100\%$. We note that OCT-H creates hyperplanes that do not consider the margin. Indeed, the objective of this approach is to minimize the misclassification cost and the number of features used across the whole tree. Thus, as it happens for the orange hyperplane $\H_1$, OCT-H tends to select axis-aligned hyperplanes to split the points. As regards the tree produced by the MM-SVM-OCT model, only the minimum margin among all the hyperplanes is maximized. Consequently, only the green hyperplane $\H_2$ lies in the centre between the partitions of points, while the others do not. Instead, the MARGOT tree is the one that most resembles the ground truth in Fig. \ref{fig: ground truth 4 parts} 
and both the $\H_1$ and $\H_2$ hyperplanes have a wider margin. 

Fig. \ref{fig: 6 parts complete} shows the more complex synthetic dataset \texttt{6-partitions} where we set $D = 3$, the minimum depth to correctly classify all samples.
MARGOT and OCT-H reach perfect classification on the whole dataset, and both Local SVM Heuristic and MM-SVM-OCT return good accuracies above 90\%.
Moreover, the solutions produced by MARGOT and MM-SVM-OCT models are optimal, while OCT-H optimization reaches the time limit with a MIP gap of 50\%.
\fnotee{In this case as well, MARGOT appears to produce the most reliable classifier among all the generated trees, as it is the one closest to the ground truth.}


Finally, we evaluated the four methods on the \texttt{fourclass} dataset. 
In this case, being the problem the most complex of the three, none of the models has been solved to certified optimality, and OCT-H, MM-SVM-OCT and MARGOT optimization procedures reach a MIP gap of 100\%, 74.5\% and 67.2\% respectively.
MARGOT and OCT-H approaches were able to correctly classify almost all the samples, outperforming MM-SVM-OCT, which reaches an accuracy of 88.4\%. 
\red{It is possible to observe how the greedy fashion of the Local SVM Heuristic may generate models not able to capture the underlying truth of the data. Indeed, when applying local SVMs after the split at the root node, it might happen that producing splits is not "convenient" in terms of the objective function. This is due to the 
  highly nonlinear separability of the dataset, which cannot be effectively handled by a single linear SVM. 
  Thus, when applying the Local SVM Heuristic, not all the possible 15 splits are generated.
   This case illustrates the drawbacks of the greedy methods compared to optimal ones: when applied to highly non-linearly separable datasets, these heuristic approaches lead to myopic decisions resulting in poor predictive performances.}

\begin{figure}[ht]
    \centering 
    \begin{tabular}{ccc}
\includegraphics[width=.264\linewidth]{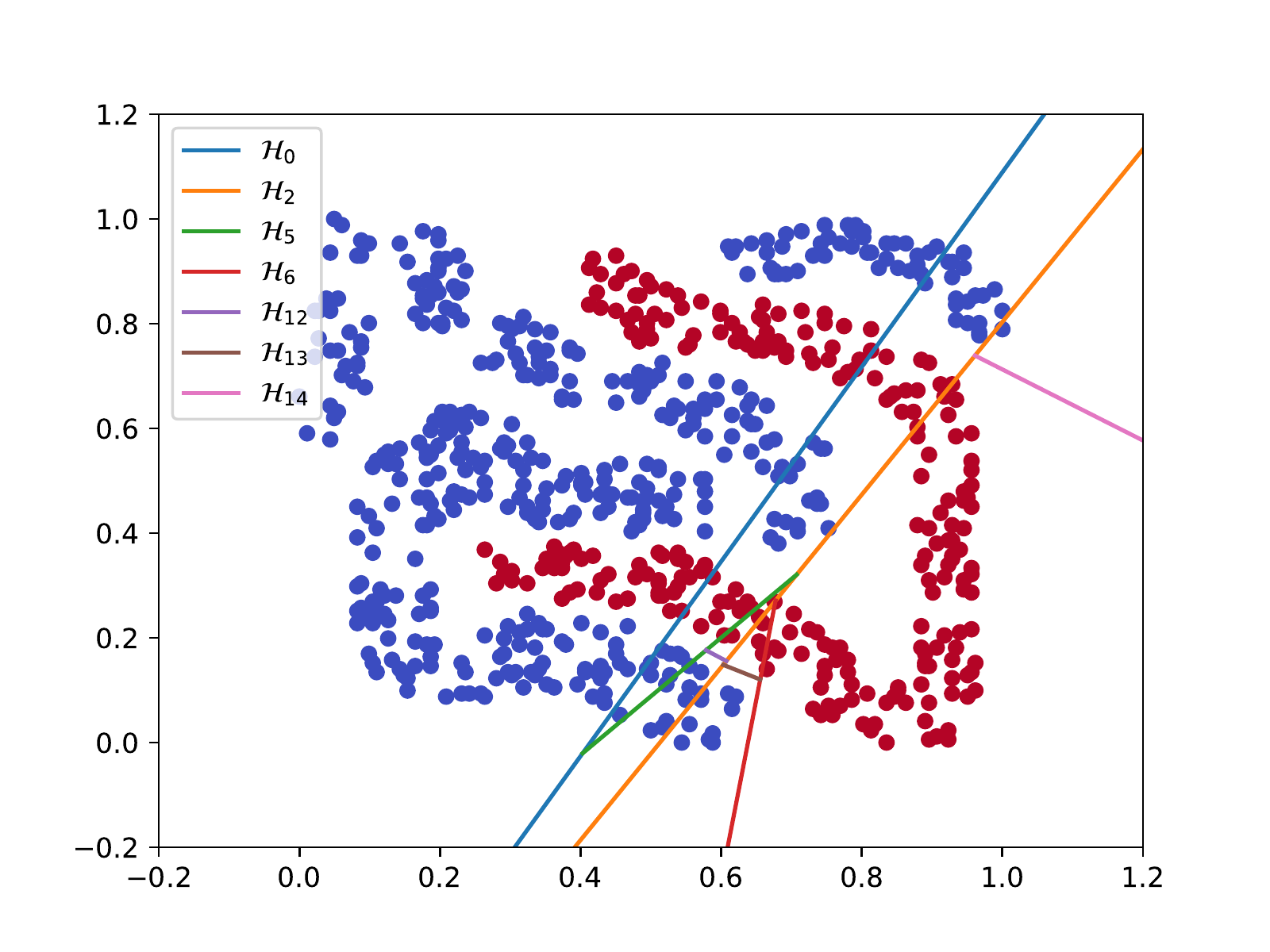}  & \qquad &
\includegraphics[width=.264\linewidth]{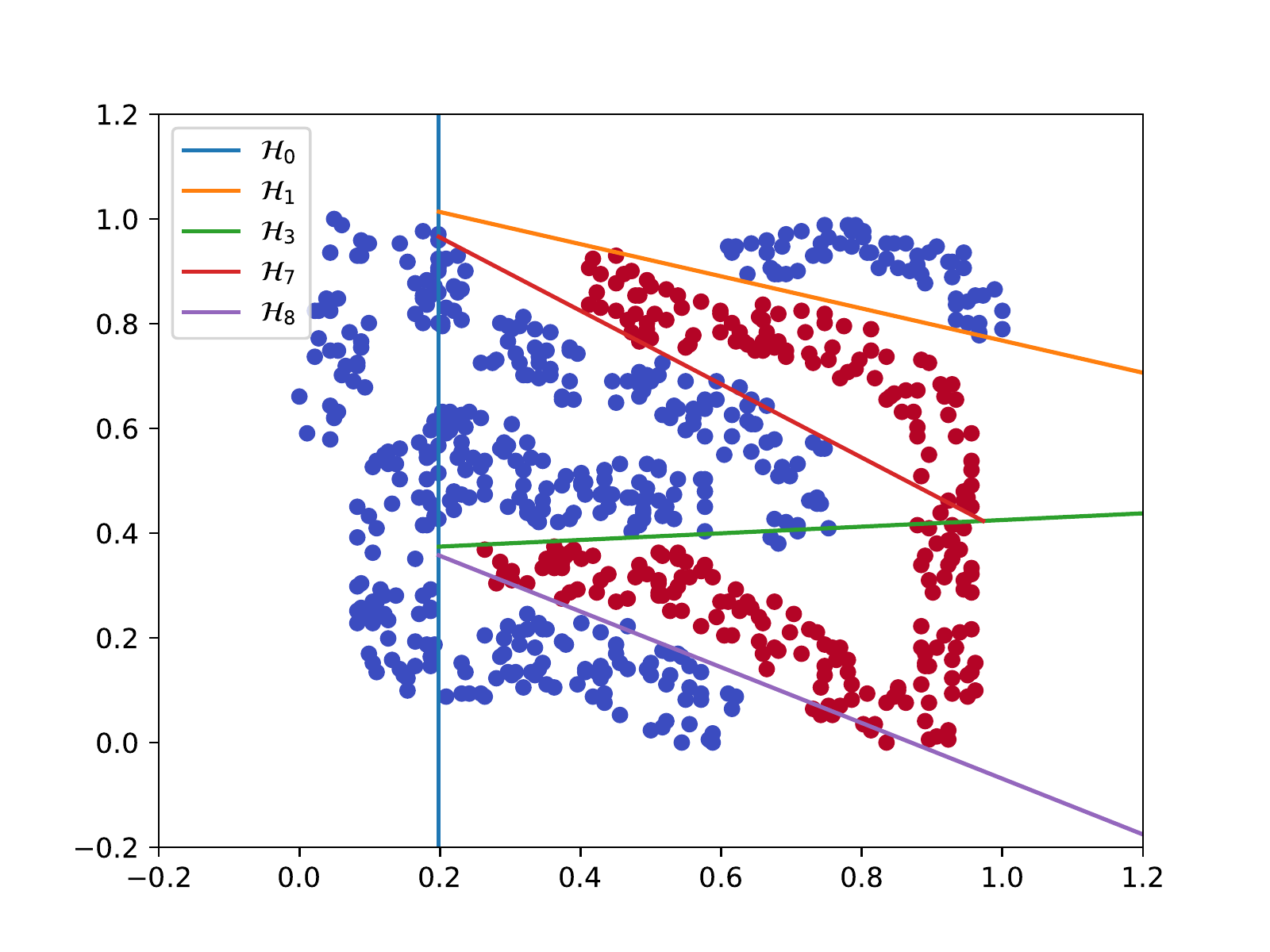}\\
        \footnotesize{(i) Local SVM Heuristic}  & & \footnotesize{(ii) OCT-H }
    \end{tabular}
        \begin{tabular}{ccc}
\includegraphics[width=.264\linewidth]{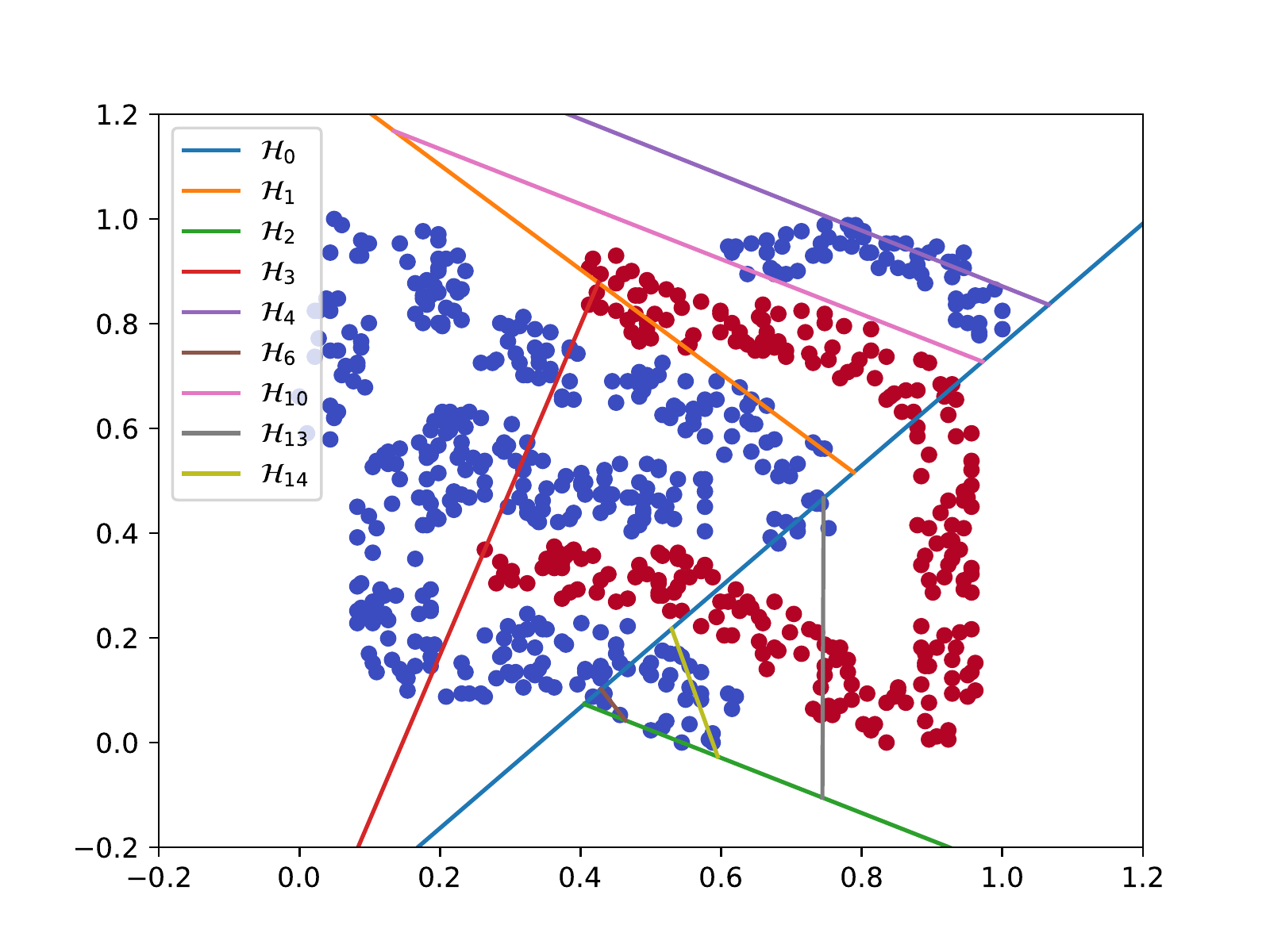} & \qquad&
\includegraphics[width=.264\linewidth]{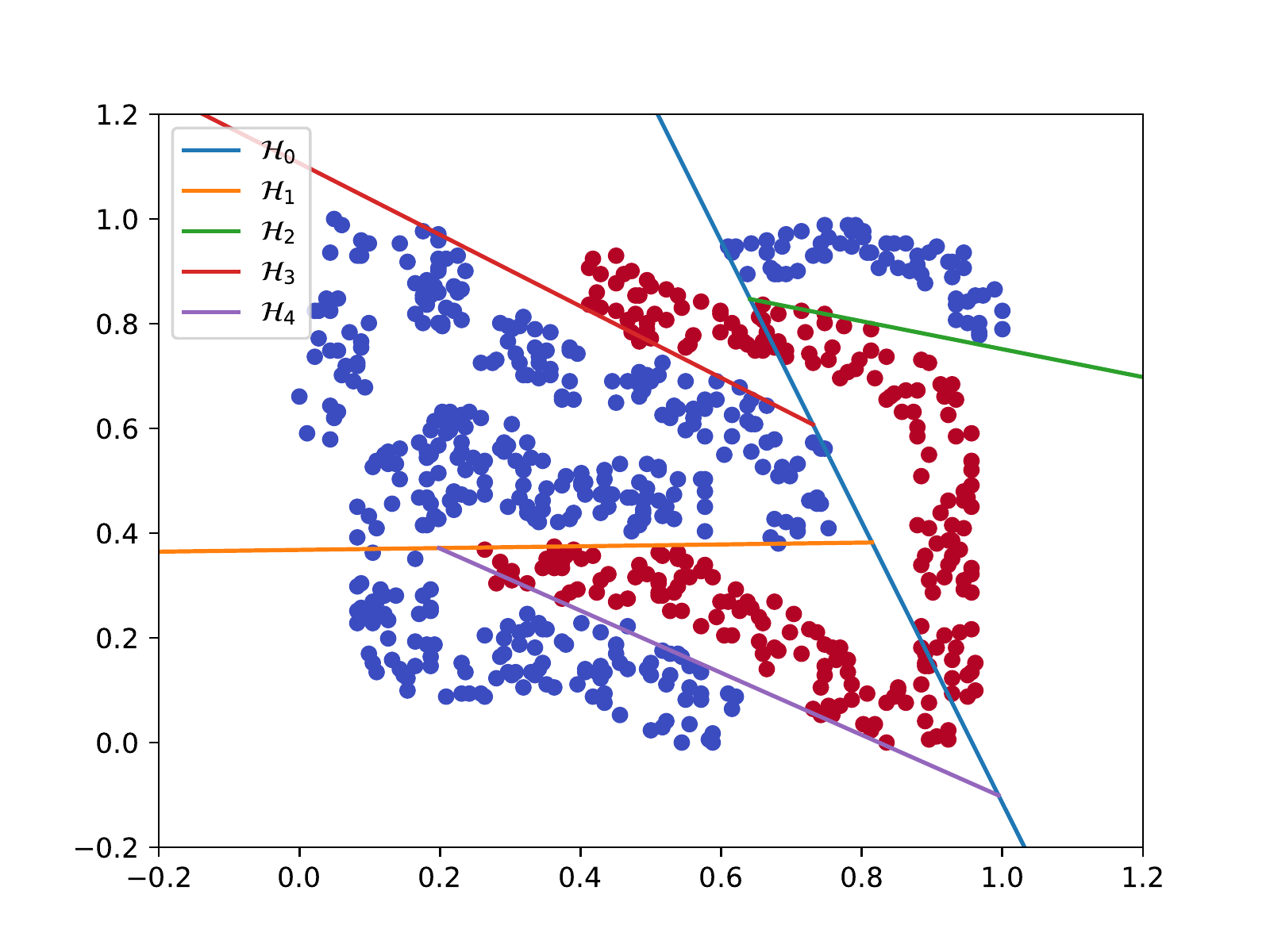}\\
    \footnotesize{(iii) MM-SVM-OCT} & & \footnotesize{(iv) MARGOT}\\
    \end{tabular}
    \caption{Results on the \texttt{fourclass} synthetic dataset.}
    \label{fig: fourclass complete}
\end{figure}

It is worth noticing that all the approaches on the \texttt{6-partitions} and \texttt{fourclass} datasets tend to minimize the complexity of the tree, i.e. the number of hyperplane splits created. In OCT-H model, this is a consequence of both the penalization of the selection of features in the objective function and the presence of specific constraints and variables which model the pruning of the tree. Similarly, MM-SVM-OCT controls the complexity of the tree with a penalization term in the objective function by introducing binary variables and related constraints. On the contrary, in MARGOT model the complexity is implicitly minimized in the objective function. Indeed, creating an hyperplane split at a node $t$ leads to new coefficients $w_t\neq 0$ and variables $\xi_t \neq 0$ 
that appear in the objective function. 



\begin{table}[ht]
\centering
\scalebox{.8}{
\renewcommand\arraystretch{1.2}
\begin{tabular}{=lcccccccccccccccccc}
                                       &     &     &  & \multicolumn{3}{c}{Local SVM} &  & \multicolumn{3}{c}{OCT-H} &  & \multicolumn{3}{c}{MM-SVM-OCT} &  & \multicolumn{3}{c}{MARGOT} \\ \hline\hline
Dataset                                & $|\I|$ & $D$ &  & Time     & Gap     & ACC      &  & Time    & Gap     & ACC   &  & Time      & Gap     & ACC      &  & Time     & Gap    & ACC    \\ \hline
\texttt{4-partitions} & 108 & 2   &  & \textbf{0.1}      & -       & 86.1     &  & 53.5    & \textbf{0}       & \textbf{100}   &  & 5.2       &\textbf{ 0}       & \textbf{100}      &  & 0.2      & \textbf{0}      & \textbf{100}    \\
\texttt{6-partitions} & 96  & 3   &  &\textbf{ 0.1  }    & -       & 91.7     &  &  \underline{14400}   & 50      & \textbf{100}   &  & 181.2    & \textbf{0}       & 96.9     &  & {4215.0}   & \textbf{0}      & \textbf{100}    \\
\texttt{fourclass}    & 689 & 4   &  & \red{\textbf{126.0}}      &     -    & \red{80.0}     &  & \underline{14400}   & {100.0}   & 97.8  &  & \underline{14400}     &    74.5     &   88.4       &  & \underline{14400}    & \textbf{67.2}   & \textbf{99.6 } 
\end{tabular}}
\caption{ Time (s), Gap (\%) and ACC (\%) performances on the synthetic datasets (time limit $= 14400$s $=4$h).}\label{tab:2D}
\end{table}

\subsection{Results on UCI datasets}\label{sec: results_UCI}

We compare the different OCT models computing two predictive measures: the accuracy (ACC) and the balanced accuracy (BACC). The ACC is the percentage of correctly classified samples, and the BACC is the mean of the percentage of correctly classified samples with positive labels and the percentage of correctly classified samples with negative labels.

Hence,
$$ACC = \frac{TP + TN}{TP +TN +FP +FN},$$

$$BACC = \frac{    
\frac{TP}{TP + FN}  + \frac{TN}{TN+FP}
}{2},$$
where TP are true positives, TN are true negatives, FP false positives and FN false negatives. 
Information about the datasets considered is provided in Table \ref{table: datasets}.
We first partitioned each dataset in training ($80\%$) and test ($20\%$) sets 
and then performed a 4-fold cross-validation (4-FCV) on the training set in order to find the best hyperparameters. 
In the first set of results where feature selection is not taken into account, the selected hyperparameters are the ones which gave the best average validation accuracy in the \mnotee{cross-validation} process. 
\mnotee{For the second set of results, we implemented a more specific tuning of the hyperparameters to address the sparsity of the hyperplanes' coefficients, as will be explained later in this section.}
Once the best hyperparameters are selected, we compute the predictive measures on the training and test sets. The results on the training dataset, as well as all the hyperparameters used in this paper, can be found in the \mnotee{\ref{app: appendix2}}.

\begin{table}[ht]
\centering
\renewcommand\arraystretch{1.2}
\begin{tabular}{lcccc}
Dataset                  &  & $|\I|$   & $n$ & Class (\%) \\ \hline \hline
Breast Cancer D. &  & 569 & 30  & 63/36      \\
Breast Cancer W.  &  & \red{683} & 9   & \red{65/35}      \\
Climate Model    &  & 540 & 18  & 9/91       \\
Heart Disease C.  &  & 297 & 13  & 54/46      \\
Ionosphere               &  & \red{351} & \red{34}  & 36/64      \\
Parkinsons               &  & 195 & 22  & 25/75      \\
Sonar                    &  & 208 & 60  & 53/47      \\
SPECTF H.                   &  & 267 & 44  & 21/79      \\
Tic-Tac-Toe              &  & 958 & \red{27}  & 35/65      \\
Wholesale                &  & 440 & 7   & 81/19   \\\hline  
\end{tabular}
\caption{Information about the datasets considered.}
\label{table: datasets}
\end{table}

For the resolution of MARGOT, HFS-MARGOT and SFS-MARGOT models, we injected warm start solutions produced by the Local SVM Heuristic. 
\red{Similarly, for OCT-1, we used as a heuristic solution the one produced by CART \cite{Breiman1984CART} using the same setting of OCT-1 for depth $D$, the minimum number of samples per leaf $N_{min}$ and complexity parameter $\alpha$.}
OCT-H was also given its starting solution, following the warm start procedure presented in \cite{Bertsimas2017OptimalClassification}.
A time limit of 30 seconds was set for every warm start procedure, and an overall time limit of 600 seconds was set for training the models.
The maximum depth of the trees generated was fixed to $D=2$, and the ranges of the hyperparameters used in the 4-FCV for the different models are the following:

\begin{itemize}
    \item For \red{OCT-1 and} OCT-H, $N_{min}$ was set to 5\% of the total number of training samples and the grid used for the hyperparameter $\alpha$ is $\{0\}\cup \{2^i: i  \in \{-8, \dots,  2\}\}$.
    \item For MM-SVM-OCT, we used the same grid as the one specified in the related paper; $c_1 \in \{10^i: i \in \{-5,\dots,5\}\}$ and the complexity hyperparameter $c_3 \in \{10^i: i \in \{-2,\dots,2\}\}$.
    \item For MARGOT, we consider all possible combinations resulting
    from $C_t \in \{10^i: i\in\{ -5,\dots, 5\}\}$ for all $t\in \Tb$, and $C_1 = C_2$, imposing the same $C_t$ values for all nodes $t$ belonging to the same branching level. 
    \item For HFS-MARGOT, we used the same grid for the $C_t$ values as in MARGOT and, concerning the budget parameters $\{B_t, \text{ for all } t \in \Tb\}$, 
   \red{we varied all possible combinations resulting from values of $B_t  \in \{1,2,3\}$, with $t \in\Tb$, considering only combinations where $B_0 \leq B_1 = B_2$, with the value 3 regarded as the maximum number of features a node can admit in order to be interpretable.}
    
    \item For SFS-MARGOT, the grid used for the hyperparameter $\alpha$ is $\{2^i: i \in \{0, \dots, 10\}\}$ and we varied $C_t$  values as in MARGOT but in a smaller grid $\{10^i: i \in \{-4,-2,0,2,4\}\}$. 
    We set all budget values $B_t = 1$ for all $t\in\Tb$, allowing the model to have full flexibility on where to use more features than the budget value.

\end{itemize}


\subsubsection{Choice of the Big-M and $\varepsilon$}


Regarding the $\varepsilon$ parameter used in our formulation, for similar reasons as the one stated in \cite{Bertsimas2017OptimalClassification}, we set $\varepsilon = 0.001$ as a compromise between choosing small values that lead to numerical issues and large values that can affect the feasible region excluding possible solutions. 
Moreover, we carefully tuned the big-M parameters through extensive computational experiments in order to find values as tight as possible. As a result, we set those values as follows: 
$$M_\xi = M_w = 50,\ M_\H = 100,$$
while the Big-M values in MM-SVM-OCT are fixed as indicated in  \cite{BlancoRobust}.

\subsubsection{First set of results}

The first set of results is shown in Table \ref{table: test_acc_BACC}, where we compare the predictive performances of MARGOT against OCT-H and MM-SVM-OCT. We can see how MARGOT takes full advantage of the generalization capabilities deriving from the maximum margin approach, resulting in much higher ACC and BACC scores on the test sets. 

\begin{table}[ht!]
\centering
\renewcommand\arraystretch{1.2}
\begin{tabular}{=l +c +c +c +c +c +c +c +c +c}
                                  &  & \multicolumn{2}{c}{OCT-H} &  & \multicolumn{2}{c}{MM-SVM-OCT} &  & \multicolumn{2}{c}{MARGOT}    \\ \hline\hline
Dataset                           &  & ACC     & BACC            &  & ACC            & BACC          &  & ACC           & BACC          \\ \hline
Breast Cancer D.                  &  & 94.7    & 94.3            &  & 93.9           & 92.7          &  & \textbf{97.4} & \textbf{96.9} \\
\rowstyle{\color{red}} Breast Cancer W.                  &  & 94.9    & 95.6            &  & \textbf{96.4}  & \textbf{96.2} &  & \textbf{96.4} & \textbf{96.2} \\
Climate Model                     &  & 93.5    & 86.4            &  & \textbf{97.2}  & \textbf{88.4} &  & \textbf{97.2} & \textbf{88.4} \\
Heart Disease C.                  &  & 80.0    & 79.7            &  & 81.7           & 81.3          &  & \textbf{83.3} & \textbf{83.0} \\
\rowstyle{\color{red}} Ionosphere                        &  & 87.3    & 85.7            &  & 85.9           & 80.9          &  & \textbf{93.0} & \textbf{90.0} \\
Parkinsons                        &  & 82.1    & \textbf{84.7}   &  & \textbf{87.2}  & 81.6          &  & 84.6          & 83.1          \\
Sonar                             &  & 66.7    & 65.9            &  & \textbf{73.8}  & \textbf{73.0} &  & \textbf{73.8}          & \textbf{73.0}          \\
SPECTF H.                         &  & 74.1    & 53.3            &  & \textbf{79.6}  & 50.0          &  & \textbf{79.6} & \textbf{56.8} \\
\rowstyle{\color{red}}Tic-Tac-Toe &  & 96.9    & 96.2            &  & 97.4           & \textbf{97.0}          &  & \textbf{97.9} & \textbf{97.0} \\
Wholesale                         &  & 84.1    & 79.8            &  & 83.0           & 74.2          &  & \textbf{87.5} & \textbf{85.1} \\ \hline
\end{tabular}
\caption{Results on the \red{test} predictive performances of the OCT models evaluated: test ACC (\%) and test BACC (\%).}
\label{table: test_acc_BACC}
\end{table}

Computational times and MIP gaps can be found in Table \ref{table: time_gap}. It is clear how MARGOT optimization problem is much easier to solve than OCT-H and MM-SVM-OCT. Indeed, 9 times out of 10, MARGOT reaches the optimal solution with a mean computational time and MIP gap of \red{121.5} seconds and 0.3\%, respectively. MM-SVM-OCT reaches the optimal solution 3 times out of 10,  with a mean running time and gap of \red{482.7} seconds and \red{18.8}\% respectively, and OCT-H reaches the optimum \fnotee{only in one case}, with a mean time and gap of \red{564.3} seconds and \red{63.7}\% respectively.

\begin{table}[ht]
\centering
\renewcommand\arraystretch{1.2}
\begin{tabular}{=l+c+c+c+c+c+c+c+c+c}
                                        &  & \multicolumn{2}{c}{OCT-H}  &  & \multicolumn{2}{c}{MM-SVM-OCT}              &  & \multicolumn{2}{c}{MARGOT}                \\ \hline\hline
Dataset                                 &  & Time        & Gap          &  & Time                & Gap                   &  & Time              & Gap                   \\ \hline
Breast Cancer D.                        &  & \underline{620.2} & 72.8         &  & \underline{600.2}   & 27.5                  &  & \textbf{7.4}      & \textbf{0.0}          \\
\rowstyle{\color{red}} Breast Cancer W. &  & \underline{615.3} & 80.1         &  & 333.0               & \textbf{\textbf{0.0}} &  & \textbf{8.7}      & \textbf{\textbf{0.0}} \\
Climate Model                           &  & \underline{620.2} & 90.6         &  & \underline{600.2}   & 0.6                   &  & \textbf{10.7}     & \textbf{0.0}          \\
Heart Disease C.                        &  & \underline{630.1} & 91.9         &  & \underline{600.0}   & 28.8                  &  & \textbf{318.1}    & \textbf{0.0}          \\
\rowstyle{\color{red}} Ionosphere       &  & 59.1        & \textbf{0.0} &  & {\underline{600.0}} & 8.4                   &  & \textbf{11.2}     & \textbf{\textbf{0.0}} \\
Parkinsons                              &  & \underline{612.3} & 90.8         &  & \underline{600.1}   & 12.5                  &  & \textbf{207.4}    & \textbf{0.0}          \\
Sonar                                   &  & \underline{620.2} & 96.1         &  & 262.7               & \textbf{0.0}          &  & \textbf{1.5}      & \textbf{0.0}          \\
SPECTF H.                               &  & \underline{620.2} & 96.9         &  & \textbf{30.3}       & \textbf{0.0}          &  & \underline{600.2} & 3.3                   \\
\rowstyle{\color{red}}Tic-Tac-Toe       &  & \underline{625.4} & 100.0        &  & \underline{600.1}   & 78.1                  &  & \textbf{3.5}      & \textbf{0.0}          \\
Wholesale                               &  & \underline{620.1} & 95.6         &  & \underline{600.0}   & 32.4                  &  & \textbf{46.4}     & \textbf{0.0}          \\ \hline
\end{tabular}
\caption{Results on the optimization performances of the OCT models evaluated: computational times (s) and MIP Gaps (\%).}
\label{table: time_gap}
\end{table}

\medskip

\subsubsection{Second set of results: feature selection}

In the following set of results, we compare HFS-MARGOT and SFS-MARGOT with \red{OCT-1 and} OCT-H. Both MARGOT and MM-SVM-OCT do not appear in this set of results because these models do not address \mnotee{the} sparsity of the \mnotee{hyperplanes' weights}. For this analysis, a different hyperparameter selection was carried out. This was done \fnotee{to take into account} that we are not just comparing the predictive performances of the methods, but we are also evaluating the feature selection aspect. The 4-FCV was \mnotee{still} conducted, but this time we did not select the hyperparameters which gave the highest mean validation accuracy. Indeed, \mnotee{the} highest mean validation accuracy values yield to models which select a high number of features, thus contrasting the aim to create more interpretable trees. At the same time, \mnotee{solely considering sparsity is not useful as} the hyperparameters which gave the best results in terms of feature selection may result in less performing classifiers.

\begin{table}[ht!]
\centering
\renewcommand\arraystretch{1.2}
\begin{tabular}{=l+c+c+c+c+c+c+c+c+c+c+c+c}
 &  & \multicolumn{2}{c}{\red{OCT-1}} &  & \multicolumn{2}{c}{OCT-H*} &  & \multicolumn{2}{c}{HFS-MARGOT*} &  & \multicolumn{2}{c}{SFS-MARGOT*} \\ \hline\hline
Dataset &  & ACC & BACC &  & ACC & BACC &  & ACC & BACC &  & ACC & BACC \\ \hline
Breast Cancer D. &  & 91.2 & 91.6 &  & 94.7 & 94.3 &  & \textbf{95.6} & \textbf{95.5} &  & 94.7 & 94.3 \\
\rowstyle{\color{red}} Breast Cancer W. &  & 92.0 & 90.9 &  & 92.0 & 91.4 &  & \textbf{94.2} & \textbf{94.5} &  & \textbf{94.2} & 93.6 \\
Climate Model &  & 91.7 & 60.1 &  & \textbf{98.1} & \textbf{93.9} &  & \red{96.3} & \red{77.8} &  & 96.3 & 82.8 \\
Heart Disease C. &  & 71.7 & 71.0 &  & 80.0 & 79.5 &  & 83.3 & 82.6 &  & \textbf{86.7} & \textbf{86.2} \\
\rowstyle{\color{red}} Ionosphere &  & \textbf{91.5} & \textbf{91.7} &  & 90.1 & 86.9 & \textbf{} & 87.3 & 83.8 &  & 84.5 & 79.8 \\
Parkinsons &  & \textbf{89.7} & {83.3} &  & {87.2} & 81.6 &  & \red{\textbf{89.7}} & \red{{83.3}} &  & {87.2} & \textbf{84.8} \\
Sonar &  & 69.0 & 68.9 &  & 71.4 & 71.1 &  & 71.4 & 70.9 &  & \textbf{73.8} & \textbf{73.6} \\
SPECTF H. &  & 77.8 & 65.8 &  & 75.9 & 51.1 &  & 75.9 & 57.8 &  & \textbf{83.3} & \textbf{65.9} \\
\rowstyle{\color{red}} Tic-Tac-Toe &  & 69.3 & 61.5 &  & 96.4 & 95.5 &  & 76.0 & 65.7 &  & \textbf{97.9} & \textbf{97.0} \\
Wholesale &  & 86.4 & 85.2 &  & 87.5 & 86.1 &  & 86.4 & 85.2 &  & \textbf{88.6} & \textbf{86.9} \\ \hline
\end{tabular}
\caption{Results on the \red{test} predictive performances of the \red{OCT models with feature selection}: test ACC (\%) and test BACC (\%).}
\label{table: test_acc_BACC_fs}
\end{table} 

\begin{table}[ht!]
\centering
\scalebox{1.0}{
\renewcommand\arraystretch{1.2}
\begin{tabular}{=l+c+c+c+c+c+c+c+c+c+c+c+c}
                                        &  & \multicolumn{2}{c}{OCT-1}           &           & \multicolumn{2}{c}{OCT-H*} &  & \multicolumn{2}{c}{HFS-MARGOT*}              &  & \multicolumn{2}{c}{SFS-MARGOT*}   \\ \hline\hline
Dataset                                 &  & Time                & Gap           &           & Time            & Gap      &  & Time                    & Gap                &  & Time              & Gap           \\ \hline
Breast Cancer D.                        &  & \underline{600.0}   & 96.3          &           & \underline{620.2}     & 72.8     &  & \textbf{313.9}          & \textbf{0.0}       &  & 314.7             & \textbf{0.0}  \\
\rowstyle{\color{red}} Breast Cancer W. &  & 138.6               & \textbf{0.0}  &           & \underline{615.4}     & 72.3     &  & \textbf{14.1}           & \textbf{0.0}       &  & 37.6              & \textbf{0.0}  \\
Climate Model                           &  & \underline{600.0}   & 100.0         &           & \underline{620.1}     & 81.4     &  & \red{\underline{600.3}} & \red{100.0}        &  & \underline{600.4} & \textbf{14.2} \\
Heart Disease C.                        &  & \textbf{93.3}       & \textbf{0.0}  & \textbf{} & \underline{630.1}     & 86.9     &  & 106.5                   & \textbf{0.0}       &  & \underline{600.3} & 28.2          \\
\rowstyle{\color{red}} Ionosphere       &  & {\underline{600.0}} & 52.0          & \textbf{} & \underline{625.2}     & 83.2     &  & 124.0                   & \textbf{0.0}       &  & \textbf{15.9}     & \textbf{0.0}  \\
Parkinsons                              &  & \underline{600.0}   & 92.9          &           & \underline{620.1}     & 82.4     &  & \red{\textbf{7.8}}      & \red{\textbf{0.0}} &  & \underline{600.2} & 62.0          \\
Sonar                                   &  & \underline{600.1}   & \textbf{30.7} &           & \underline{620.1}     & 91.4     &  & {\underline{601.2}}     & 98.3               &  & \underline{610.8} & 42.2          \\
SPECTF H.                               &  & \underline{600.0}   & 98.0          &           & \underline{626.1}     & 94.5     &  & {\underline{601.3}}     & \textbf{93.8}      &  & \underline{610.1} & 99.9          \\
\rowstyle{\color{red}}Tic-Tac-Toe       &  & 268.1               & \textbf{0.0}  &           & \underline{630.0}     & 94.7     &  & \underline{604.6}       & 88.2               &  & \textbf{47.2}     & \textbf{0.0}  \\
Wholesale                               &  & \underline{600.0}   & 86.5          &           & \underline{620.2}     & 73.5     &  & \textbf{36.4}           & \textbf{0.0}       &  & \underline{600.2} & 80.1          \\ \hline
\end{tabular}}
\caption{Results on the optimization performances of the \red{OCT models with feature selection}: computational times (s) and MIP Gap values (\%). }
\label{table: time_gap_fs}
\end{table}

Thus, we proceeded as follows. For each dataset, after performing the standard 4-FCV, we highlighted the combinations of hyperparameters which resulted in a mean validation accuracy in the range $[0.975\gamma, \gamma]$, where $\gamma$ is the best mean validation accuracy value scored. This way, we selected the combinations of hyperparameters corresponding to "good" classifiers. Among these combinations, we chose the ones corresponding to the \mnotee{lower} number of \mnotee{features} used, and, among these last ones, we picked the combination corresponding to the best validation accuracy. 
We denote by OCT-H*, HFS-MARGOT* and SFS-MARGOT* the tree models generated with this \mnotee{feature selection driven} hyperparameter tuning. \red{For OCT-1, we performed the standard hyperparameter search considering that the univariate splits of the model are sparse by definition.} To the best of our knowledge, this is the first time a tailored \mnotee{cross-validation} was carried out
to fairly compare optimization-based ML models that embed feature selection.

{
Tables \ref{table: test_acc_BACC_fs} and \ref{table: time_gap_fs} report both the predictive and optimization performances of the compared models. \red{Tables \ref{table: features} and \ref{table: features_fs}},
together with \red{their} graphical representation in Figure \ref{fig: features_oct}, show the difference in the features selected among the
\mnotee{analysed} OCT models. We denote by $F$ the set of distinct features used overall in the tree, and by $F_t$
the set of features selected at node $t$. 
As expected, MARGOT and MM-SVM-OCT
models tend to use all the available features because they have no feature selection constraints or
penalization terms in the objective function. One thing to notice is that, in many cases, MM-SVM-OCT tends to
activate more branch splits than MARGOT, each involving all the features. This might be due to the difference in the objective
functions of the two models.}
Moreover, the OCT-H model for which a standard 4-FCV is carried out, \fnote{tends to produce sparser tree models since the number of features selected is penalized in its formulation}. 


OCT-1, OCT-H*, HFS-MARGOT* and SFS-MARGOT* generate models selecting a lower number of features, maintaining good prediction performances, as shown in Table \ref{table: test_acc_BACC_fs}. In particular, we can notice how SFS-MARGOT* presents better ACC and BACC values, and this is probably due to the combination of the maximum margin approach and the feature selection penalization. 
Indeed, we can see how SFS-MARGOT*, where violation of the budget constraints is allowed, has more freedom in the selection of features compared to the more restrictive approach of HFS-MARGOT*,  which at each node cannot exceed a predefined number of selected features.
One thing to notice is that, on these results, both OCT-H and OCT-H* tend to use the selected features just in the first branch node of the tree. \fnotee{A similar behaviour is exhibited} by the SFS-MARGOT$^*$ model. \fnotee{In contrast,} the features selected by HFS-MARGOT$^*$ tend to be more spread throughout the branch nodes. Using hard budget constraints limiting the number of features selected at each branch node has indeed two consequences: on the one hand, it spreads more evenly the features selected among the tree structure, while on the other, this restriction  
might result insufficient in order achieve the best performances. This is the case of the Tic-Tac-Toe dataset where the HFS-MARGOT* model clearly did not perform well, and this is most likely because the values we adopted for budget parameters $B_t$ were too limiting. \red{In this sense, OCT-1 represents the extreme case, in that it is limited to select only one feature per split, generally worsening the predictive quality of the classifier.} Of course, for HFS-MARGOT*, we could have used higher budget values, but this, apart from leading to less interpretable tree structures and time-consuming hyperparameters tuning, does not attain the scope of our computational results. \red{We finally note how OCT-1 is easier to solve than OCT-H*, closing the gap in 3 datasets, but it still results to be computationally more expensive than the MARGOT feature selection models.}


In Figures \ref{fig:tree_graphs}, we give more insight into how the features selected by the different OCT models divide the training samples. We focused on two datasets, the Parkinsons and the \red{Breast Cancer D.} ones, in that we found them explicative of the behaviour of all the trees \mnotee{generated} by \red{OCT-1}, OCT-H*, HFS-MARGOT* and SFS-MARGOT* models.

\begin{table}[ht!]
\centering
\scalebox{0.89}{
\renewcommand\arraystretch{1.2}
\begin{tabular}{=l+c+c+c+c+c+c+c+c+c+c}
                                        &     &  & \multicolumn{2}{c}{OCT-H}   &  & \multicolumn{2}{c}{MM-SVM-OCT} &  & \multicolumn{2}{c}{MARGOT}  \\ \hline\hline
Dataset                                 & $n$ &  & $|F|$ & $|F_0|,|F_1|,|F_2|$ &  & $|F|$   & $|F_0|,|F_1|,|F_2|$  &  & $|F|$ & $|F_0|,|F_1|,|F_2|$ \\ \hline
Breast Cancer D.                        & 30  &  & 3     & 3, 0, 0             &  & 30      & 30, 30, 30           &  & 30    & 30, 30, 0           \\
\rowstyle{\color{red}} Breast Cancer W. & 9   &  & 4     & 4, 0, 0             &  & 9       & 9, 9, 0              &  & 9     & 9, 8, 0             \\
Climate Model                           & 18  &  & 11    & 11, 0, 0            &  & 18      & 18, 18, 18           &  & 18    & 18, 0, 18           \\
Heart Disease C.                        & 13  &  & 4     & 4, 0, 0             &  & 13      & 13, 0, 0             &  & 13    & 13, 0, 6            \\
\rowstyle{\color{red}} Ionosphere       & 34  &  & 32    & 28, 22, 23          &  & 33      & 33, 33, 33           &  & 33    & 33, 33, 31          \\
Parkinsons                              & 22  &  & 13    & 7, 0, 8             &  & 22      & 22, 22, 22           &  & 22    & 22, 22, 22          \\
Sonar                                   & 60  &  & 23    & 23, 0, 0            &  & 60      & 60, 60, 60           &  & 51    & 51, 0, 0            \\
SPECTF H.                               & 44  &  & 25    & 10, 4, 17           &  & 0       & 0, 0, 0              &  & 44    & 44, 0, 42           \\
\rowstyle{\color{red}} Tic-Tac-Toe      & 27  &  & 26    & 18, 18, 1           &  & 27      & 26, 26, 22           &  & 18    & 18, 0, 0            \\
Wholesale                               & 7   &  & 7     & 6, 0, 5             &  & 7       & 7, 7, 7           &  & 7     & 7, 0, 6             \\ \hline
\end{tabular}
}
\caption{Comparison on the number of features selected by the OCT models; $F$ is is the set of features used in the tree, and $F_t$ are the features selected at the node $t=0,1,2$.}
\label{table: features}
\end{table}

\begin{table}[ht!]
\hspace{-0.25cm}
\scalebox{0.89}{
\renewcommand\arraystretch{1.2}
\begin{tabular}{=l+c+c+c+c+c+c+c+c+c+c+c+c+c}
                                        &     &  & \multicolumn{2}{c}{OCT-1}   &  & \multicolumn{2}{c}{OCT-H*}  &  & \multicolumn{2}{c}{HFS-MARGOT*} &  & \multicolumn{2}{c}{SFS-MARGOT*} \\ \hline\hline
Dataset                                 & $n$ &  & $|F|$ & $|F_0|,|F_1|,|F_2|$ &  & $|F|$ & $|F_0|,|F_1|,|F_2|$ &  & $|F|$    & $|F_0|,|F_1|,|F_2|$  &  & $|F|$   & $|F_0|,|F_1|,|F_2|$   \\ \hline
Breast Cancer D.                        & 30  &  & 2     & 1, 0, 1             &  & 3     & 3, 0, 0             &  & 4        & 2, 0, 2              &  & 3       & 3, 0, 0               \\
\rowstyle{\color{red}} Breast Cancer W. & 9   &  & 3     & 1, 1, 1             &  & 2     & 2, 0, 0             &  & 4        & 2, 3, 0              &  & 3       & 3, 0, 0               \\
Climate Model                           & 18  &  & 2     & 1, 1, 1             &  & 7     & 7, 0, 0             &  & \red{6}  & \red{3, 3, 0}        &  & 4       & 4, 0, 0               \\
Heart Disease C.                        & 13  &  & 1     & 1, 0, 0             &  & 3     & 3, 0, 0             &  & 3        & 1, 2, 2              &  & 6       & 6, 0, 0               \\
\rowstyle{\color{red}} Ionosphere       & 34  &  & 2     & 1, 0, 1             &  & 3     & 2, 0, 1             &  & 7        & 2, 2, 3              &  & 2       & 1, 0, 1               \\
Parkinsons                              & 22  &  & 2     & 1, 1, 1             &  & 5     & 3, 3, 0          &  & \red{3}  & \red{1, 2, 0}        &  & 12      & 9, 2, 4               \\
Sonar                                   & 60  &  & 1     & 1, 0, 0             &  & 15    & 15, 0, 0            &  & 5        & 1, 2, 2              &  & 8       & 7, 0, 1               \\
SPECTF H.                               & 44  &  & 2     & 1, 0, 1             &  & 5     & 5, 0, 0             &  & 6        & 2, 1, 3              &  & 10      & 9, 1, 0               \\
\rowstyle{\color{red}} Tic-Tac-Toe      & 27  &  & 2     & 1, 0, 1             &  & 18    & 18, 0, 0            &  & 5        & 2, 3, 0              &  & 18      & 18, 0, 0              \\
Wholesale                               & 7   &  & 2     & 1, 1, 1             &  & 2     & 2, 0, 0             &  & 5        & 1, 2, 2              &  & 3       & 3, 0, 0               \\ \hline
\end{tabular}
}
\caption{Comparison on the number of features selected by the \red{OCT models with feature selection}; $F$ is is the set of features used in the tree, and $F_t$ are the features selected at the node $t=0,1,2$.}
\label{table: features_fs}
\end{table}

\begin{sidewaysfigure}[ph!]     \includegraphics[width=\textwidth]{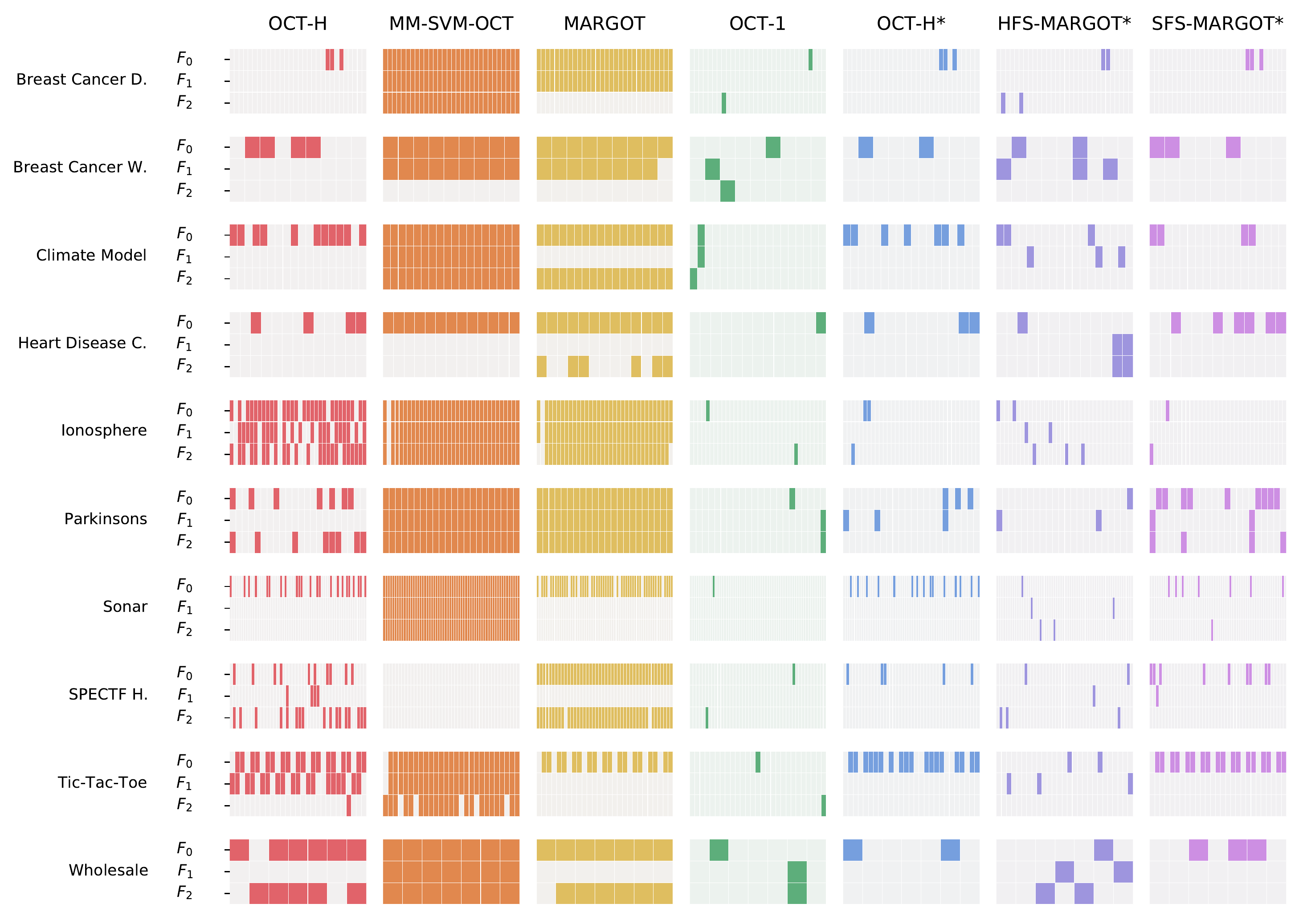}
\caption{Graphical representation of \red{Tables \ref{table: features}, \ref{table: features_fs}}: for each branch node $t$, the darker coloured cells correspond to the features selected $F_t$.}
\label{fig: features_oct}
\end{sidewaysfigure}


\begin{figure}[ht!]
\begin{subfigure}{1\textwidth}
\centering 
\begin{subfigure}[t]{0.49\textwidth}
\centering 
\scalebox{0.5}{\input{parkinsons_tree_oct_cut}}\\
\footnotesize{\red{OCT-1}}
\end{subfigure}%
\hfill
\begin{subfigure}[t]{0.49\textwidth}
\centering 
\scalebox{0.5}{\input{parkinsons_tree_octh_cut}}\\
\footnotesize{OCT-H*}
\end{subfigure}\bigskip

\begin{subfigure}[t]{0.49\textwidth}
  \centering 
\scalebox{0.5}{\input{parkinsons_tree_margot_hardb}}\\
\footnotesize{HFS-MARGOT*}
\end{subfigure}
\hfill
\begin{subfigure}[t]{0.49\textwidth}
  \centering 
\scalebox{0.5}{\input{parkinsons_tree_margot_softb}}\\
\footnotesize{SFS-MARGOT*}
\end{subfigure}
\medskip
\caption{Parkinsons}
\label{fig:tree_graphs park}

\end{subfigure}\bigskip\smallskip

\begin{subfigure}{1\textwidth}
\centering 
\footnotesize
\begin{subfigure}[t]{0.50\textwidth}
\centering 
\scalebox{0.5}{\input{bcd_tree_oct_cut}}\\
\footnotesize{\red{OCT-1}}
\end{subfigure}%
\hfill
\begin{subfigure}[t]{0.50\textwidth}
\centering 
\scalebox{0.5}{\input{bcd_tree_octh_cut}}\\
\footnotesize{OCT-H*}
\end{subfigure}\bigskip

\begin{subfigure}[t]{0.50\textwidth}
  \centering
\scalebox{0.5}{\input{bcd_tree_margot_hardb_cut}}\\
\footnotesize{HFS-MARGOT*}
\end{subfigure}%
\hfill
\begin{subfigure}[t]{0.50\textwidth}
 \centering 
\scalebox{0.5}{\input{bcd_tree_margot_softb_cut}}\\
\footnotesize{SFS-MARGOT*}
\end{subfigure}
\medskip
\caption{\red{Breast Cancer D.}}
\label{fig:tree_graphs bcw}

\end{subfigure}\\

\caption{Trees generated by models in Table \ref{table: test_acc_BACC_fs} on the Parkinsons and \red{Breast Cancer D.} datasets. For each branch node, we report the number of positive and negative training samples in the squared brackets and below the features selected at each node. For each leaf node, the number of positive and negative samples is indicated, together with the assigned class label.}
\label{fig:tree_graphs}
\end{figure}

\subsection{Warm Start}
\label{sec: warm start results}




{During the branch and bound process implemented by Gurobi, after finding an initial incumbent solution, the solver applies heuristics to improve the quality of the incumbent solution before further exploring the branch and bound tree}. In general, a warm start input solution is accepted as the first incumbent if its value is better than the initial solution found by the solver.
In Table \ref{table: margot ws}, we analyse the quality of the warm start input solution produced by the Local SVM Heuristic. To this aim, we introduce the following definitions: $f_0$ refers to the value of the first incumbent solution, and $f_1$ is the value of the best incumbent solution after the root node of the branch and bound tree has been explored. 
We report these values in two different cases. In the first one, no input solution was injected, while in the second, the solver was given the warm start computed with the Local SVM Heuristic.
From the results in the tables, we can see how generally the value of the Local SVM solution is better than the one of the first incumbent solution found by Gurobi. Only for the SPECTF H. dataset, our warm start for HFS-MARGOT model did not produce a better first solution.
Similarly, in most cases, $f_1$ values are better when the Local SVM solution is injected. Moreover\mnotee{,} we can notice how, almost every time the Local SVM input solution was given, values $f_0$ and $f_1$ are equal.

\begin{table}[ht!]
\centering
\scalebox{.9}{
\renewcommand\arraystretch{1.2}
\begin{tabular}{=l+r+r+r+r+r+r}
MARGOT                   &  & \multicolumn{2}{c}{$f_0$}                        &  & \multicolumn{2}{c}{$f_1$}                \\ \hline\hline
Dataset                  &  & No warm start              & Local SVM           &  & No warm start      & Local SVM           \\ \hline
Breast Cancer D. &  & 412.69                     & \textbf{71.33}      &  & 73.31              & \textbf{71.33}      \\
\rowstyle{\color{red}} Breast Cancer W.  &  & 41939.75                   & \textbf{5250.53}    &  & 5357.29            & \textbf{5250.53}    \\
Climate Model    &  & 2592000000398.93           & \textbf{94.14}      &  & 113.17             & \textbf{94.14}      \\
Heart Disease C.  &  & 29.56                      & \textbf{18.48}      &  & \textbf{18.48}     & \textbf{18.48}      \\
\rowstyle{\color{red}} Ionosphere               &  & 1680000000000.00           & \textbf{768.96}     &  & \textbf{726.56}    & 768.96              \\
Parkinsons               &  & 306298.02                  & \textbf{196369.04}  &  & 213242.43          & \textbf{196369.04}  \\
Sonar                    &  & 11.02                      & \textbf{11.01}      &  & \textbf{0.17}      & 0.35                \\
SPECTF H.             &  & 127800000000.00            & \textbf{17.47}      &  & \textbf{17.47}     & \textbf{17.47}      \\
\rowstyle{\color{red}} Tic-Tac-Toe              &  & 346.00                     & \textbf{84.00}      &  & \textbf{84.00}     & \textbf{84.00}      \\
Wholesale                &  & 349170.61                  & \textbf{104787.61}  &  & \textbf{104787.61} & \textbf{104787.61} \\ \hline \\

HFS-MARGOT               &  & \multicolumn{2}{c}{$f_0$}                        &  & \multicolumn{2}{c}{$f_1$}                \\ \hline\hline
Dataset                  &  & No warm start              & Local SVM           &  & No warm start      & Local SVM           \\ \hline
Breast Cancer D. &  & 569527.09                  & \textbf{78042.45}   &  & 569527.09          & \textbf{78042.45}   \\
\rowstyle{\color{red}} Breast Cancer W.  &  & 45466685.23                & \textbf{7278653.51} &  & 45466685.23        & \textbf{7278653.51} \\
\rowstyle{\color{red}} Climate Model    &  & 7400074.00                 & \textbf{6624789.91} &  & 7400073.21         & \textbf{6624789.91} \\
Heart Disease C.  &  & 107570.28                  & \textbf{98193.89}   &  & 107570.28          & \textbf{98193.89}   \\
\rowstyle{\color{red}} Ionosphere               &  & 2820.49                    & \textbf{1477.29}    &  & 2820.49            & \textbf{1477.29}    \\
\rowstyle{\color{red}} Parkinsons               &  & 119951.06                  & \textbf{88203.38}   &  & 119951.06          & \textbf{88203.38}   \\
Sonar                    &  & 1434.07                    & \textbf{990.21}     &  & 1434.07            & \textbf{990.21}     \\
SPECTF H.             &  & \textbf{0.89}              & \textbf{0.89}       &  & \textbf{0.89}      & \textbf{0.89}       \\
\rowstyle{\color{red}} Tic-Tac-Toe              &  & 58300.00                   & \textbf{46528.00}   &  & 54822.00           & \textbf{46528.00}   \\
Wholesale                &  & 7532.43                    & \textbf{7256.24}    &  & 7532.43            & \textbf{7256.24}    \\ \hline \\

SFS-MARGOT               &  & \multicolumn{2}{c}{$f_0$}                        &  & \multicolumn{2}{c}{$f_1$}                \\ \hline\hline
Dataset                  &  & No warm start              & Local SVM           &  & No warm start      & Local SVM           \\ \hline
Breast Cancer D. &  & 125004.07                  & \textbf{7382.54}    &  & 24619.24           & \textbf{7382.54}    \\
\rowstyle{\color{red}} Breast Cancer W.  &  & 810.54                     & \textbf{117.26}     &  & 206.87             & \textbf{117.26}     \\
Climate Model    &  & 61588.84                   & \textbf{14800.00}   &  & \textbf{14800.00}  & \textbf{14800.00}   \\
Heart Disease C.  &  & 414.07                     & \textbf{215.00}     &  & 299.28             & \textbf{215.00}     \\
\rowstyle{\color{red}} Ionosphere               &  & 1120000004951.79           & \textbf{247.86}     &  & 404.00             & \textbf{247.86}     \\
Parkinsons               &  & 376661.06                  & \textbf{194921.43}  &  & 344525.19          & \textbf{194921.43}  \\
Sonar                    &  & 902.33                     & \textbf{228.87}     &  & 311.75             & \textbf{228.87}     \\
SPECTF H.             &  & {85200000033077.50} &  \red{\textbf{18794.98}}            &  & \textbf{8800.01}   & \textbf{8800.01}    \\
\rowstyle{\color{red}} Tic-Tac-Toe              &  & 3064000691010.00           & \textbf{101.00}     &  & \textbf{101.00}    & \textbf{101.00}     \\
Wholesale                &  & 40023.52                   & \textbf{12177.22}   &  & 12438.92           & \textbf{12177.22} \\  \hline
\end{tabular}}
\caption{Warm start analysis for MARGOT, HFS-MARGOT and SFS-MARGOT}
\label{table: margot ws}
\end{table}

\section{Conclusions and future research}

In this paper, we propose a novel MIQP model, MARGOT, to train multivariate optimal classification trees which employ maximum margin hyperplanes by following the soft SVM paradigm. The proposed model presents fewer binary variables and constraints than other OCT methods by exploiting the SVM approach and the binary classification setting, resulting in a much more compact formulation. 
The computational experience shows that MARGOT results in a much easier model to solve compared to state-of-the-art OCT models, and, thanks to the statistical properties inherited by the SVM approach, it reaches better predictive performances.
In the case sparsity of the hyperplane splits is a \mnotee{desirable} requirement, HFS-MARGOT and SFS-MARGOT represent two valid interpretable alternatives, which model feature selection with hard budget constraints and soft penalization, respectively.
Both the feature selection versions are comparable to OCT-H and MM-SVM-OCT approaches in terms of prediction quality, though \mnotee{they are} easier to solve. On the one hand, HFS-MARGOT, 
results in a more interpretable model where the selected features are evenly spread among tree branch nodes without losing too much prediction quality. 
On the other, SFS-MARGOT presents better out-of-sample performances than HFS-MARGOT,  
though the selection of the features does not exploit the hierarchical tree structure of the classifier as much.

Plenty of future directions of this work are of interest.
Firstly, the method can be extended to deal with the multi-class case. In addition, being SVMs widely used for regression tasks, a similar version of MARGOT to learn optimal regression trees can be further addressed.
Lastly, the development of a tailored optimization algorithm for the resolution of the proposed models can be investigated to improve computational times on real-world instances. 

\red{\section*{Acknowledgements}
This research has been partially carried out in the framework of the CADUCEO project (No. F/180025/01-05/X43), supported by the Italian Ministry of Enterprises and Made in Italy. This support is gratefully acknowledged. 
Laura Palagi acknowledges financial support from Progetto di Ricerca Medio Sapienza Uniroma1 - n. RM1221816BAE8A79.
Marta Monaci acknowledges financial support from Progetto Avvio alla Ricerca Sapienza Uniroma1 - n. AR1221816C6DC246 and Federico D'Onofrio acknowledges financial support from Progetto Avvio alla Ricerca Sapienza Uniroma1 - n. AR1221816C78D963.}

\bibliographystyle{apalike}
\bibliography{biblio}  

\appendix
\section{Additional results}\label{app: appendix1}

In this section, we present additional computational results regarding the OCT models\mnotee{.} 

\subsection{Scalability of MARGOT and OCT-H models with respect to the depth}

\red{In order to \mnotee{assess} how the resolution of MARGOT and OCT-H optimization models scales with increasing depth values, we compare the two models with depths $D \in \{2,3,4\}$. The hyperparameters used are the same as the ones selected for $D = 2$ \mnotee{(reported in Table \ref{table: hyperparameters})} to evaluate how the optimization problems scale only with respect to the depth hyperparameter. \fnotee{The} same time limit of $3600$ seconds was set for all experiments, and a time limit of $120$ seconds was set for the warm start procedure.
}
\red{Results on both predictive and computational performances are presented in an aggregated form using box plots \ref{fig: box_plot_testacc_gap}. We can notice how, regarding the predictive performances, both models perform similarly for every depth value, with a slight decrease at $D = 4$. Regarding the computational performances, both optimization models become harder to solve, with OCT-H reaching higher MIP Gap values than MARGOT. We notice that, 
\mnotee{though not deducible from the boxplots}, OCT-H with depths $D \in \{3,4\}$ was able to achieve a MIP Gap \mnotee{value} of 0 only for the Ionosphere dataset, while MARGOT with $D=3$ \mnotee{certified} the optimal solution only for the Breast Cancer D. and the Sonar datasets, while it always reached the time limit for depth $D=4$. 
In general, scaling with respect to the depth is a \mnotee{challenging} aspect in MIP-based OCT models as they become more complex to solve, presenting a high computational complexity depending both on the data dimensionality and on the depth of the tree. 
In particular, as shown in Table \ref{tab: dimension}, in MARGOT, both variables and constraints grow exponentially with the depth $D$ of the tree, and the same can be stated for OCT-H (\cite{Bertsimas2017OptimalClassification}).
Nonetheless, solving MARGOT for increasing values of $D$ seems to be more tractable than solving OCT-H. 
}

\begin{figure}[ht!]
\centering
    \begin{tabular}{ccc}
    \includegraphics[width=.5\linewidth]{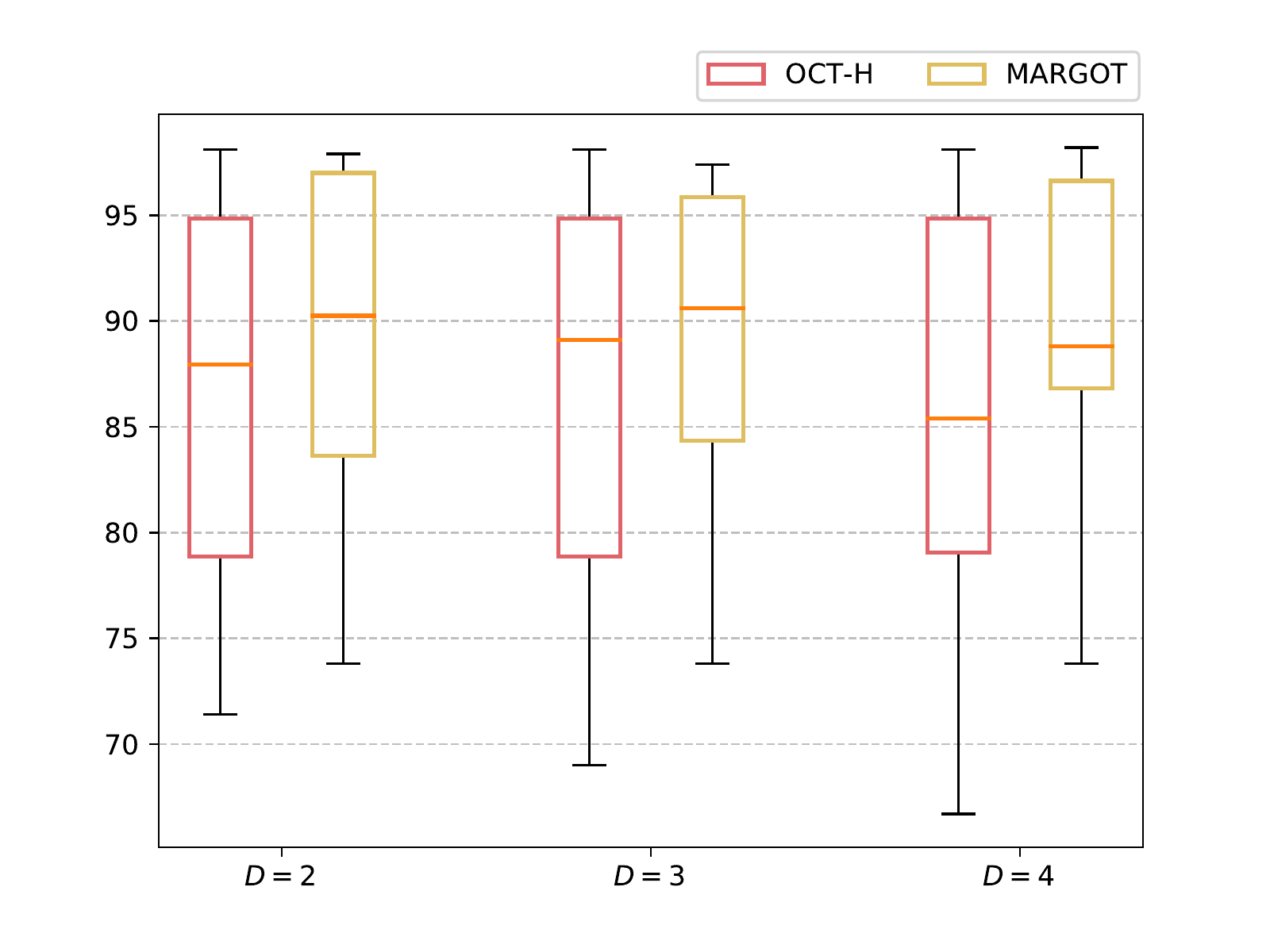}& 
    \includegraphics[width=.5\linewidth]{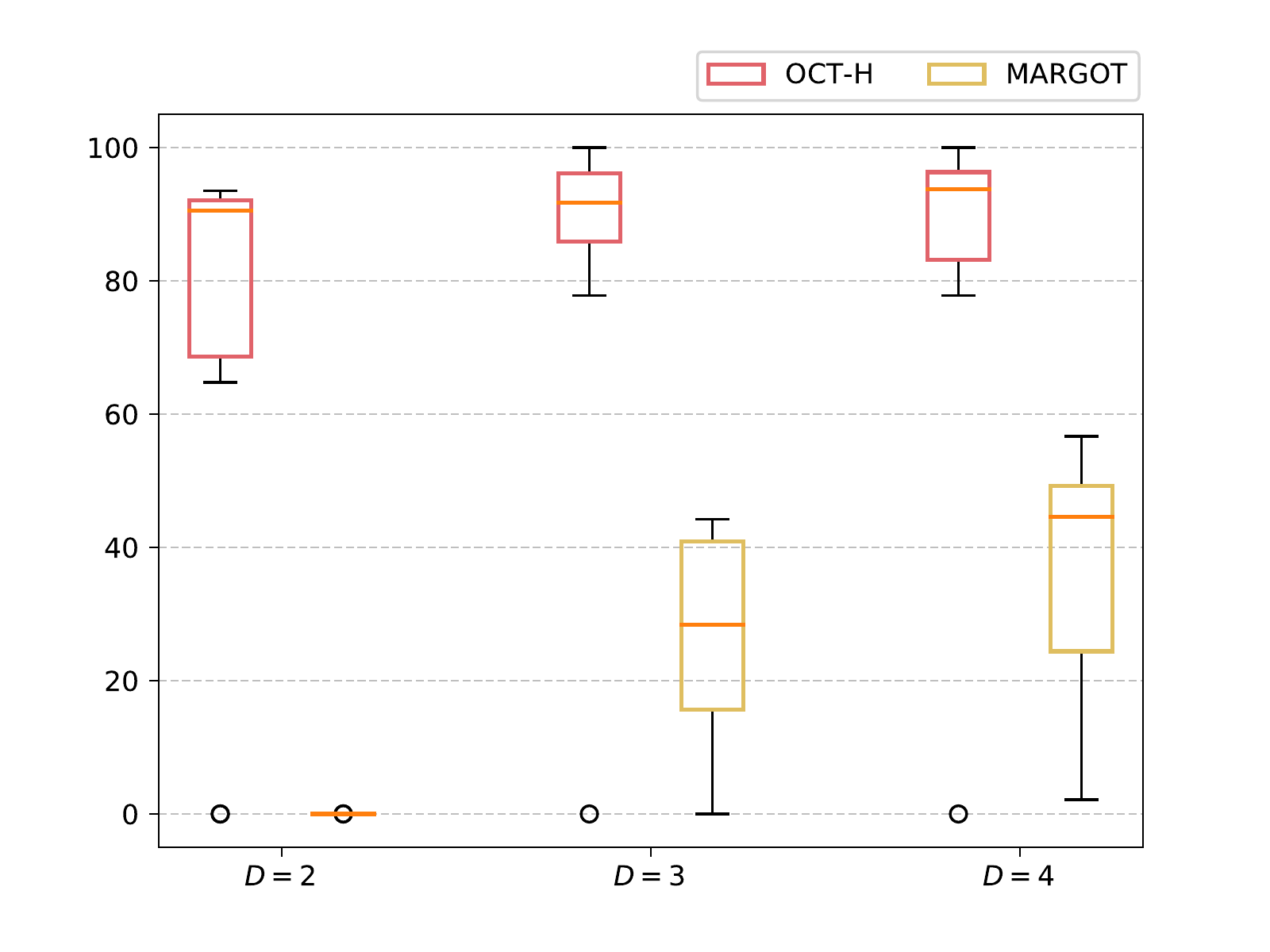}\\
    {Test ACC (\%)} & {MIP Gap (\%)}
    \end{tabular}
 \caption{Aggregated representation of the Test ACC (\%) and MIP Gap (\%) values of OCT-H and MARGOT model with depths $D\in\{2,3,4\}$.}
\label{fig: box_plot_testacc_gap}
\end{figure}


\subsection{Comparison of HFS-MARGOT against univariate optimal trees with \mnotee{higher} depths}

\red{
In the following analysis, we compare HFS-MARGOT* with $D=2$ against the univariate model OCT-1 with $D\in\{2, 3,4\}$ in order to assess if employing shallow and sparse multivariate trees ultimately pays off even if univariate ones are allowed higher depths. {In general, multivariate splits are less interpretable but allow for more flexibility, yielding shallow trees with less branching
levels than univariate ones.} In this comparison, we consider only HFS-MARGOT* in that it represents the most sparse multivariate MARGOT model. Thus, if in this extreme case\mnotee{,} HFS-MARGOT* performs better, we can deduce that MARGOT and SFS-MARGOT* will perform similarly, if not much better, in terms of predictive capabilities, but at the expense of interpretability.}

\red{For OCT-1 with $D\in\{3,4\}$, we maintained the same values of $\alpha$ and $N_{min}$ chosen for OCT-1 with $D=2$. We set a time limit of $1800$ and $3600$ seconds for $D=3$ and $D=4$ respectively. For the computation of the warm start solutions, the time limit was set to $60$ and $120$ seconds for the two depth values, respectively}. 
\red{Table \ref{table: test_acc_F_gap_BACC_uni} reports the predictive test performances and the number of all the (\mnotee{non-distinct}) features used computed as $\sum_{t\in\Tb} |F_t|$.
Moreover, the mean MIP Gap value and its standard deviation are provided for each model evaluated.}
\red{We can notice how HFS-MARGOT* with a shallow depth of 2 favourably compares \fnotee{to OCT-1} in terms of prediction quality 
even \fnotee{in the case where the latter is allowed a depth of 4}, sometimes even using far fewer features.
It is also evident that OCT-1 struggles to scale efficiently \fnotee{to} higher depths, showing larger mean gaps even if provided with higher time limit values. Though not explicitly deducible in the table, \fnotee{we} 
highlight that for OCT-1 with depths $D=\{3,4\}$, the solver never certified the optimal solution except for just a single dataset (Heart Disease C).}
\mnotee{The train predictive performances are reported in Table \ref{table: train_acc_BACC_uni}}.



\begin{table}[ht!]
\centering
\scalebox{0.82}{
\renewcommand\arraystretch{1.2}
\begin{tabular}{lcccccccccccccccc}
                 &  & \multicolumn{3}{c}{OCT-1 ($D=2$)}           &  & \multicolumn{3}{c}{OCT-1 ($D=3$)}           &  & \multicolumn{3}{c}{OCT-1 ($D=4$)}           &           & \multicolumn{3}{c}{HFS-MARGOT* ($D=2$)}     \\ \hline\hline
Dataset          &  & ACC           & BACC          & $\sum|F_t|$ &  & ACC           & BACC          & $\sum|F_t|$ &  & ACC           & BACC          & $\sum|F_t|$ & \textbf{} & ACC           & BACC          & $\sum|F_t|$ \\ \hline
Breast Cancer D. &  & 91.2          & 91.6          & 2           &  & 94.7          & 94.3          & 3           &  & \textbf{95.6} & \textbf{95.5} & 4           &           & \textbf{95.6} & \textbf{95.5} & 4           \\
Breast Cancer W. &  & 92.0          & 90.9          & 3           &  & 92.7          & 92.0          & 6           &  & 92.7          & 92.5          & 11          &           & \textbf{94.2} & \textbf{94.5} & 5           \\
Climate Model    &  & 91.7          & 60.1          & 3           &  & 91.7          & 60.1          & 6           &  & 93.5          & 76.3          & 15          &           & \textbf{96.3} & \textbf{77.8} & 6           \\
Heart Disease C. &  & 71.7          & 71.0          & 1           &  & 71.7          & 71.0          & 1           &  & 71.7          & 71.0          & 1           &           & \textbf{83.3} & \textbf{82.6} & {5}         \\
Ionosphere       &  & 91.5          & 91.7          & 2           &  & \textbf{93.0} & \textbf{92.7} & {2}         &  & \textbf{93.0} & \textbf{92.7} & 2           &           & 87.3          & 83.8          & 7           \\
Parkinsons       &  & 89.7          & 83.3          & 3           &  & 82.1          & 81.4          & 7           &  & \textbf{92.3} & \textbf{91.6} & 15          &           & 89.7          & 83.3          & 3           \\
Sonar            &  & 69.0          & 68.9          & 1           &  & 61.9          & 61.8          & 1           &  & 69.0          & 68.9          & 1           &           & \textbf{71.4} & \textbf{70.9} & 5           \\
SPECTF H.        &  & 77.8          & 65.8          & 2           &  & \textbf{83.3} & \textbf{69.2} & 3           &  & 79.6          & 63.5          & 6           &           & 75.9          & 57.8          & 6           \\
Tic-Tac-Toe      &  & 69.3          & 61.5          & 2           &  & 74.5          & 63.4          & 3           &  & 72.9          & 65.0          & 5           &           & \textbf{76.0} & \textbf{65.7} & 5           \\
Wholesale        &  & \textbf{86.4} & \textbf{85.2} & 3           &  & 85.2          & 83.5          & 2           &  & 84.1          & 82.6          & 4           &           & \textbf{86.4} & \textbf{85.2} & 5           \\ \hline
Mean Gap         &  & \multicolumn{3}{c}{55.6 ($\pm$ 44.3)}       &  & \multicolumn{3}{c}{75.7 ($\pm$ 36.5)}       &  & \multicolumn{3}{c}{77.9 ($\pm$ 35.6)}       &           & \multicolumn{3}{c}{\textbf{38.0} ($\pm$ 49.2)}       \\ \hline
\end{tabular}}
\caption{\red{Results on the test predictive performances (\mnotee{test ACC (\%) and test BACC (\%)}) of OCT-1 model with $D\in\{2,3,4\}$ and HFS-MARGOT* model with $D=2$ and comparison on the number of (\mnotee{non-distinct}) features used $\sum_{t\in\Tb}|F_t|$. The last row is the Mean MIP Gap values ($\pm$ standard deviations) (\%) among all datasets for each model. }}
\label{table: test_acc_F_gap_BACC_uni}
\end{table}


\subsection{Results with a fixed set of hyperparameters}
\fnotee{Here} we provide a set of experiments \fnotee{on} MARGOT, HFS-MARGOT and SFS-MARGOT, with depth $D=2$, using a fixed set of hyperparameters across all datasets. The time limit was set to $600$ seconds, while the warm start optimization was given a limit of $30$ seconds. These experiments allow us to evaluate the robustness of the proposed methods without conducting an ad-hoc hyperparameter search for each dataset.
We have chosen a set of hyperparameters that was never selected by the 4-FCV and that could be considered a reasonable choice. 
The fixed set of hyperparameters chosen for the MARGOT models \mnotee{is} reported in Table \ref{table: hyperparameters_fixed}.

\red{Tables \ref{table: test_acc_BACC_fixed}, \ref{table: train_acc_BACC_fixed}, \ref{table: time_gap_fixed} provide the predictive and optimization performances obtained with the fixed set of hyperparameters. 
For each \fnotee{performance measure }
(ACC, BACC, Time or Gap), we provide the difference $d$ defined as:
$$ d = \pm \left|\text{best value - actual value}\right|,$$
where the $\text{best value}$ is the value of the \fnotee{performance measure}
obtained using the hyperparameters in Tables \ref{table: hyperparameters} and \ref{table: hyperparameters_FS} obtained with the 4-FCV, whereas  $\text{actual value}$ is the value 
obtained} with the hyperparameters set as in Table \ref{table: hyperparameters_fixed}. The sign is "$+$" when the best value obtained by 4-FCV is better, \mnotee{otherwise it is} "$-$".}
\red{Table \ref{table: test_acc_BACC_fixed} reports the predictive performances on the test sets. These results show that, 
despite using the same set of hyperparameters, all the methods achieve good generalization performances across all the datasets. 
Additionally, Table \ref{table: time_gap_fixed} provides the computational times and MIP Gap values. It is evident that the MARGOT model is easier to solve, certifying the optimal solutions in 8 out of 10 datasets. 
Performances on the train sets can be found in Table \ref{table: train_acc_BACC_fixed}, where it is possible to notice that MARGOT and SFS-MARGOT tend to overfit on some datasets. This result was expected in that, without conducting \mnotee{a} cross-validation procedure, the risk of overfitting is higher.}
\red{Overall, these results show that MARGOT models are robust, reporting overall good prediction quality, without the need for a dataset-specific hyperparameter search. However, a tailored hyperparameter tuning can be advisable to mitigate the risk of overfitting and to maximize the model's potential.}

\begin{table}[ht!]
\centering
\scalebox{0.95}{
\renewcommand\arraystretch{1.2} 
\begin{tabular}{lcccccccccccc}
        &  & \multicolumn{2}{c}{MARGOT} &  & \multicolumn{4}{c}{HFS-MARGOT}           &  & \multicolumn{3}{c}{SFS-MARGOT}  \\ \hline\hline
Dataset &  & $C_0$     & $C_1 = C_2$    &  & $C_0$  & $C_1 = C_2$ & $B_0$ & $B_1=B_2$ &  & $C_0$  & $C_1 = C_2$ & $\alpha$ \\ \hline
all     &  & $10^2$    & $10^3$         &  & $10^2$ & $10^3$      & 2     & 2         &  & $10^2$ & $10^3$      & $2^{5}$ \\ \hline
\end{tabular}}
\caption{\red{Hyperparameters selected for results in Table \ref{table: test_acc_BACC_fixed}, Table \ref{table: train_acc_BACC_fixed} and Table \ref{table: time_gap_fixed}.}}
\label{table: hyperparameters_fixed}
\end{table}

\begin{table}[ht!]
\centering
\scalebox{0.89}{
\renewcommand\arraystretch{1.2} 
\begin{tabular}{lcccrcccrcccr}
                 &  & \multicolumn{3}{c}{MARGOT}         &  & \multicolumn{3}{c}{HFS-MARGOT}     &  & \multicolumn{3}{c}{SFS-MARGOT}     \\ \hline\hline
Dataset          &  & ACC  & BACC & $d$ &  & ACC  & BACC & $d$ &  & ACC  & BACC & $d$ \\ \hline
Breast Cancer D. &  & 94.7 & 94.3 & 2.6 / 2.6            &  & 95.6 & 95.5 & 0.0 / 0.0            &  & 95.6 & 94.5 & -0.9 / -0.2          \\
Breast Cancer W. &  & 93.4 & 93.0 & 2.9 / 3.2            &  & 93.4 & 93.5 & 0.7 / 1.0            &  & 93.4 & 93.0 & 0.7 / 0.6            \\
Climate Model    &  & 94.4 & 81.8 & 2.8 / 6.6            &  & 92.6 & 65.7 & 3.7 / 12.1           &  & 95.4 & 87.4 & 0.9 / -4.5           \\
Heart Disease C. &  & 81.7 & 81.7 & 1.7 / 1.3            &  & 73.3 & 72.3 & 10.0 / 10.3          &  & 83.3 & 83.3 & 3.3 / 2.9            \\
Ionosphere       &  & 87.3 & 83.8 & 5.6 / 6.2            &  & 85.9 & 81.8 & 1.4 / 2.0            &  & 83.1 & 78.7 & 1.4 / 1.1            \\
Parkinsons       &  & 89.7 & 93.1 & -5.1 / -10.0         &  & 87.2 & 81.6 & 2.6 / 1.7            &  & 92.3 & 94.8 & -5.1 / -10.0         \\
Sonar            &  & 73.8 & 73.0 & 0.0 / 0.0            &  & 69.0 & 68.9 & 2.4 / 2.0            &  & 76.2 & 75.5 & -2.4 / -1.8          \\
SPECTF H.        &  & 75.9 & 61.2 & 3.7 / -4.4           &  & 79.6 & 60.1 & -3.7 / -2.3          &  & 74.1 & 56.7 & 9.3 / 9.2            \\
Tic-Tac-Toe      &  & 95.3 & 95.4 & 4.7 / 2.4            &  & 70.3 & 62.0 & 5.7 / 3.7            &  & 96.4 & 95.8 & 1.6 / 1.2            \\
Wholesale        &  & 87.5 & 85.1 & 0.0 / 0.0            &  & 88.6 & 86.9 & -2.3 / -1.7          &  & 87.5 & 85.1 & 1.1 / 1.8            \\ \hline
Mean             &  & 87.4 & 84.2 & 1.9 / 0.8            &  & 83.6 & 76.8 & 2.1 / 2.9            &  & 87.7 & 84.5 & 1.0 / 0.0            \\ \hline
\end{tabular}}
\caption{\red{Results on the predictive test performances \mnotee{(test ACC (\%) and test BACC (\%))} of MARGOT models with $D=2$ and the hyperparameters set as in Table \ref{table: hyperparameters_fixed} and difference $d = \pm |\text{best value - actual value}|$; a value $d>0$ indicates an advantage for 4-FCV selection, while $d<0$ denotes an advantage for the fixed set of hyperparameters.}}
\label{table: test_acc_BACC_fixed}
\end{table}

\begin{table}[ht!]
\centering
\scalebox{0.89}{
\renewcommand\arraystretch{1.2} 
\begin{tabular}{lcccrcccrcccr}
                 &  & \multicolumn{3}{c}{MARGOT}  &  & \multicolumn{3}{c}{HFS-MARGOT} &  & \multicolumn{3}{c}{SFS-MARGOT} \\ \hline\hline
Dataset          &  & ACC   & BACC  & $d$         &  & ACC    & BACC   & $d$          &  & ACC    & BACC   & $d$          \\ \hline
Breast Cancer D. &  & 100.0 & 100.0 & -0.7 / -0.9 &  & 98.0   & 97.4   & 0.0 / 0.0    &  & 100.0  & 100.0  & -1.8 / -2.2  \\
Breast Cancer W. &  & 100.0 & 100.0 & -1.1 / -0.9 &  & 97.6   & 97.7   & 0.7 / 0.9    &  & 100.0  & 100.0  & -2.2 / -2.4  \\
Climate Model    &  & 100.0 & 100.0 & -0.7 / -4.1 &  & 95.1   & 74.1   & 1.6 / 7.0    &  & 100.0  & 100.0  & -3.2 / -15.2 \\
Heart Disease C. &  & 93.7  & 93.6  & -4.6 / -4.7 &  & 83.5   & 83.0   & -0.8 / -0.9  &  & 89.5   & 89.4   & -3.0 / -3.3  \\
Ionosphere       &  & 100.0 & 100.0 & -1.1 / -1.3 &  & 92.5   & 90.5   & 1.4 / 1.8    &  & 100.0  & 100.0  & -9.6 / -12.9 \\
Parkinsons       &  & 100.0 & 100.0 & 0.0 / 0.0   &  & 92.3   & 86.0   & -1.3 / -3.5  &  & 100.0  & 100.0  & 0.0 / 0.0    \\
Sonar            &  & 100.0 & 100.0 & 0.0 / 0.0 &  & 85.5   & 85.3   & -6.6 / -7.1  &  & 100.0  & 100.0  & -9.0 / -9.7  \\
SPECTF H.        &  & 100.0 & 100.0 & -3.3 / -8.0 &  & 86.9   & 75.7   & -1.9 / -1.2  &  & 100.0  & 100.0  & -8.9 / -19.1 \\
Tic-Tac-Toe      &  & 100.0 & 100.0 & -1.6 / -2.3 &  & 72.6   & 64.6   & 5.5 / 3.7    &  & 100.0  & 100.0  & -1.6  / -2.3 \\
Wholesale        &  & 96.3  & 95.2  & 0.0 / 0.0   &  & 93.8   & 93.3   & 0.0 / 1.4    &  & 96.3   & 95.2   & -1.1 / -0.6  \\ \hline
Mean             &  & 99.0  & 98.9  & -1.3 / -2.2 &  & 89.8   & 84.8   & -0.1 / 0.2   &  & 98.6   & 98.5   & -4.0 / -6.8  \\ \hline
\end{tabular}}
\caption{
\red{Results on the predictive train performances \fnotee{(train ACC (\%) and train BACC (\%))} of MARGOT models with $D=2$ and the hyperparameters set as in Table \ref{table: hyperparameters_fixed} and difference $d = \pm |\text{best value - actual value}|$; a value $d>0$ indicates an advantage for 4-FCV selection, while $d<0$ denotes an advantage for the fixed set of hyperparameters.}
}
\label{table: train_acc_BACC_fixed}
\end{table}

\begin{table}
\centering
\scalebox{0.89}{
\renewcommand\arraystretch{1.2} 
\begin{tabular}{lcccrcccrcccr}
                 &  & \multicolumn{3}{c}{MARGOT}   &  & \multicolumn{3}{c}{HFS-MARGOT} &  & \multicolumn{3}{c}{SFS-MARGOT} \\ \hline\hline
Dataset          &  & Time  & Gap  & $d$           &  & Time   & Gap   & $d$           &  & Time   & Gap   & $d$           \\ \hline
Breast Cancer D. &  & 3.9   & 0.0  & -3.5 / 0.0    &  & \underline{600.8}  & 79.6  & 286.9 / 79.6  &  & \underline{608.0}  & 11.1  & 293.4 / 11.1  \\
Breast Cancer W. &  & 149.2 & 0.0  & 140.5 / 0.0   &  & 126.1  & 0.0   & 112.1 / 0.0   &  & 451.9  & 0.0   & 414.3 / 0.0   \\
Climate Model    &  & 10.5  & 0.0  & -0.2 / 0.0    &  & \underline{600.2}  & 78.6  & -0.1 / -21.4  &  & \underline{600.3}  & 3.7   & 0.2 / -10.5   \\
Heart Disease C. &  & \underline{600.0} & 81.2 & 281.9 / 81.2  &  & \underline{600.2}  & 80.3  & 493.7 / 80.3  &  & \underline{600.2}  & 83.5  & 0.2 / 55.2    \\
Ionosphere       &  & 13.3  & 0.0  & 2.1 / 0.0     &  & \underline{600.8}  & 85.1  & 476.8 / 85.1  &  & \underline{601.8}  & 24.1  & 585.8 / 24.1  \\
Parkinsons       &  & 3.0   & 0.0  & -204.4 / 0.0  &  & \underline{600.4}  & 74.1  & 592.2 / 74.1  &  & \underline{600.6}  & 4.4   & 0.6 / -57.6   \\
Sonar            &  & 0.2   & 0.0  & -1.4 / 0.0    &  & \underline{602.7}  & 91.2  & 1.8  / -7.1   &  & \underline{620.1}  & 4.9   & 20.1 / -37.3  \\
SPECTF H.        &  & 10.2  & 0.0  & -590.0 / -3.3 &  & \underline{600.8}  & 88.7  & 0.8 / -5.1    &  & \underline{615.7}  & 20.8  & 15.6 / -79.1  \\
Tic-Tac-Toe      &  & 172.8 & 0.0  & 169.3 / 0.0   &  & \underline{602.5}  & 90.1  & 2.4 / 1.9     &  & \underline{604.3}  & 14.8  & 557.1 / 14.8  \\
Wholesale        &  & \underline{600.0} & 31.6 & 553.6 / 31.6  &  & \underline{600.1}  & 22.1  & 563.8 / 22.1  &  & \underline{600.1}  & 31.6  & 0.1 / -48.5   \\ \hline
Mean             &  & 156.3 & 11.3 & 34.8 / 11.0   &  & 553.5  & 69.0  & 253.1 / 30.9  &  & 590.3  & 19.9  & 188.7 / -12.8 \\ \hline
\end{tabular}}
\caption{
\red{Results on the optimization performances \mnotee{(computational times (s) and MIP Gaps (\%))} of MARGOT models with $D=2$ and the hyperparameters set as in Table \ref{table: hyperparameters_fixed} and difference $d = \pm |\text{best value - actual value}|$; a value $d>0$ indicates an advantage for 4-FCV selection, while $d<0$ denotes an advantage for the fixed set of hyperparameters.}
}
\label{table: time_gap_fixed}
\end{table}
\clearpage

\section{Additional tables}\label{app: appendix2}

In this section, we present additional tables. Table \ref{tab: notation} presents a summary \mnotee{of} all the notation of sets, parameters\mnotee{,} and hyperparameters adopted in the paper. ACC and BACC performances on training samples are reported in Tables \ref{table: train_acc_BACC}, \ref{table: train_acc_BACC_FS} \red{and \ref{table: train_acc_BACC_uni}}. 
Finally, for the sake of replicability, we present in Tables \ref{table: hyperparameters} and \ref{table: hyperparameters_FS} all the hyperparameters that were chosen to carry out our computational experiments. 

\begin{table}[ht]
\centering
\scalebox{0.98}{
\renewcommand\arraystretch{1.5}
\begin{tabular}{lll}
Notation                                    & \multicolumn{1}{c}{} & Description                                                                      \\ \hline\hline
\textbf{Sets}                               &                      &                                                                                  \\

$\Tb$                                       &                      & Branch nodes                                                                     \\
$\Tl$                                       &                      & Leaf nodes                                                                       \\

$\Tbf $                &                      & Branch nodes excluded the ones in the last branching level                                  \\

$\Tbl$                                      &                      & Branch nodes of the last branching level                                                \\

$\S(t)$                     &                      & Branch nodes of the subtree rooted at node $t\in\Tb$                           \\
$\S^{\prime\prime}(t) $     &                      & Nodes of $\S(t)$ in the last branching level $\Tbl$                              \\
$\S^{\prime\prime}_L(t)$    &                      & Nodes in $\Tbl$ under the left branch of $t \in \Tbf$                                     \\
$\S^{\prime\prime}_R(t)$  &                      & Nodes in $\Tbl$ under the right branch of $t \in \Tbf$                                    \\
$\I$                                        &                      & Index set of data samples                                                        \\
$\I_t$                                      &                      & Index set of data samples assigned to node $t\in\Tb$                                   \\
$\I_{L(t)}$                                      &                      & Index set of data samples assigned to the left child node of $t\in\Tb$                                                      \\
$\I_{R(t)}$                                      &                      & Index set of data samples assigned to the right child node of $t\in\Tb$                                                     \\
\textbf{Parameters}  &                      &                                                                                  \\
$n$                                         &                      & Number of features                                                               \\
$\varepsilon$                                         &                      & Parameter to model the strict inequality in routing constraints                                                               \\
$\{M_{w}, M_{\xi}, M_{\H}\}$                                         &                      & Set of Big-M parameters used in MARGOT formulations                                                             \\
\textbf{Hyperparameters}                    &                      &                                                                                  \\
$D$                                         &                      & Maximum depth  of the tree                                                                  \\
$C_t$                                           &                      & Penalty parameter on the misclassification error at node $ t\in\Tb $\\
$B_t$                       &                      & Budget value on the number of features at node $ t\in\Tb $         \\
$\alpha$                     &                      & Penalty parameter for the soft feature selection       \\
   \hline         
\end{tabular}}

\caption{Notation: sets, parameters and hyperparameters of MARGOT models.}
\label{tab: notation}
\end{table}

\begin{table}[ht]

\centering
\scalebox{1}{
\renewcommand\arraystretch{1.2}
\begin{tabular}{=l+c+c+c+c+c+c+c+c+c}
                                       &  & \multicolumn{2}{c}{OCT-H}       &  & \multicolumn{2}{c}{MM-SVM-OCT} &  & \multicolumn{2}{c}{MARGOT}    \\ \hline\hline
Dataset          &  & ACC           & BACC          &  & ACC            & BACC           &  & ACC            & BACC           \\ \hline
Breast Cancer D. &  & 98.2          & 97.8          &  & 98.7           & 98.2           &  & \textbf{99.3}  & \textbf{99.1}  \\
\rowstyle{\color{red}}Breast Cancer W. &  & 98.4           & 98.6           &  & 98.2           & 98.0          &  & \textbf{98.9} & \textbf{99.2} \\
Climate Model                          &  & 98.6           & 91.9           &  & \textbf{99.3}  & \textbf{95.9} &  & \textbf{99.3} & \textbf{95.9} \\
Heart Disease C. &  & 85.7          & 85.6          &  & 88.6           & 88.5           &  & \textbf{89.0}  & \textbf{88.9}  \\
\rowstyle{\color{red}}Ionosphere       &  & \textbf{100.0} & \textbf{100.0} &  & 98.6           & 98.2          &  & 98.9          & 98.7          \\
Parkinsons       &  & 99.4          & 98.7          &  & 96.8           & 93.4           &  & \textbf{100.0} & \textbf{100.0} \\
Sonar            &  & 98.8          & 98.8          &  & \textbf{100.0} & \textbf{100.0} &  & \textbf{100.0}           & \textbf{100.0}           \\
SPECTF H.        &  & \textbf{98.6} & \textbf{96.6} &  & 79.3           & 50.0           &  & 96.7           & 92.0           \\
\rowstyle{\color{red}}Tic-Tac-Toe      &  & \textbf{99.6}  & \textbf{99.4}  &  & 98.6           & 98.4          &  & 98.4          & 97.7          \\
Wholesale        &  & 96.0          & \textbf{95.5} &  & 83.5           & 75.0           &  & \textbf{96.3}  & 95.2           \\ \hline
\end{tabular}}
\caption{Results on the \red{train} predictive performances of the OCT models evaluated: train ACC (\%) and train BACC (\%).}
\label{table: train_acc_BACC}
\end{table}

\begin{table}

\centering
\scalebox{1}{
\renewcommand\arraystretch{1.2}
\begin{tabular}{=l+c+c+c+c+c+c+c+c+c+c+c+c}
 &  & \multicolumn{2}{c}{\red{OCT-1}} &  & \multicolumn{2}{c}{OCT-H*} &  & \multicolumn{2}{c}{HFS-MARGOT*} &  & \multicolumn{2}{c}{SFS-MARGOT*} \\ \hline\hline
Dataset &  & ACC & BACC &  & ACC & BACC &  & ACC & BACC &  & ACC & BACC \\ \hline
Breast Cancer D. &  & 94.3 & 93.4 &  & \textbf{98.2} & \textbf{97.8} &  & 98.0 & 97.4 &  & \textbf{98.2} & \textbf{97.8} \\
\rowstyle{\color{red}} Breast Cancer W. &  & 96.3 & 95.7 &  & 97.3 & 96.9 &  & \textbf{98.4} & \textbf{98.6} &  & 97.8 & 97.6 \\
Climate Model &  & 93.8 & 70.9 &  & \textbf{98.1} & \textbf{91.6} &  & \red{96.8} & \red{81.1} &  & 96.8 & 84.8 \\
Heart Disease C. &  & 77.6 & 77.5 &  & 85.7 & 85.4 &  & 82.7 & 82.1 &  & \textbf{86.5} & \textbf{86.1} \\
\rowstyle{\color{red}} Ionosphere &  &90.7 & 90.1 &  & {90.7} & 87.6 &  & \textbf{93.9} & \textbf{92.2} &  & 90.4 & 87.1 \\
Parkinsons &  & 91.0 & 83.4 &  & 96.2 & 94.8 &  & \red{91.0} & \red{82.5} &  & \textbf{100.0} & \textbf{100.0} \\
Sonar &  & 77.7 & 77.2 &  & \textbf{97.0} & \textbf{96.8} &  & 78.9 & 78.2 &  & 91.0 & 90.3 \\
SPECTF H. &  & 83.1 & 70.0 &  & 88.3 & 73.3 &  & 85.0 & 74.6 &  & \textbf{91.1} & \textbf{80.9} \\
\rowstyle{\color{red}} Tic-Tac-Toe &  & 70.9 & 61.9 &  & \textbf{99.1} & \textbf{98.7} &  & 78.1 & 68.3 &  & 98.4 & 97.7 \\
Wholesale &  & \textbf{95.5} & \textbf{95.5} &  & 91.5 & 89.1 &  & 93.8 & {94.7} &  & {95.2} & 94.6 \\ \hline
\end{tabular}}
\caption{Results on the \red{train} predictive performances of the \red{OCT models with feature selection}: train ACC (\%) and train BACC (\%).}
\label{table: train_acc_BACC_FS}
\end{table}

\begin{table}[ht]
\centering
\scalebox{1}{
\renewcommand\arraystretch{1.2}
\begin{tabular}{lcccccccccccc}
                 &  & \multicolumn{2}{c}{OCT-1 ($D=2$)} &  & \multicolumn{2}{c}{OCT-1 ($D=3$)} &  & \multicolumn{2}{c}{OCT-1 ($D=4$)} &  & \multicolumn{2}{c}{HFS-MARGOT* ($D=2$)} \\ \hline\hline
Dataset          &  & {ACC}    & {BACC}   &  & {ACC}    & {BACC}   &  & {ACC}    & {BACC}   &  & {ACC}       & {BACC}     \\ \hline
Breast Cancer D. &  & 94.3            & 93.4            &  & 95.6            & 94.5            &  & 96.9            & 96.2            &  & \textbf{98.0}      & \textbf{97.4}     \\
Breast Cancer W. &  & 96.3            & 95.7            &  & 97.1            & 96.8            &  & 97.3            & 97.0            &  & \textbf{98.4}      & \textbf{98.6}     \\
Climate Model    &  & 93.8            & 70.9            &  & 93.8            & 70.9            &  & 94.2            & 74.8            &  & \textbf{96.8}      & \textbf{81.1}     \\
Heart Disease C. &  & 77.6            & 77.5            &  & 77.6            & 77.5            &  & 77.6            & 77.5            &  & \textbf{82.7}      & \textbf{82.1}     \\
Ionosphere       &  & 90.7            & 90.1            &  & 90.7            & 90.1            &  & 90.7            & 90.1            &  & \textbf{93.9}      & \textbf{92.2}     \\
Parkinsons       &  & 91.0            & 83.4            &  & 96.8            & 94.3            &  & \textbf{98.1}   & \textbf{97.8}   &  & 91.0               & 82.5              \\
Sonar            &  & 77.7            & 77.2            &  & 77.7            & 77.3            &  & 77.7            & 77.2            &  & \textbf{78.9}      & \textbf{78.2}     \\
SPECTF H.        &  & 83.1            & 70.0            &  & \textbf{86.4}   & \textbf{74.6}   &  & 85.4            & 72.3            &  & 85.0               & \textbf{74.6}              \\
Tic-Tac-Toe      &  & 70.9            & 61.9            &  & 75.1            & 64.0            &  & 75.3            & 66.7            &  & \textbf{78.1}      & \textbf{68.3}     \\
Wholesale        &  & \textbf{95.5}   & \textbf{95.5}   &  & 94.6            & 94.9            &  & 94.9            & \textbf{95.5}            &  & 93.8               & 94.7              \\ \hline
\end{tabular}}
\caption{\red{Results on the train predictive performances of OCT-1 model with $D\in\{2,3,4\}$ and HFS-MARGOT model with $D=2$: train ACC (\%) and train BACC (\%).}}
\label{table: train_acc_BACC_uni}
\end{table}

\begin{table}
\centering
\scalebox{1.0}{
\renewcommand\arraystretch{1.2} 
\begin{tabular}{=l+c+c+c+c+c+c+c+c}
                                   &  & OCT-H    &  & \multicolumn{2}{c}{MM-SVM-OCT} &  & \multicolumn{2}{c}{MARGOT} \\ \hline\hline
Dataset                            &  & $\alpha$ &  & $c_1$          & $c_3$         &  & $C_0$       & $C_1 = C_2$  \\ \hline
Breast Cancer D.                   &  & $2^{-5}$ &  & $10^{4}$       & $10^0$        &  & $10^0$      & $10^0$       \\
\rowstyle{\color{red}} Breast Cancer W. &  & \textbf{$2^{-7}$} &  & \textbf{$10^{4}$} & \textbf{$10^{1}$}  &  & \textbf{$10^{2}$} & \textbf{$10^{2}$} \\
Climate Model                      &  & $2^{-6}$ &  & $10^{2}$       & $10^{-2}$     &  & $10^0$      & $10^0$       \\
Heart Disease C.                   &  & $2^{-5}$ &  & $10^1$         & $10^{-2}$     &  & $10^{-1}$   & $10^{-1}$    \\
\rowstyle{\color{red}} Ionosphere       &  & \textbf{$0$}      &  & \textbf{$10^2$}   & \textbf{$10^{-2}$} &  & \textbf{$10^1$}   & \textbf{$10^1$}   \\
Parkinsons                         &  & $2^{-8}$ &  & $10^3$         & $10^{1}$      &  & $10^0$      & $10^{4}$     \\
Sonar                              &  & $2^{-7}$ &  & $10^0$         & $10^0$        &  & $10^{-3}$   & $10^{-1}$    \\
SPECTF H.                          &  & $2^{-6}$ &  & $10^{-5}$      & $10^{-2}$     &  & $10^{-1}$   & $10^{-1}$    \\
\rowstyle{\color{red}} Tic-Tac-Toe &  & 0        &  & $10^5$         & $10^{1}$      &  & $10^0$      & $10^0$       \\
Wholesale                          &  & $2^{-7}$ &  & $10^3$         & $10^{-1}$     &  & $10^{3}$    & $10^{3}$     \\ \hline
\end{tabular}}
\caption{Hyperparameters selected for results in Table \ref{table: test_acc_BACC}, Table \ref{table: time_gap}, \mnotee{Table \ref{table: features}} and Table \ref{table: train_acc_BACC}.}
\label{table: hyperparameters}
\end{table}

\begin{table}
\centering
\scalebox{1.0}{
\renewcommand\arraystretch{1.2} 
\begin{tabular}{=l+c+c+c+c+c+c+c+c+c+c+c+c+c}
 &  & \red{OCT-1} &  & OCT-H* &  & \multicolumn{4}{c}{HFS-MARGOT*} &  & \multicolumn{3}{c}{SFS-MARGOT*} \\ \hline\hline
Dataset &  & $\alpha$ &  & $\alpha$ &  & $C_0$ & $C_1 = C_2$ & $B_0$ & $B_1=B_2$ &  & $C_0$ & $C_1 = C_2$ & $\alpha$ \\ \hline
Breast Cancer D. &  & $2^{-8}$ &  & $2^{-5}$ &  & $10^{3}$ & $10^{3}$ & 2 & 2 &  & $10^{2}$ & $10^{2}$ & $2^{10}$ \\
\rowstyle{\color{red}}  Breast Cancer W. &  & 0 & & $2^{-5}$ &  & $10^{5}$ & $10^{5}$ & 2 & 3 &  & $10^0$ & $10^0$ & $2^{4}$ \\
Climate Model &  & 0 &  & $2^{-4}$ &  & \red{$10^{0}$} & \red{$10^5$} & \red{3} & \red{3} &  & $10^{2}$ & $10^{2}$ & $2^{10}$ \\
Heart Disease C. &  & $2^{-2}$ &  & $2^{-4}$ &  & $10^1$ & $10^{3}$ & 1 & 2 &  & $10^0$ & $10^0$ & $2^{2}$ \\
\rowstyle{\color{red}}  Ionosphere &  & $2^{-3}$ &  & $2^{-4}$ &  & $10^{1}$ & $10^{1}$ & 2 & 3 &  & $10^0$ & $10^0$ & $2^{8}$ \\
Parkinsons &  & 0 &  & $2^{-5}$ &  & \red{$10^{3}$} & \red{$10^{3}$} & \red{1} & \red{2} &  & $10^{2}$ & $10^4$ & $2^{10}$ \\
Sonar &  & $2^{-2}$ &  & $2^{-7}$ &  & $10^{-4}$ & $10^1$ & 1 & 2 &  & $10^0$ & $10^0$ & $2^{2}$ \\
SPECTF H. &  & $2^{-7}$ &  & $2^{-5}$ &  & $10^{-4}$ & $10^{-2}$ & 2 & 3 &  & $10^{-4}$ & $10^{2}$ & $2^{8}$ \\
\rowstyle{\color{red}} Tic-Tac-Toe &  & $2^{-5}$ &  & $2^{-5}$ &  & $10^{1}$ & $10^{2}$ & 2 & 3 &  & $10^{0}$ & $10^{0}$ & $2^{0}$ \\
Wholesale &  & $2^{-8}$ &  & $2^{-6}$ &  & $10^{-2}$ & $10^{2}$ & 1 & 2 &  & {$10^{-2}$} & {$10^{0}$} & {$2^{2}$} \\ \hline
\end{tabular}}
\caption{Hyperparameters selected for results in Table \ref{table: test_acc_BACC_fs}, Table \ref{table: time_gap_fs}, \mnotee{Table \ref{table: features_fs}} and Table \ref{table: train_acc_BACC_FS}.}
\label{table: hyperparameters_FS}
\end{table}

\end{document}

%% file: SVMplot.tex
  \tikzset{
    leftNode/.style={circle, scale=0.9pt, fill=none,draw, thick},
    rightNode/.style={circle,minimum width=0.2cm, fill=black,thick,draw,scale=0.9pt},
    rightNodeInLine/.style={solid,circle, scale=1pt, fill=black,thick,draw=white},
    leftNodeInLine/.style={solid,circle,scale=1pt, fill=none,thick,draw},
  }
  \begin{tikzpicture}[
        scale=2,
        important line/.style={thick}, dashed line/.style={dashed, thin},
        every node/.style={color=black},
    ]
    \draw[dashed line, yshift=.7cm]
       (.1,.1) coordinate (sls) -- (2.7,2.7) coordinate (sle)
       node[leftNodeInLine,draw=pred,scale=1pt] (name) at (2.3,2.3){}
       node[solid,circle,draw=pred,scale=0.69pt,fill=none,thick,draw] (name) at (2.3,2.3){}
       node [above right] {$w^Tx + b = 1$};

    \draw[important line]
       (.3,.3) coordinate (lines) -- (3,3) coordinate (linee)
       node [above right] {$w^Tx + b = 0$};

    \draw[dashed line, xshift=.7cm]
       (.1,.1) coordinate (ils) -- (2.7,2.7) coordinate (ile)
       node[solid,circle,draw=pblue,scale=1pt,fill=none,thick,draw] (name) at (1.8,1.8){}
       node[rightNodeInLine, draw=pblue, fill=white, scale=0.69pt] (name) at (1.8,1.8){}
       node [above right] {$w^Tx + b = -1$};

    \draw[very thick,<->, fill = white] ($(sls)+(.2,.2)$) -- ($(ils)+(.2,.2)$) node[sloped, midway, fill=white,text opacity=1,fill opacity=0.7, above] {$\bm{ \frac {2}{||w||_2} } $} 
    ;
       

    \foreach \Point in {(.5,2.4), (2.5,1), (1.1,3.1), (1.1,2.6), (1.1,1.5)}{
      \draw \Point node[leftNode, draw=pred, fill=white, scale=0.9pt]{};
    } 
      \draw (2.5,1) node[leftNode, draw=pred, fill=pred, scale=0.9pt]{};    
      
      \foreach \Point in {(2.9,1.4), (1.8,.5), (3.5,1.5), (2.3,0.6), (1.6,3)}{
      \draw \Point node[rightNode, draw=pblue, fill=white, scale=0.9pt]{};
        \draw (1.6,3)node[leftNode, draw=pblue, fill=pblue, scale=0.9pt]{};  
    }
    
    \draw[thick, dashed] (1.6,3) -- (2.67,1.98) node[anchor=north west, pos=-0.003cm, right=0.3cm] {$\xi_1$};
    \draw[thick, dashed] (1.1,1.5) -- (0.942,1.651) node[anchor=north west, pos=0.02cm, right=0.3cm] {$\xi_2$};
    \draw[thick, dashed] (2.5,1) -- (1.4,2.04) node[anchor=north west, pos=0.003cm, right=0.3cm] {$\xi_3$};
    
    \draw[thick,->] (0,0) -- (4,0) node[anchor=north west] {$x_1$};
    \draw[thick,->] (0,0) -- (0,4) node[anchor=south east] {$x_2$};

  \end{tikzpicture}

%% file: margottona_simple_def.tex
\tikzset{
  branchnode/.style = {very thick, circle, draw=dblue, text width=1.5em, text centered, minimum height=1cm, minimum width = 1cm, font=\large},
  1dummy/.style       = {font=\Large},
  rect/.style = {rectangle, align = left, draw = white, font=\huge, text width=4cm}
}

\tikzset{
    ncbar angle/.initial=90,
    ncbar/.style={
        to path=(\tikztostart)
        -- ($(\tikztostart)!#1!\pgfkeysvalueof{/tikz/ncbar angle}:(\tikztotarget)$)
        -- ($(\tikztotarget)!($(\tikztostart)!#1!\pgfkeysvalueof{/tikz/ncbar angle}:(\tikztotarget)$)!\pgfkeysvalueof{/tikz/ncbar angle}:(\tikztostart)$)
        -- (\tikztotarget)
    },
    ncbar/.default=0.3cm,
}

\tikzset{square left brace/.style={ncbar=0.3cm}}
\tikzset{square right brace/.style={ncbar=-0.3cm}}

\begin{tikzpicture}
  [
    grow                    = down,
    level 1/.style          = {sibling distance=10.2cm},
    level 2/.style          = {sibling distance=5.5cm},
    level 3/.style          = {sibling distance=2.9cm},
    level distance          = 3.3cm,
    edge from parent/.style = {thick, draw, edge from parent path={(\tikzparentnode) -- (\tikzchildnode)}},
  ]
\node [branchnode] (0) {0}           
    child {
    node[branchnode] (1) {1}
    child { node [branchnode] (3) {3} 
    child { node [branchnode, draw = mygreen] (7) {7}
         edge from parent node[1dummy, pos=0.65, left=0.4cm] (e) {$<$} }
    child { node [branchnode, draw = mygreen] (8) {8}
      edge from parent node[1dummy, pos=0.65, right=0.4cm] (k) {$\geq$}}
    edge from parent node[1dummy, pos=0.4, left=0.4cm] (c) {${\dblue{w}}^T_1 x^i + \dblue{b}_1 < 0$}}
    child { node [branchnode] (4) {4} 
    child { node [branchnode, draw = mygreen] (9) {9}
    edge from parent node[1dummy, pos=0.65, left=0.4cm] (h) {$<$} }
    child { node [branchnode, draw = mygreen] (10) {10}
      edge from parent node[1dummy, pos=0.65, right=0.4cm] (k) {$\geq$}}
    edge from parent node[1dummy, pos=0.4, right=0.4cm] {$\geq$} }
  edge from parent node [1dummy, pos=0.3, left=0.4cm] (a) {${\dblue{w}}^T_0 x^i + \dblue{b}_0 < 0\ \ $}
}
    child {
     node[branchnode] (2) {2}
    child { node  [branchnode] (5) {5} 
    child { node [branchnode, draw = mygreen] (11) {11}
        edge from parent node[1dummy, pos=0.65, left=0.4cm] (j) {$<$} }
    child { node [branchnode, draw = mygreen] (12) {12}
        edge from parent node[1dummy, pos=0.65, right=0.4cm] (k) {$\geq$} }
     edge from parent 
     node [1dummy, pos=0.4, left=0.4cm] {$<$}}
    child { node [branchnode] (6) {6} 
    child { node [branchnode, draw = mygreen] (13) {13}
      edge from parent node[1dummy, pos=0.65, left=0.4cm] (l) {$<$}}
    child { node [branchnode, draw = mygreen] (14) {14}
        edge from parent node [1dummy, pos=0.65, right=0.4cm] (f) {$\geq$}}
    edge from parent 
     node[1dummy, pos=0.4, right=0.4cm, draw=none] (d) {$\dblue{w}^T_2 x^i + \dblue{b}_2 \geq 0$}}
    edge from parent node [1dummy, pos=0.3, right=0.4cm] (b) {$\ \ {\dblue{w}}^T_0 x^i + \dblue{b}_0 \geq 0$}
};

\draw [loosely dashed] ($(e.west) + (-0.5cm,0.75cm)$) to  ($(f.east) + (+0.5cm,0.75cm)$);




\node[rect,  above right = 1 cm and 6.8 cm of 2] (Tbf) {\dblue{$\Tbf$}}; 
\node[rect, below = 4 cm of Tbf] (Tbl) {\dblue{$\Tbl$}}; 
\node[rect, below = 2.3 cm of Tbl] {\mygreen{$\Tl$}}; 


\draw [black, thick] (11,1) to [square left brace] (11,-4);
\draw [black, thick] (11,-5.5) to [square left brace] (11,-7.5);
\draw [black, thick] (11,-8.7) to [square left brace] (11,-10.7);


\node[1dummy, below = 0.03cm of 7]  {\large class $-1$};
\node[1dummy, below = 0.03cm of 8]  {\large class $1$};
\node[1dummy, below = 0.03cm of 9]  {\large class $-1$};
\node[1dummy, below = 0.03cm of 10]  {\large class $1$};
\node[1dummy, below = 0.03cm of 11]  {\large class $-1$};
\node[1dummy, below = 0.03cm of 12]  {\large class $1$};
\node[1dummy, below = 0.03cm of 13]  {\large class $-1$};
\node[1dummy, below = 0.03cm of 14]  {\large class $1$};

\end{tikzpicture}

%% file: parkinsons_tree_oct_cut.tex
\begin{tikzpicture}[->,>=stealth',auto,node distance=1.5cm,
  thick,
  node blue/.style={align=center,minimum size=2.2cm,draw,circle, draw=dblue,font=\large},
  node green/.style={align=center,minimum size=1.8cm,draw=mygreen,circle,font=\large},
  text node/.style = {font=\large},
  every edge/.style={draw=black,thick}] 
\node[node blue] (0) at (0,0) {$[38,118]$\\ $F_0$: $16$};  
    \node[node blue] (1) at (-2.5,-2.6) {$[16,21]$
    \\ $F_1$: $21$};
    \node[node blue] (2) at (2.5,-2.6) {$[22, 97]$
    \\ 
    $F_2$: $21$};  
    \node[node green] (3) at (-4,-5.2) {$[8, 1]$};  
    \node[node green] (4) at (-1,-5.2) {$[8, 20]$};  
    \node[node green] (5) at (1,-5.2) {$[18, 1]$};  
    \node[node green] (6) at (4,-5.2) {$[4, 96]$};  

    \path [-] (0) edge (1)  ;
    \path [-] (0) edge (2);  
    \path [-] (1) edge (3);  
    \path [-] (1) edge (4);  
    \path [-] (2) edge (5);  
    \path [-] (2) edge (6);  

\node[text node] (a) [below=0.2cm of 3] {class -1};
\node[text node] (b) [below=0.2cm of 4] {class 1};
\node[text node] (c) [below=0.2cm of 5] {class -1};
\node[text node] (d) [below=0.2cm of 6] {class 1};

\end{tikzpicture}  

%% file: parkinsons_tree_octh_cut.tex
\begin{tikzpicture}[->,>=stealth',auto,node distance=1.5cm,
  thick,
  node blue/.style={align=center,minimum size=2.2cm,draw,circle, draw=dblue, font=\large},
  node green/.style={align=center,minimum size=1.8cm,draw=mygreen,circle, font=\large},
  text node/.style = {font=\large},
  every edge/.style={draw=black,thick}] 
    \footnotesize
\node[node blue] (0) at (0,0) {$[38,118]$ \\

$F_0$: $16, 18, 20$};  
    \node[node blue] (1) at (-2.5,-2.6) {$[11,115]$
    \\ $F_1$: $0, 5, 16$};
    \node[node green] (2) at (2.5,-2.6) {$[27,3]$};
    \node[node green] (3) at (-4,-5.2) {$[8,0]$};  
    \node[node green] (4) at (-1,-5.2) {$[3,115]$};  

    \path [-] (0) edge (1)  ;
    \path [-] (0) edge (2);  
    \path [-] (1) edge (3);  
    \path [-] (1) edge (4);  

\node[text node] (a) [below=0.2cm of 3] {class -1};
\node[text node] (b) [below=0.2cm of 4] {class 1};
\node[text node] (d) [below=0.2cm of 2] {class -1};

\end{tikzpicture}  

%% file: parkinsons_tree_margot_hardb.tex
\begin{tikzpicture}[->,>=stealth',auto,node distance=1.5cm,
  thick,
  node blue/.style={align=center,minimum size=2.2cm,draw,circle, draw=dblue, font=\large},
  node green/.style={align=center,minimum size=1.8cm,draw=mygreen,circle, font=\large},
  text node/.style = {font=\large},
  every edge/.style={draw=black,thick}] 
\node[node blue] (0) at (0,0) {$[38,118]$\\ \red{$F_0$: $21$}};  
 \node[node blue] (1) at (-2.5,-2.6) {\red{$[26,6]$}
    \\ \red{$F_1$: $0, 16$}};
\node[node green] (2) at (2.8,-2.6) {\red{$[12,112]$}};  
    \node[node green] (3) at (-4,-5.2) {\red{$[25,1]$}};  
    \node[node green] (4) at (-1,-5.2) {\red{$[1,5]$}};  

    \path [-] (0) edge (1)  ;
    \path [-] (0) edge (2);  
    \path [-] (1) edge (3);  
    \path [-] (1) edge (4);  

\node[text node] (a) [below=0.2cm of 3] {class -1};
\node[text node] (b) [below=0.2cm of 4] {class 1};
\node[text node] (c) [below=0.2cm of 2] {class 1};

\end{tikzpicture}  

%% file: parkinsons_tree_margot_softb.tex
\begin{tikzpicture}[->,>=stealth',auto,node distance=1.5cm,
  thick,
  node blue/.style={align=center,minimum size=2.2cm,draw,circle, draw=dblue, font=\large},
  node green/.style={align=center,minimum size=1.8cm,draw=mygreen,circle, font=\large},
  text node/.style = {font=\large},
  every edge/.style={draw=black,thick}] 
\node[node blue] (0) at (0,0) {$[38,118]$\\ $F_0$: $1, 2, 5, 6,12,$\\
$ 17, 18, 19, 20$};  
    \node[node blue] (1) at (-2.5,-2.6) {$[26,30]$
    \\ $F_1$: $0, 16$};
    \node[node blue] (2) at (2.5,-2.6) {$[12,111]$
    \\ 
    $F_2$: $0, 5, 16, 21$};  
    \node[node green] (3) at (-4,-5.2) {$[26,0]$};  
    \node[node green] (4) at (-1,-5.2) {$[0,7]$};  
    \node[node green] (5) at (1,-5.2) {$[12,0]$};  
    \node[node green] (6) at (4,-5.2) {$[0,111]$};  

    \path [-] (0) edge (1)  ;
    \path [-] (0) edge (2);  
    \path [-] (1) edge (3);  
    \path [-] (1) edge (4);  
    \path [-] (2) edge (5);  
    \path [-] (2) edge (6);  

\node[text node] (a) [below=0.2cm of 3] {class -1};
\node[text node] (b) [below=0.2cm of 4] {class 1};
\node[text node] (c) [below=0.2cm of 5] {class -1};
\node[text node] (d) [below=0.2cm of 6] {class 1};

\end{tikzpicture}  

%% file: bcd_tree_oct_cut.tex
\begin{tikzpicture}[->,>=stealth',auto,node distance=1.5cm,
  thick,
  node blue/.style={align=center,minimum size=2.2cm,draw,circle, draw=dblue, font=\large},
  node green/.style={align=center,minimum size=1.8cm,draw=mygreen,circle, font=\large},
  text node/.style = {font=\large},
  every edge/.style={draw=black,thick}] 
\node[node blue] (0) at (0,0) {$[285,170]$\\ $F_0$: $26$};  
    \node[node green] (1) at (-2.5,-2.6) {$[219, 9]$};
    \node[node blue] (2) at (2.5,-2.6) {$[66, 161]$
    \\ 
   $F_2$: $7$};  
    \node[node green] (5) at (1,-5.2) {$[57, 8]$};  
    \node[node green] (6) at (4,-5.2) {$[9, 153]$};


    \path [-] (0) edge (1)  ;
    \path [-] (0) edge (2);  
    \path [-] (2) edge (5);  
    \path [-] (2) edge (6);  

\node[text node] (a) [below=0.2cm of 1] {class -$1$};
\node[text node] (c) [below=0.2cm of 5] {class -$1$};
\node[text node] (d) [below=0.2cm of 6] {class $1$};

\end{tikzpicture}  

%% file: bcd_tree_octh_cut.tex
\begin{tikzpicture}[->,>=stealth',auto,node distance=1.5cm,
  thick,
  node blue/.style={align=center,minimum size=2.2cm,draw,circle, draw=dblue, font=\large},
  node green/.style={align=center,minimum size=1.8cm,draw=mygreen,circle, font=\large},
  text node/.style = {font=\large},
  every edge/.style={draw=black,thick}] 
\node[node blue] (0) at (0,0) {$[285,170]$\\ $F_0$: $21, 22, 24$};  
    \node[node green] (1) at (-2.5,-2.5) {$[1,163]$};
    \node[node green] (2) at (2.5,-2.5) {$[284,7]$};
    


    \path [-] (0) edge (1)  ;
    \path [-] (0) edge (2);  

\node[text node] (b) [below=0.2cm of 1] {class $1$};
\node[text node] (d) [below=0.2cm of 2] {class $-1$};

\end{tikzpicture}  

%% file: bcd_tree_margot_hardb_cut.tex
\begin{tikzpicture}[->,>=stealth',auto,node distance=1.5cm,
  thick,
  node blue/.style={align=center,minimum size=2.2cm,draw,circle, draw=dblue, font=\large},
  node green/.style={align=center,minimum size=1.8cm,draw=mygreen,circle, font=\large},
  text node/.style = {font=\large},
  every edge/.style={draw=black,thick}] 
\node[node blue] (0) at (0,0) {$[285,170]$\\ $F_0$: $23, 24$};  
    \node[node green] (1) at (-2.5,-2.6) {$[278,9]$};
    \node[node blue] (2) at (2.5,-2.6) {$[7,161]$
    \\ 
    $F_2$: $1, 5$};  
    \node[node green] (5) at (1,-5.2) {$[7,0]$};  
    \node[node green] (6) at (4,-5.2) {$[0,161]$};  

    \path [-] (0) edge (1)  ;
    \path [-] (0) edge (2);  
    \path [-] (2) edge (5);  
    \path [-] (2) edge (6);  

\node[text node] (a) [below=0.2cm of 1] {class -1};
\node[text node] (c) [below=0.2cm of 5] {class -1};
\node[text node] (d) [below=0.2cm of 6] {class 1};

\end{tikzpicture}  

%% file: bcd_tree_margot_softb_cut.tex
\begin{tikzpicture}[->,>=stealth',auto,node distance=1.5cm,
  thick,
  node blue/.style={align=center,minimum size=2.2cm,draw,circle, draw=dblue, font=\large},
  node green/.style={align=center,minimum size=1.8cm,draw=mygreen,circle, font=\large},
  text node/.style = {font=\large},
  every edge/.style={draw=black,thick}] 
\node[node blue] (0) at (0,0) {$[285,170]$\\ $F_0$: $21, 22, 24$};  
    \node[node green] (1) at (-2.5,-2.6) {$[284,7]$};
    \node[node green] (2) at (2.5,-2.6) {$[1,163]$};

    \path [-] (0) edge (1)  ;
    \path [-] (0) edge (2);  

\node[text node] (a) [below=0.2cm of 1] {class -1};
\node[text node] (b) [below=0.2cm of 2] {class 1};

\end{tikzpicture}  